\documentclass[english,12pt,leqno]{article}
\usepackage{tikz}
\usetikzlibrary{decorations.markings}
\usetikzlibrary{calc}

\usepackage[T2A]{fontenc}

    \usepackage{graphicx}
    \usepackage{wrapfig, framed, caption}

\usepackage[latin1, russian, utf8]{inputenc}

\usepackage{amsxtra}
\usepackage{amsmath}
\usepackage{amssymb}
\usepackage{amsfonts}
\usepackage[all]{xy}
\usepackage{mathrsfs}
\usepackage{amsthm}
\usepackage{float}
 \usepackage{caption}
 \usepackage{epigraph}
 \usepackage{wrapfig}

\newtheorem{thm}[equation]{Theorem}

\newtheorem{cor}[equation]{Corollary}
\newtheorem{lem}[equation]{Lemma}
\newtheorem{prop}[equation]{Proposition}

\newtheoremstyle{example}{\topsep}{\topsep}%
     {}
     {}
     {\bfseries}
     {.}
     {2pt}
     {\thmname{#1}\thmnumber{ #2}\thmnote{ #3}}

   \theoremstyle{example}
   
   \newtheorem{Defi}[equation]{Definition}
   
   \newtheorem{rem}[equation]{Remark}
   \newtheorem{rems}[equation]{Remarks}

   \newtheorem{ex}[equation]{Example}

\newtheorem{conj}[equation]{Conjecture}
 
 \newtheorem{prpsl}[equation]{Proposal}

\newtheoremstyle{example}{\topsep}{\topsep}%
     {}
     {}
     {\bfseries}
     {.}
     {2pt}
     {\thmname{#1}\thmnumber{ #2}\thmnote{ #3}}

  \numberwithin{equation}{section}

\def\la{{\langle}}
\def\ra{{\rangle}}

\def\AAA{\mathbb{A}}

\def\CC{\mathbb{C}}

\def\GG{\mathbb{G}}
\def\PP{\mathbb{P}}
\def\RR{\mathbb{R}}
\def\ZZ{\mathbb{Z}}
\def\QQ{\mathbb{Q}}
\def\HH{\mathbb{H}}

\def\aen{{\mathfrak{a}}}
\def\ben{\mathfrak{b}}

\def\gen{\mathfrak{g}}

\def\hen{\mathfrak{h}}
\def\len{\mathfrak{l}}
\def\men{\mathfrak{m}}
\def\nen{\mathfrak{n}}
\def\pen{\mathfrak{p}}

\def\sen{\mathfrak{s}}

\def\Ac{\mathcal{A}}
\def\Bc{\mathcal{B}}
\def\Cc{\mathcal{C}}
\def\Kc{\mathcal{K}}
\def\Dc{\mathcal{D}}
\def\Ec{\mathcal{E}}
\def\Fc{\mathcal{F}}
\def\Gc{\mathcal{G}}
\def\Ic{\mathcal{I}}

\def\Lc{\mathcal{L}}
\def\Mc{\mathcal{M}}
\def\Nc{\mathcal{N}}
\def\Hc{\mathcal{H}}
\def\Oc{\mathcal{O}}

\def\Rc{\mathcal{R}}
\def\Sc{\mathcal{S}}
\def\Tc{\mathcal{T}}

\def\Fen{{\mathfrak {F}}}
\def\Len{{\mathfrak {L}}}

\def\bb{\mathbf{b}}

\def\hb{{\mathbf{h}}}

\def\<{\langle}
\def\>{\rangle}

\def\be{\begin{equation}}
\def\ee{\end{equation}}

\def\bef{\begin{figure}[H]\centering}
\def\enf{\end{figure}}
\def\Bl{{\on{Bl}}}
\def\Br{\on{Br}}

\def\btp{\begin{tikzpicture}}
\def\etp{\end{tikzpicture}}

\def\coh{{\on{coh}}}
\def\Coh{{\on{Coh}}}

\def\Coker{\on{Coker}}

\def\dgCat{{\tt dgCat}}
\def\dgVect{{\tt dgVect}}

\def\ds{{\Delta_{\on{sim}}}}

\def\End{{\on{End}}}
\def\Ext{\on{Ext}}

\def\Fl{{\on{Fl}}}

\def\FM{{\on{FM}}}

\def\gene{{\on{gen}}}
\def\Gr{{\on{Gr}}}

 \def\hocolim{{  \underrightarrow {\on{holim}} }}
 \def\holim{{  \underleftarrow {\on{holim}} }}
\def\Hom{\on{Hom}}

\def\Imm{\on{Im}}
\def\Id{{\on{Id}}}
\def\Irr{{\on{Irr}}}

\def\k {\mathbf k}
\def\Ker{\on{Ker}}

\def\lla{\longleftarrow}
\def\lra{\longrightarrow}

\def\Mat{{\on{Mat}}}

\def\min{{\on{min}}}

\def\Mrt{{\tt Mrt}}

\def\lra{\longrightarrow}

\def\nr{{\on{nr}}}

\def\Ob{{\on{Ob}}}
\def\ol{\overline}
\def\on{\operatorname}
\def\op{{\on{op}}}

\def\Perf{\on{Perf}}
\def\Perv{\on{Perv}}
\def\Pic{\on{Pic}}
\def\pr{{\on{pr}}}
\def\Proj{{\on{Proj}}}

\def\Qeq{{\tt Qeq}}

\def\rel{{\on{rel}}}

\def\reg{{\on{reg}}}
\def\RHom{{\on{RHom}}}
\def\rss{{\on{rss}}}

\def\Sch{{\on{Sch}}}
\def\sel{{\sen\len}}

\def\sing{{\on{sing}}}
\def\smth{{\on{sm}}}

\def\Span{{\on{Span}}}
\def\Spec{{\on{Spec}}}

\def\underbar{{\underbar}}
\def\ul{\underline}

\def\Vect{\on{Vect}}

\def\wh{\widehat}
\def\wt{\widetilde}

\title{ Perverse schobers and birational geometry}

\author{ Alexey Bondal, Mikhail Kapranov, Vadim Schechtman }

\begin{document}

 \maketitle
 
  \thanks {\em To Sasha Beilinson on his 60th birthday}

 \setlength\epigraphwidth{.5\textwidth}
 \epigraph{  \`A travers le brouillard, il contemplait
 des clochers, des \'edifices,  dont il ne savait pas les noms.
  }{Flaubert, {\em  L'\'Education sentimentale.}}

\begin{abstract}
Perverse schobers are conjectural categorical analogs of perverse sheaves.
We show that such structures  appear naturally in Homological Minimal Model Program
which studies the effect of birational transformations such as flops, on the coherent derived
categories. More precisely, the flop data are analogous to hyperbolic stalks of a perverse sheaf.

In the first part of the paper we 
study
schober-type diagrams of categories
 corresponding to flops of relative dimension 1, in particular we determine the categorical analogs of
 the (compactly supported) cohomology with coefficients in such schobers.
 
 In the second part we consider the example of a ``web of flops" provided by the Grothendieck
 resolution associated to a reductive Lie algebra $\gen$ and study the
 corresponding schober-type  diagram. For $\gen=\sel_3$ we relate this diagram to the 
 classical  {\em space of complete triangles} studied by Schubert, Semple and others.

\end{abstract}

 \addtocounter{section}{-1}
 
 \vfill\eject
 \tableofcontents

 \section{Introduction}

 \paragraph{A. Goals of the paper.} Perverse schobers are conjectural categorical analogs of perverse sheaves.
  The possibility of   a meaningful categorified theory  of perverse sheaves
  was suggested in \cite {KS-schobers}. 
 It gradually becomes clear that such a   theory   must  indeed exist  
 and have applications to various areas of mathematics.
  In some simple cases a precise definition of  perverse schobers can be given using
 quiver description of perverse sheaves as a starting point and then replacing quivers
 by  analogous diagrams of triangulated categories and exact functors.

 One important role of perverse schobers is that they can serve as natural
 coefficient data for forming Fukaya categories \cite{DKSS}, just 
  like sheaves can serve as coefficient data for forming cohomology. 
 
  \vskip .3cm
 
 The goal of this paper is to investigate a different, perhaps ``mirror dual''  to the above,
 appearance    of pervese schobers: 
  in birational geometry. We are talking especially about 
  the Homological Minimal Model Program which studies derived categories of
 coherent sheaves on algebraic varieties related by flops.  Our starting point was the observation
 that the basic framework of the derived equivalence corresponding to a flop \cite{BO}
 matches very precisely the ``hyperfunction'' description (given in \cite{KS} and recalled in
 Proposition \ref{prop:dirac} below)
 of $\Perv(\CC, 0)$,
 the category  of  perverse sheaves on  $\CC$
  smooth outside $0$. 
  We do not know  any {\em a priori} reason for this  remarkable match. 
 
 \paragraph{B. Brief summary of the paper.} 
 In the first part of the paper, we study various features of the schober-type diagrams of categories
 corresponding to flops of relative dimension 1. In this case the corresponding perverse schobers $\Fen$
 (we call it them {\em flobers})  can be seen as categorifications of objects from $\Perv(\CC,0)$. 
 We describe the categories having the meaning of $\HH^i(\CC,\Fen)$ and $\HH^i_c(\CC,\Fen)$. 
 In particular, we see an appearance of the ``categorical Poincar\'e duality'' between $\HH^0$
 and $\HH^2_c$: it corresponds to  To\"en's  Morita duality \cite{toen} between the coherent derived
 category $D^b(Z)$ and the category of perfect complexes $\Perf(Z)$ of a singular projective variety $Z$. 
 
 \vskip .2cm
 
 At the same time, the  description of $\Perv(\CC,0)$ mentioned above,   is a particular case of a much more general
 classification result \cite{KS} for
 perverse sheaves on $\CC^n$  smooth with respect  to a stratification given by an arrangement of hyperplanes $\Hc$
 with real equations. That result (recalled as Theorem \ref{thm:perv-arr} below) is formulated in terms of combinatorics
 of the {\em chamber structure} of 
  the real part $\RR^n\subset \CC^n$   (including walls between chambers and, more generally, faces
  of all dimensions)
 given by $\Hc$. Again, we note that similar kinds of chamber structures appear in birational geometry \cite{kaw} 
 in the study of the cone of movable
 divisors and the corresponding flops. \vskip .2cm
 
 In the second part of the paper, we consider the prototypical example where the movable cone with
 its chamber structure can be investigated in great detail: that  of the simultaneous Grothendieck resolution
 associated to a reductive Lie algebra $\gen$. In this case the chamber structure is given by the Weyl
 chambers in  the lattice of integer weights $\hen_\ZZ^*$. 
 
 \paragraph{C. Relation to earlier work.} 
Whatever possible approach one chooses, a perverse schober on a stratified variety $(M,\{M_\alpha\}))$ 
amounts to the following levels of structure:
\begin{itemize}
\item[(0)] A local system $\Len$  of (enhanced) triangulated categories on the open stratum $M_0$. 

\item [(1, 2,$\cdots$)] ``Singularity data'' extending $\Len$ to strata of higher codimension. 
\end{itemize}

In several situations, the local system $\Len$ is known; our observation is  that (some of) the singularity
data are also naturally present. 

\vskip .2cm

  We note the recent work of  Donovan and Wemyss \cite{DW}  who  associated,  to a system 
 of 3-dimensional flops, 
 a local system $\Len$  of triangulated categories on the complement of a hyperplane arrangement closely related
 to the root arrangements of types A-D-E. 
 
 \vskip .2cm
 
 In the context of Grothendieck resolutions, the work of 
  Bezrukavnikov \cite{bez}   and Bezrukavnikov-Riche \cite{BR} has, quite some time ago, constructed the
  braid group action on the  corresponding coherent derived category. 
 This gives a  ($W$-equivariant) local system $\Len$  of triangulated categories on the open stratum $\hen^*_\reg$
 in the space of complex weights.
 These results
  can be seen as ``quasi-classical'' analogs
 of the following fundamental fact: the braid group associated to a semisimple Lie algebra $\gen$, acts 
 by Schubert correspondences,
 in $D^b_{\on{constr}}(G/B)$, the constructible derived category of the flag variety $G/B$. 
 Our point of view is that this local system $\Len$ should admit a natural extension to
 a perverse schober on entire $\hen^*$ smooth with respect to the stratification given by the arrangement of the  coroot
 hyperplanes. Indeed, the transition functors of our flober diagram along real codimension 1 faces
 give precisely this local system. 
 
 \vskip .2cm
 
In fact,   \cite{bez}  \cite{BR} have constructed the action of the {\em affine Weyl group}, not just the
Artin braid group. From our point of view, this should be extended to a schober on the  dual  torus $H^\vee$ 
(instead of the
Cartan subalgebra) smooth with respect to the stratification given by the arrangement of the  ``coroot subtori''.
This indicates that the types of perverse schobers appearing in  birational geometry
 go beyond the realm of  hyperplane arrangements.

\vskip .2cm
 
 Further, our approach seems closely related to the general philosophy of Bezrukavnikov \cite{bez-talk} concerning derived
 categories of symplectic resolutions. That context is more general than  that of the  
 Grothendieck resolution (which is  the ``universal deformation'' of 
  $T^*(G/B)$,  a symplectic resolution of the nilpotent cone)
 and includes, for instance, all quiver varieties of Nakajima. 
  In that case, the chamber structure also plays a crucial role. 
  Moreover,  the point of view based on quantization in characteristic $p$, provides additional insight into
  the existence and nature of derived equivalences between resolutions corresponding to adjoining chambers
  (``symplectic flops''). 
  
  \vskip .2cm

 If we go one categorical level down, we get actions of affine Hecke algebras  on various  Grothendieck groups
 related to  $T^*(G/B)$ and Grothendieck resolutions which have been, for some time,  a classical subject in 
 representation theory
 \cite{chriss-ginz}.  (They are, of course, modern manifestations of the perennial
 {\em principal series intertwiners}, familiar in representation theory since the 1940's.)
  However, 
 such questions end up being quite subtle, due to the fact that various composite correspondences become 
 very singular as algebraic varieties. In our context, the analogs of such correspondences are given by various
 fiber products $X_C$ of the flopped Grothendieck resolutions.  They are labelled by cells $C$ of the root arrangements
 in the space of real weights.  
 
 \vskip.2cm
 
 We observe that for $\gen=\sel_n$ the highest fiber product $X_0$ s related to the classical {\em variety of simplices}  in $\PP^{n-1}$
 \cite{magyar} \cite{babson1} \cite{babson2}.
 In particular, for $n=3$, the variety of triangles $T$ has been studied  in the XIX century 
 by Schubert \cite{schubert} who has constructed its desingularization,  the {\em variety of complete triangles}
  $\wh T$,
 see  \cite{semple}\cite{roberts}.
  Using $\wh T$ allows us to lift the diagram of fiber products for $\sel_3$  to  a diagram containing only smooth varieties.  The relation of 
  Schubert's  complete triangles with    correspondences  (also bearing
  Schubert's name) on  the Grothendieck resolution has, 
 apparently,  not been noticed before. 
 
 \vskip .2cm
 
 In this paper we consider only ``B-model'' aspects of perverse schobers, i.e., schobers formed 
 by coherent derived categories only. Relations of such schobers to homological mirror symmetry,
 cf.  the work of Harder-Katzarkov  \cite{harder}, seem very interesting.

 \paragraph{D. Organization and results of the paper. }
   \S 1 explains an analogy between flopping diagrams
 \be\label{eq:intro-flop}
 X_- \lla X_0 \lra X_+
 \ee
 and perverse sheaves on the disk (or on the complex plane $\CC$)
 smooth outside of $0$. We denote by $\Fen$ and call {\em flober} the diagram of coherent
 derived categories associated to a flopping diagram as above. 
 
 In \S 2 we make sense and identify the categories
 $\HH^i(\CC, \Fen)$ and $\HH^i_c(\CC,\Fen)$ in terms of the singular
 variety $Z$ (the common flopping contraction of $X_\pm$).  
 This is done in Theorems \ref{thm:H^0-flober}  and \ref{thm:H^2-c}). 
 
 \S 3  explains the role of hyperplane arrangements in HMMP and presents our approach to perverse schobers
 smooth with respect to an arrangement with real equations. This approach is based in
 the description of perverse sheaves given in \cite{KS}. 
 
 In \S 4 we give a detailed treatment of the Grothendieck resolution from the point of view of flops.
This material is known but we could not find an existing  presentation suitable for our purposes. 

\S 5 is devoted to constucting the ``anti-flober'', the diagram of partial blowdowns $Z_C$ of various
flopped Grothendieck resolutions These blowdowns are labelled by faces $C$ of
the coroot arrangement. We pay special attension to the singular nature of $Z_C$ for $C$
of small codimension and their relation to similar varieties for Levi subalgebras in $\gen$. 

In \S 6 we construct fiber products $X_C$ of flopped Grothendieck resolutions (the analogs of the
variety $X_0$ in \eqref{eq:intro-flop}). These varieties, although extremely
beautiful,  have not been traditionally considered
in representation theory. The closest concept discussed in the literature, is the
{\em variety of simplices}  \cite{babson1} defined as the space of maps of the Coxeter complex
of $\gen$ to the spherical Bruhat-Tits building.  We denote the principal component of
this variety by $T_\gen$. 
 As for the $Z_C$ above, we establish relations
of $X_C$ for $C$ of small codimension, with similar varieties for Levi subalgebras in $\gen$
but corresponding to the smallest face $0$. . 
In fact, the variety $X_0$ appears as the total space of a natural coherent sheaf on $T_\gen$
whose fibers, generically  are   Cartan subalgebras in $\gen$ but can have higher dimension
over certain degenerate points. This gives rise to a natural partial desingularization
which we call {\em Cartanization}. 

In \S 7 we study the first new case $\gen=\sel_3$. Remarkably, in this case,
 Cartanization of $T_\gen$
gives a smooth variety which is identified with Schubert's {\em variety of complete triangles}
(Proposition \ref{prop:schb-cart}). Therefore,  the entire flop diagram for $\sel_3$ can be
interpreted in terms of Schubert's space. We perform a detailed geometric study
of various intermediate fiber products (2-ray and 3-ray varieties)
 which appear in composing the elementary correspondences \cite{riche}.
  These partial fiber products
 are related 
 to {\em varieties of partial triangles} which form classical stepping stones to $\wh T$,
 see \cite{roberts} and references therein. 
 
 Finally, in \S 8 we establish the main properties of the flober diagram for $\sel_3$
 (Theorem \ref{thm:main}). 
 We prove all the conditions of the Definition \ref{def:H-flober} defining $\Hc$-schobers except
 one which concerns a diagram which we call the 
   {\em Schubert transform}. It connects two $1$-ray fiber products of flopped Grothendieck resolutions
   via the Cartanized variety $Y_0$  (the central term in the  $\sel_3$-flop diagram,
   which is the total space of a vector bundle over $\wh T$).
  Studying the effect of the Schubert transform on derived categories 
   might help to understand $\sel_n$-schobers for higher $n$.

 \paragraph {E. Acknowledgements.} 
 We are grateful to R. Bezrukavnikov, A. Bodzenta, T. Bridgeland, W. Donovan,  A. Efimov, 
 A. Kuznetsov, Yu. Prokhorov, B. Toen and M. Reid for useful 
 discussions and correspondence.  We are also grateful to the referees for numerous remarks which
 have improved the paper. 
 
 The work of A.B. and M.K. was  supported by the World Premier International Research Center Initiative (WPI Initiative), MEXT, Japan. A.B. is partially supported by Laboratory of Mirror Symmetry NRU HSE, RF Government grant, ag. № 14.641.31.0001.
 
  V. S.  thanks Kavli IPMU for the support of the visits during which
a part of this work was performed. 
 
 \vfill\eject
 
 \section{Flops and schobers on the disk}\label{sec:flops-and-schobers}

\noindent{\bf A. Reminder on perverse sheaves. } Let $M$ be a
connected  complex manifold and $S=(M_\alpha)$ be a complex 
 analytic stratification of $M$. We assume that each $M_\alpha$ is smooth but the closure $\ol M_\alpha$ can be a singular complex
 analytic space.  In particular,  there is an open dense stratum $M_0$.

  \vskip .2cm

  We denote by $\Perv(M,S)$  the abelian category of $S$-constructible perverse sheaves 
  of $\QQ$-vector
 spaces on $M$, with respect to the middle perversity \cite{BBD}.    That is, we view $\Perv(M,S)$ 
 as an abelian subcategory
 in the triangulated category of $S$-constructible complexes of sheaves on $M$,
 defined as the heart of
 the perverse t-structure. We further normalize the shift of the t-structure  by the following requirement: 
   for $\Fc\in \Perv(M,S)$
the restriction  $\Fc|_{M_0}$ is a local system in degree $0$. 

Thus a perverse sheaf can be viewed as a particular way to extend a local system on $X_0$ to strata
of higher codimension. 

\vskip .2cm

$\Perv(M,S)$ is a Noetherian and Artinian abelian category. In many cases it can be identified with the
category of representations of an explicit quiver with relations.

The simplest case is  $M=\CC$, stratified by $0$ and $\CC-\{0\}$. Denote the corresponding category of perverse
      sheaves by
        $\Perv(\CC,0)$.   Its  descrption by linear algebra
data was first found in \cite{GGM}.  The following alternative description \cite{KS} was the starting point of
the present paper.
 
      \begin{prop}  \label{prop:dirac}
      $\Perv(\CC,0)$ is equivalent to the
      category of diagrams of 
    $\QQ$-vector spaces
     \[
      \xymatrix{
E_-\ar@<-.7ex>[r]_{\delta_-}& E_0
\ar@<-.7ex>[l]_{\gamma_-}
\ar@<.7ex>[r]^{\gamma_+}& E_+
\ar@<.7ex>[l]^{\delta_+}   
}
\]
satisfying the two following conditions:
\begin{enumerate}
\item[(1)]  $\gamma_-\delta_- = \on{Id}_{E_-}, \,\,\,\gamma_+
\delta_+ = \on{Id}_{E_+}$.
\item[(2)] The maps $\gamma_-\delta_+: E_+\to E_-, \,\,\, \gamma_+\delta_-: E_-\to E_+ $
are invertible. \qed
\end{enumerate}

\end{prop}

Explicitly,  for $\Fc\in\Perv(\CC,0)$ the corresponding diagram is found as follows. 
 Let
$\Rc= \ul\HH^1_\RR(\Fc)$,  the sheaf of  the $1$st hypercohomology with support on $\RR\subset \CC$.

 This is a  sheaf on the real  line $\RR$,
constructible with respect to the stratification into $\RR_{<0}$, $0$ and $\RR_{>0}$.
The spaces $E_-, E_0, E_+$ are the stalks of 
  $\Rc$ at these strata. The maps $\gamma_{\pm}$ are the generalization maps for $\Rc$. 
  The maps $\delta_{\pm}$ can be obtained by considering the dual perverse
  sheaf $\Fc^*$, see \cite{KS}. 
  
  \vskip .2cm
  
  The isomorphism
  \[
  T \,\,=\,\, (\gamma_+\delta_-) (\gamma_-\delta_+): E_+\lra E_+
  \]
  is the monodromy of the local system $\Fc|_{\CC^*}$ calculated at any point of $\RR_{>0}$. 
  
  \vskip .2cm
  
  Further, the global
  (hyper)cohomology    and the (hyper)cohomology with compact support   of $\Fc$
are given by the following two complexes:
\be\label{eq:rgamma}
\begin{gathered}
R\Gamma(\CC,\Fc)\,\,=\,\, \bigl\{  
E_-\oplus E_+ \overset{(\delta_-, \delta_+)}\lra  E_0\bigr\}, 
 \\
R\Gamma_c(\CC, \Fc) \,\,=\,\,\bigl\{  E_0 \overset{(\gamma_-, \gamma_+)} \lra 
E_-\oplus E_+ \bigr\},
\end{gathered}
\ee
concentrated in degrees  $[0,1]$ and $[1,2]$, respectively.  The canonical map
\[
R\Gamma_c(\CC, \Fc) \lra R\Gamma(\CC, \Fc),
\]
is the unique morphism of the above complexes, identical on $E_0$.

     \paragraph {\bf B. The Atiyah flop and its flober.} We want to compare Proposition   \ref{prop:dirac}
     with the simplest example of a {\em flop} in birational geometry. 
     We will work over an algebraically closed field $\k$. For an  quasi-projective scheme $X$   over
     $\k$ we denote by $D^b(X)$ the bounded derived category of {\em coherent} sheaves on $X$,
     and by $\Perf_X$ the triangulated category of perfect complexes, so $\Perf_X\simeq D^b(X)$
     if $X$ is smooth and quasi-projective.  
     When needed, we will consider $D^b(X)$ and $\Perf_X$ as dg-categories with their standard
     enhancements by the $\RHom$-complexes.

     \vskip .2cm
     
     The classical (Atiyah) flop is the diagram
 \be\label{eq:a-flop}
 (X_-, C_-)  \buildrel f_- \over\lra (Z,z) \buildrel f_+\over\lla (X_+, C_+)
 \ee
formed by two birational desingularizations of a 3-fold $Z$ with one quadratic singular point $z$. 
    The desingularizations $X_\pm$ 
     are {\em small}, i.e., each preimage $C_\pm = f_\pm^{-1}(z)$ is identified with $\PP^1$, in contrast
     with   the blowup
     $q: X_0 = \Bl_z Z\to Z$ which has $q^{-1}(z)=\PP^1\times\PP^1$. 
     In fact, $X_0 = X_-\times_Z X_+$ is the fiber product of $X_-$ and $X_+$, so 
     we have a diagram
\be\label{eq:a-fib-prod}
X_- \buildrel p_-\over\lla X_0 \buildrel p_+\over\lra X_+. 
\ee

\begin{ex}\label{ex:a-flop-loc}
In the simplest  case, $Z$ is the 3-dimensional quadratic cone,
 which we can write in the determinantal form
\[
Z \,\,=\,\,\bigl\{ z\in\Mat_2(\k)| \det(z)=0\bigr\}. 
\]
The varieties $X_{\pm}$ are given explicitly by:
\[
X_+ =\bigl\{ (z,l)\in\Mat_2(\k)\times \PP^1 |\,\, z(l)=0\bigr\}, \quad
X_- =\bigl\{ (z,l)\in\Mat_2(\k)\times \PP^1 |\,\, z^t(l)=0\bigr\}.
\]
Here $\PP^1$ is viewed as consisting of $1$-dimensional subspaces $l$ in $\k^2$, on which
$\Mat_2(\k)$ acts. So $X_\pm$, as abstract varieties,  are
  each identified with the total space of the vector bundle 
$\Oc(-1)^{\oplus 2}$ on $\PP^1$. 
The variety $X_0$ is the total space of the line bundle $\Oc(-1,-1)$ on $\PP^1\times\PP^1$. 
\end{ex}

It was proved in \cite{BO} that $D^b(X_+)$ and $D^b(X_-)$ are equivalent. 
More precisely, each of
the functors (now known as the {\em flop functors}) 
\be\label{eq:flop-flop}
\begin{gathered}
T_{+,-} = Rp_{-*}\circ p_+^*: D^b(X_+)\lra D^b(X_-), 
\\
T_{-,+} = Rp_{+*}\circ p_-^*: D^b(X_-)\lra D^b(X_+), 
\end{gathered}
\ee
 is an equivalence. Therefore the diagram of categories
  \be\label{eq:atiyah-schob}
      \xymatrix{
D^b(X_-) \ar@<-.7ex>[r]_{Lp^*_-}& D^b(X_0)
\ar@<-.7ex>[l]_{Rp_{-*}}
\ar@<.7ex>[r]^{Rp_{+*}}& D^b(X_+)
\ar@<.7ex>[l]^{Lp^*_+}   
}
\ee
satisfies the categorical analogs of the properties (1) and (2) of Proposition \ref{prop:dirac}. 
We  say that \eqref{eq:atiyah-schob} {\em represents a perverse schober on $\CC$}.

We denote this schober by  $\Fen$. Since $\Fen$ comes from a flop, we refer to it as the
{\em flober} associated to \eqref{eq:a-flop}.

\paragraph{C. General flops of relative dimension 1.} 
The Atiyah flop is a particular case of a more general situation \cite{vdB} \cite{BB} which also leads to
  schobers on the disk. We recall  this context  in a partial generality, see \cite{BB} for general case.

 We consider a Cartesian diagram of quasi-projective schemes 
  (also referred to as a {\em  flop})
 \be\label{eq:general-flop}
\xymatrix{&X_0 =  {X_-}\times_Z {X_+} \ar[dl]_{p_-} \ar[dr]^{p_+} & \\ {X_-} \ar[dr]_{f_-} && {X_+} \ar[dl]^{f_+} \\ & Z &}
\ee
 with the following properties:
 \begin{itemize}
 \item[(F1)]  $X_\pm$ and $Z$ are irreducible algebraic   varieties (i.e., reduced schemes),
 with $X_\pm$ smooth and $Z$ normal.
 
 \item[(F2)]   $Z$ is  allowed to have 
  canonical hypersurface singularities. Further, the multiplicity  of each singular point $z\in Z$ is at most 2. 
  That is, in $\wh z$,  the formal neighborhood of   $z$,  the equation of $Z$ can be  brought to the form
  $t_0^2=f(t_1, \cdots, t_n)$.

 \item[(F3)]  The morphisms $f_\pm$ are birational, and each fiber of  both $f_\pm$ has dimension $\leq 1$.
 Further, for a singular point $z\in Z$ as in (F2),  the two preimages $f_\pm^{-1}(\wh z)$ can be identified 
 as formal schemes (but not as formal schemes over $\wh z$) via the involution $t_0\mapsto (-t_0)$.
 
 \item[(F4)]  The exceptional loci of $f_\pm$ have codimension $\geq 2$.

 \end{itemize}
  
 \noindent 
 
 Note the scheme $X_0$ can be singular and even non-reduced. 
    \vskip .2cm

In this setup we have  the flop-flop functors $T_{+, -}$ and $T_{-,+}$  similar to \eqref{eq:flop-flop}. More precisely, 
 It was shown in 
  \cite{BB} that  $T_{+, -}$ and $T_{-,+}$   are equivalences.
 In particular,   for this situation we have a
  schober $\Fen$ as in  \eqref{eq:atiyah-schob}. This schober further  was studied in \cite{BB}.      
    \vskip .2cm
    
    Another important feature of a flop \eqref{eq:general-flop} is the 
    (triangulated) {\em null-category} of the morphism $f_-$  which is the full subcategory in $D^b(X_-)$ defined as 
\be\label{eq:nullcat}
{\cal C}_-  \,\, := \,\, \bigl\{A_-\in D^{b}(X_-)| \,\, R(f_{-})_* (A_-) =0 \bigr\}, 
\ee
and similarly for ${\cal C}_+$.  
\begin{rems}
(a)  For an extension to the case when $X_\pm$ have Gorenstein singularities,
 see \cite{BB}.
 
 \vskip .2cm
(b) If one relaxes the assumption that the dimensions of the fibers of $f_\pm$ are $\leq 1$, then the flop functors
\eqref{eq:flop-flop}
may no longer  be equivalences,  see \cite{namikawa}. 
\end{rems}

\paragraph {\bf D. Flober decategorified.}
For a triangulated category $\Dc$ we denote  by $K(\Dc)= K_0(\Dc)\otimes\QQ$ the
Grothendieck group of $\Dc$ made into a $\QQ$-vector space. 
We denote by $[\Fc]\in K(\Dc)$ the class of an object $\Fc\in\Dc$. 
For an algebraic variety $X$
over $\k$, we write $K(X)=K(D^b(X))$; this is the rational Grothendieck group of coherent
sheaves on $X$.

\vskip .2cm

Applying the functor $K$ to the diagram \eqref{eq:atiyah-schob} of the kind considered in
\S C, we get a diagram of $\QQ$-vector
spaces
\be\label{eq:k-flob}
     \xymatrix{
E_-  =  K(X_-) \ar@<-.7ex>[rr]_{\delta_- = p^*_-}&& E_0 = K(X_0)
\ar@<-.7ex>[ll]_{\gamma_- = p_{-*}}
\ar@<.7ex>[rr]^{\gamma_+ = p_{+*}}&& E_+ = K(X_+)
\ar@<.7ex>[ll]^{\delta_+ = p^*_+}.   
}
\ee
This diagram satisfies the properties of  Proposition \ref{prop:dirac} and so represents a perverse
sheaf on $(\CC,0)$ which we denote $K(\Fen)$. 

\begin{ex}\label{ex:atiyah-decat}
Consider the situation of Example \ref{ex:a-flop-loc}. 
In this case the diagram \eqref{eq:k-flob} has the form
\[
 \xymatrix{
\QQ^2 \ar@<-.7ex>[r]_{\delta_-}& \QQ^4
\ar@<-.7ex>[l]_{\gamma_-}
\ar@<.7ex>[r]^{\gamma_+}& \QQ^2
\ar@<.7ex>[l]^{\delta_+}.   
}
\]
Indeed,  we have
 the diagram
\be\label{eq:Xpm&P1}
\xymatrix{
X_- 
\ar[d]_{\pi_-} & \ar[l]_{p_-} X_0\ar[d]^{\pi_0} \ar[r]^{p_+}& X_+\ar[d]^{\pi_+}
\\
\PP^1 & \ar[l]^{\pr_-} \PP^1\times\PP^1 \ar[r]_{\pr_+}& \PP^1
}
\ee
where $\pi_\pm$ represents $X_\pm$ as the total space of $\Oc_{\PP^1}(-1)^{\oplus 2}$ and
$\pi_0$ represents $X_0$ as the total space of $\Oc_{\PP^1\times\PP^1}(-1,-1)$. 
Therefore 
\[
K(X_\pm) \simeq 
 K(\PP^1)\simeq \QQ^2, \quad K(X_0) \simeq  K(\PP^1\times\PP^1) \simeq \QQ^4. 
\]
 This means that the perverse sheaf $K(\Fen)$ has generic rank $2$.
  \end{ex}
 
 \begin{prop}\label{prop:flober-decat}
 Let $j: \CC^*\hookrightarrow \CC$ be the embedding. Then, in the situation of Example  
 \ref{ex:atiyah-decat}, 
 \[
 K(\Fen) \,\,\simeq\,\,  \QQ_0[1] \oplus \ul\QQ_\CC
 \oplus  j_!(\ul\QQ_{\CC^*}).
 \]
 \end{prop}

  \begin{rems} 
 (a) In particular, the proposition says that the monodromy of the local system $j^* K(\Fen)$ on $\CC^*$
 is trivial. This triviality is a feature of our local case: the varieties $X_\pm$ are non-proper, being
 just neighborhoods of the corresponding $(-1)$-curves $C_\pm$. 
 The self-equivalence  (flop-flop functor) of $D^b(X_\pm)$ underlying the monodromy is, by general results,
 see   \cite{BB}, the
 spherical reflection with respect to the spherical object  $\Oc_{C_\pm}(-1)$. In our local case
 the class  $[\Oc_{C_\pm}(-1)]\in K (X_\pm)$  vanishes, see \eqref{eq:ideal} below. 
 
 \vskip .2cm
 
(b)  The situation will be different in the case of a flop relating proper smooth varieties. 
 For instance, consider the natural compactification of Example \ref{ex:a-flop-loc} that is, the 
  flop between the two resolutions of the {\em projective} 3-dimensional quadratic cone. 
  In this case the class   $[\Oc_{C_\pm}(-1)]$ and therefore the monodromy of local system is
  non-trivial. Indeed, as $X_\pm$ is proper, $K(X_\pm)$ carries    the {\em Euler form} 
  \[
  \langle [\Fc], [\Gc]\rangle \,\,=\,\,\sum _i (-1)^i\dim \Ext^i(\Fc, \Gc). 
  \]
  Let  $D_\pm\subset X_\pm$  be an effective divisor such that  the intersection number $(C_{\pm} . D_{\pm})$
  is equal to $1$. Then 
 $\langle [\Oc_{C_{\pm}}(-1)], [\Oc_{D_{\pm}}]\rangle =1$ and so $[\Oc_{C_{\pm}}(-1)]\neq 0$. 
 
 \vskip .2cm
 
 (c) 
 A way to get perverse sheaves with nontrivial monodromy
  in our local situation is to pass to equivariant derived
 categories with respect to the group $G=GL_2\times GL_2$ acting on the variety $Z\subset \Mat_2(\k)$
 (and thus on the entire flop diagram) 
 by left and right multiplications. 
 The monodromy on  corresponding equivariant $K$-groups  
 is the simplest instance of the Hecke algebra action on the equivariant $K$-theory of
 the Grothendieck resolutions, see \S \ref{sec:groth} below and \cite{chriss-ginz}. 
   \end{rems}

 \paragraph{E. Proof of Proposition \ref {prop:flober-decat}:} 
  We denote
\[
\Lc_\pm(i) = \pi_\pm^* \,\Oc_{\PP^1}(i), \,\, i\in\ZZ; \quad \Lc_0(i,j) = \pi_0^* \, \Oc_{\PP^1\times\PP^1}(i,j), 
\,\, i,j\in\ZZ
\]
the line bundles on $X_\pm$ and $X_0$. 
A basis of $K(X_\pm)$ is given by the classes $[\Lc_\pm(i)]$ with  $i=0,1$ while a basis of $K(X_0)$ is given by
the  classes $[\Lc_0(i,j)]$ with  $i,j=0,1$. 

We start with identifying all the maps in \eqref{eq:k-flob}. 

\begin{prop}\label{prop:flober-mat}
(a) The action of $\delta_\pm$ in our bases is given by
\[
\delta_- [\Lc_-(i)] = [\Lc_0(i,0)], \quad \delta_+[\Lc_+(i)] = [\Lc_0(0,i)], \quad i=0,1. 
\]
(b) The action of $\gamma_\pm$ in our bases is given by the matrices
    \begin{table}[H]
 {\scriptsize
  \renewcommand{\arraystretch}{2}
  \begin{tabular}{ p{1cm} || p{1cm}| p{1cm} | p{1cm}  | p{1cm}  }
 $\gamma_-$  & $[\Lc_0(0,0]$  & $[\Lc_0(1,0)] $ & $[\Lc_0(0,1)]$ & $[\Lc_0(1,1)]$
 \\
 \hline\hline
 $[\Lc_-(0)] $& $1$ &  $0$   &   $2$ & $1$
  \\
  \hline
  $[\Lc_-(1)] $ &   $0$ &  $1$
    &  $-1$ &$0$
  \\
   \end{tabular} \quad\quad 
    \begin{tabular}{ p{1cm} || p{1cm}| p{1cm} | p{1cm}  | p{1cm}  }
 $\gamma_+$  & $[\Lc_0(0,0]$  & $[\Lc_0(1,0)] $ & $[\Lc_0(0,1)]$ & $[\Lc_0(1,1)]$
 \\
 \hline\hline
 $[\Lc_+(0)] $& $1$ & $2$ &  $0$  &  $1$
  \\
  \hline
  $[\Lc_+(1)] $ &  $0$  & $-1$
   &  $1$ &  $0$
  \\
   \end{tabular}

}\end{table}
\end{prop}

\noindent {\sl Proof of Proposition \ref{prop:flober-mat}.} Part (a) reflects the isomorphisms
\[
p_-^* \, \Lc_-(i) \,\simeq \,\Lc_0(i,0), \quad p_+^*\, \Lc_+(i) \,\simeq \,\Lc_0(0,i),
\]
which follow
 from the commutativity of the squares in
\eqref {eq:Xpm&P1}.

\vskip .2cm

Let us focus on (b).  We note, first of all,  that the standard exact Koszul complex on $\PP^1$
  \[
  0\to \Oc_{\PP^1} \to \Oc_{\PP^1}(1)^{\oplus 2} \to \Oc_{\PP^1}(2)\to 0
  \]
  implies, after tensoring with $\Oc(i)$ and pulling back,  the  relation
    \be\label{eq:recursion}
  [\Lc_\pm(i)] -2[\Lc_\pm(i+1)] + [\Lc_\pm(i+2)] \,\,=\,\ 0
  \ee
  in $K(X_\pm)$.

\vskip .2cm

As before, we denote by $C_\pm \simeq \PP^1 \subset X_\pm$ the image of the zero section of the bundle $\Oc(-1)^{\oplus 2}$. 
 We have the tautological section of  the pullback bundle
  $\pi_\pm^*\Oc(-1)^{\oplus 2} = \Lc_\pm(-1)^{\oplus 2}$ on $X_\pm$
  which  vanishes on $C_\pm$. The corresponding Koszul resolution of $\Oc_{C_\pm}$  together with 
  \eqref{eq:recursion}
  implies that $[\Oc_{C_\pm}(j)]=0$ for any $j$.  Denoting by $\Ic_\pm\subset \Oc_{X_\pm}$ the sheaf of ideals
  of $C_\pm$, we conclude that 
  \be\label{eq:ideal}
  [\Ic_\pm\otimes \Lc_\pm(j) ] = [\Lc_{\pm}(j)], \quad j\in\ZZ. 
  \ee
  
  \vskip .2cm
  Now, $X_0$ is identified with the blowup of $X_\pm$ along $C_\pm$ that is,
  \[
  X_0 \,\,\simeq \,\, \Proj \biggl( \bigoplus_{n\geq 0} \Ic_+ ^n\biggr)  \,\,\simeq \,\,  \Proj \biggl( \bigoplus_{n\geq 0} \Ic_- ^n\biggr).  
  \]
 Therefore it carries  line bundles $\Oc(1)_\rel^+$ and $\Oc(1)^-_\rel$ coming from these two 
 representations as $\Proj$ so that
 \be\label{eq:rpO(1)}
 Rp_{\pm *} \Oc(1)^\pm_\rel \,\,\simeq \,\, \Ic_\pm. 
 \ee

\begin{lem}\label{lem:O-1}
  Both line bundles $\Oc(1)_\rel^\pm$ are isomorphic to $\Lc_0(1,1)$.
\end{lem}

\noindent{\sl Proof of the lemma:}  For definiteness, let us consider the ``$-$'' case.  It will be more
convenient for us to identify the dual bundles, i.e., to
identify the  ``tautological''  bundle $\Oc(-1)_\rel^-$ with $\Lc_0(-1,-1)$. 
 Any line bundle on $X_0$ is isomorphic to some $\Lc_0(i,j)$
and  $(i,j)$ can be found by restricting to the zero section: $\Lc_0(i,j)|_{\PP^1\times\PP^1} = \Oc_{\PP^1\times\PP^1}(i,j)$. 
We need to find $(i,j)$ for $\Oc(-1)^-_{\rel}$, i.e., to  find  the restriction of $\Oc(-1)^-_{\rel}$ to   $\PP^1\times\PP^1 =  C_- \times \PP^1$
which is the exceptional divisor of the blowup. 
 The  ``vertical'' slices $\{c\}\times\PP^1$ are the projectivizations
of the fibers of the normal bundle to $C_-$ in $X_-$, so the degree of the restriction of $\Oc(-1)^-_{\rel}$ to them is, by 
definition, equal to $(-1)$. 

Let us look at  a ``horizontal'' slice $C_-\times \{p\}$, where $p\in\PP^1$. It is represented by a rank $1$ subbundle
$\Mc_p\subset \Oc_{C_-}(-1)^{\oplus 2}$. This subbundle is isomorphic to $\Oc(-1)$. Now, 
by definition, the fiber of the ``tautological'' bundle  $\Oc(-1)_{\rel}^-$
at a point $(c,p)$ is identified with the fiber of $\Mc_p$ at $c$. This shows that the restriction of $\Oc(-1)_{\rel}^-$
to $C_-\times \{p\}$ is isomorphic to $\Oc_{C_-}(-1)$. Lemma \ref {lem:O-1} is proved. 

 Applying the projection formula and \eqref{eq:rpO(1)}, we get:

 \begin{cor}
 We have, for any  $i\in\ZZ$:
\[
\begin{gathered}
 Rp_{-*} \, \Lc_0(i,0)\, \simeq \, \Lc_-(i), \quad Rp_{-*} \, \Lc_0(i,1) \,\simeq \, \Ic_-\otimes\Lc_-(i-1),
\\
Rp_{+*} \, \Lc_0(0,i)\,\simeq \, \Lc_+(i), \quad Rp_{+*} \, \Lc_0(1,i)\,\simeq \, \Ic_+ \otimes \Lc_+(i-1). \qed
\end{gathered}
\] 
 
 \end{cor}
 
 Applying  this corollary and invoking 
the identity \eqref{eq:ideal} and the recursion \eqref{eq:recursion}, we obtain the proof of 
   Proposition  \ref{prop:flober-mat}. 
   
   \begin{cor}\label{cor:mon-triv}
   The monodromy of $K(\Fen)$ around $0$ is trivial.
   \end{cor}
   
   \noindent{\sl Proof:}    Indeed,  each of the  half-monodromy maps
\[
\varphi_{-+} = \gamma_+\delta_- : E_- \lra  E_+, \quad \varphi_{+-} = \gamma_-\delta_+ :  E_+  \lra E_-
\]
 is given, in the above bases, by the same matrix
  $\begin{pmatrix} 1& 2  \\ 0&-1
 \end{pmatrix}$. Since this matrix squares to the identity, the monodromy is trivial. \qed
   
   \vskip .2cm
   
   We now recall \cite{beil-gluing} \cite{GGM} the more standard description of $\Perv(\CC,0)$ by linear algebra
   data.
   
   \begin{prop}\label{prop:schrod}
(a)    $\Perv(\CC,0)$ is equivalent to the category of diagrams of finite-dimensional $\QQ$-vector spaces
 \[
     \xymatrix{
\Phi   \ar@<-.7ex>[r]_{v}& \Psi
\ar@<-.7ex>[l]_{u}
    }
\]
such that $T_\Phi= \Id_\Phi-uv$ and $T_\Psi = \Id_\Psi-vu$ are invertible.

\vskip .2cm

(b)  If $\Fc$ corresponds to such a diagram, then
$R\Gamma(\CC, \Fc)$ and $R\Gamma_c(\CC, \Fc)$ are identified with the two term complexes
\[
\bigl\{\Psi \buildrel u\over\lra \Phi\bigl\}, \quad \bigl\{\Phi \buildrel v\over\lra \Psi\bigl\},
\]
situated in degrees $[0,1]$ and $[1,2]$ respectively. \qed
   \end{prop}
   
     In particular, we note the following $(\Phi, \Psi)$-diagrams and     
 the corresponding perverse sheaves:
  \be\label{eq:4-indec} 
  \begin{gathered}
  \bigl\{      \xymatrix{
\QQ   \ar@<-.7ex>[r]& 0
\ar@<-.7ex>[l]
    }\bigr\} \,\sim \,\QQ_0[-1], \quad 
      \bigl\{      \xymatrix{
0   \ar@<-.7ex>[r]& \QQ
\ar@<-.7ex>[l]
    }\bigr\} \,\sim \,\QQ_\CC, 
    \\
   \bigl\{      \xymatrix{
\QQ   \ar@<-.7ex>[r]_{\Id} & \QQ
\ar@<-.7ex>[l]_0
    }\bigr\} \,\,\sim \,\, Rj_* (\ul\QQ_{\CC^*}), 
    \quad 
  \bigl\{      \xymatrix{
\QQ   \ar@<-.7ex>[r]_{0} & \QQ
\ar@<-.7ex>[l]_{\Id}
    }\bigr\} \,\,\sim \,\, j_! (\ul\QQ_{\CC^*}).   
    \end{gathered}
  \ee
  
  \begin{lem}\label{lem:hic-dirac-flober} We have the following identifications of the dimensions: 
\[
 \begin{gathered}
 \HH^0(\CC, K(\Fen))\simeq \HH^1(\CC, K(\Fen)) \,\simeq \QQ, 
 \\
 \HH^1_c(\CC, K(\Fen)) \simeq \HH^2_c(\CC, K(\Fen)) \,\simeq \, \QQ^2. 
 \end{gathered} 
 \]
  \end{lem} 
  
  \noindent{\sl Proof of the lemma:}   Indeed, the identifications of the first line follow immediately from part (a) of 
  Proposition  \ref{prop:flober-mat} together  with \eqref{eq:rgamma}. The identifications of the second
 line  follow from part (b) of  Proposition  \ref{prop:flober-mat}, 
  because the $4\times 4$ matrix obtained by stacking the matrices
  of $\gamma_+$ and $\gamma_-$,  is easily found to have rank $2$. \qed
  
 \vskip .2cm

   As explained in \cite{KS}, \S 9A, to pass from the description of Proposition
   \ref{prop:dirac}  (i.e., from the diagram as in \eqref{eq:k-flob}) to that of
   Proposition \ref{prop:schrod}, we need to form
   \[
   \Phi = \Ker(\gamma_-) \subset E_0, \quad  \Psi = E_+.
   \]
    This shows that in  our case both $\Phi$ and $\Psi$ are $2$-dimensional. 
    
    \vskip .2cm
    
    Comparing  now the second line of the identifications in Lemma  
 \ref{lem:hic-dirac-flober} with part (b) of Proposition \ref{prop:schrod}, 
 we find that the complex $\{ \Psi\buildrel v\over\to \Phi\}$
 formed by two $2$-dimensional vector spaces, has $2$-dimensional cohomology spaces, so  $v=0$. Further,
 comparing the first line of   identifications in Lemma  \ref{lem:hic-dirac-flober} with  the same part (b) of Proposition \ref{prop:schrod}
 we find that the complex $\{ \Phi\buildrel u\over\to \Psi\}$ has $1$-dimensional cohomology spaces, so
   $u$ has rank $1$. This means that  our
 $(\Phi, \Psi)$-diagram is  isomorphic to the direct sum of the first, second and fourth diagrams in 
 \eqref{eq:4-indec}. This finishes the proof of Proposition  \ref{prop:flober-decat}. 

\vfill\eject

\section{Cohomology   of 1-dimensional flobers}
\label{sec:1d-flob}

\paragraph{A. Setup of the problem.}
 In this section we study  categorical analogs of the 2-term complexes
\eqref{eq:rgamma}  for the flober
  \be\label{eq:1d-flob}
  \Fen \,\,=\,\,  \bigl\{  \xymatrix{
D^b(X_-) \ar@<-.7ex>[r]_{Lp^*_-}& D^b(X_0)
\ar@<-.7ex>[l]_{Rp_{-*}}
\ar@<.7ex>[r]^{Rp_{+*}}& D^b(X_+)
\ar@<.7ex>[l]^{Lp^*_+}   
}
\bigr\}
\ee
  discussed in  \S \ref{sec:flops-and-schobers}C
from which we keep the assumptions and notations.  

\vskip .2cm

A literal interpretation would give
 2-term complexes 
of triangulated categories  (i.e., just  exact functors)
\be\label{eq:RGamma-flobers}
\begin{gathered}
 R\Gamma(\CC, \Fen) \,\,=\,\,\bigl\{ D^b (X_-)\oplus D^b (X_+) 
\buildrel (Lp_-^*, Lp_+^*)\over\lra D^b (X_0)\bigr\},
\\
R\Gamma_c(\CC, \Fen) \,\,=\,\,\bigl\{ D^b(X_0) \buildrel (Rp_{-*}, Rp_{+*})\over\lra
D^b(X_+) \oplus D^b(X_-)\bigr\},
\end{gathered}
 \ee
 with the sources and targets formally put in degrees $0,1$, resp. $1,2$. 
 We are interested in the  ``cohomology" of these complexes, a concept which
 needs to be defined. 
 
  For a 2-term complex $V^\bullet = \{V^0\buildrel d\over\to V^1\}$ of vector spaces, the cohomology
 of $V^\bullet$ are $\Ker(d)$ and $\Coker(d)$. In our case, the correct replacement of
 ``kernel'' for the ``differential'' $ (Lp_-^*, Lp_+^*)$ in $R\Gamma(\CC, \Fen)$ would be the
 (homotopy) fiber product. Similarly, the ``cokernel'' of $ (Rp_{-*}, Rp_{+*})$ in $R\Gamma_c(\CC,\Fen)$
 is naturally understood as the (homotopy) cofiber product, or pushout. We recall these concepts.
 
 \paragraph{B. Model structures in $\dgCat$.} 
 The proper framework for homotopy limits and colimits is that of model categories  see, e.g., \cite{DHKS}
 \cite{hirschhorn} for background.

 \vskip .2cm
 
  We denote by $\dgCat$ the category of small $\k$-linear dg-categories.
  For $\Ac\in\dgCat$ we denote by $H^0(\Ac)$ the $\k$-linear category with the same objects as $\Ac$
  and the Hom-spaces obtained by taking $H^0$ of the Hom-complexes in $\Ac$. 
  We also denote by
   $\Perf_\Ac$ the category of perfect 
  contravariant dg-functors $\Ac\to\dgVect$. 
  
 We will use two model structures
  on $\dgCat$ introduced by Tabuada
    \cite{tabuada}  (see also  a short summary in
 \cite{DK-triang} \S 1.1).
 The first is the {\em quasi-equivalence model structure},
 where weak equivalences are quasi-equivalences. We denote the category $\dgCat$ equipped
  with the quasi-equivalence model structure
 by $\Qeq$. 
 
 \vskip .2cm
 
  The second is the {\em Morita model structure}, where weak equivalences are
 Morita equivalences. Recall that a dg-functor $f:\Ac\to\Bc$ is called a {\em Morita equivalence},
 if $f_*: \Perf_\Ac\to\Perf_\Bc$ is a quasi-equivalence.

 We denote $\dgCat$  equipped with the Morita model structure
 by $\Mrt$. We recall that fibrant objects in $\Mrt$ are {\em perfect dg-categories} $\Ac$,
 i.e., pre-triangulated \cite{BK-enhanced} dg-categories $\Ac$ such that the  associated triangulated category 
 $H^0(\Ac)$ is closed under direct summands. 
 In particular, for a quasi-projective scheme $X$,  the dg-categories
 $D^b(X)$ and $\Perf_X$ are perfect.  
 
 \vskip .2cm
 
 As shown in \cite{tabuada}, 
 the  model category $\Mrt$ is the left Bousfield localization of $\Qeq$ with respect to the class
 of Morita equivalences. This implies:
 
 \begin{prop}
 (a) Cofibrations in $\Mrt$ are the same as in $\Qeq$.
 
 (b) A dg-functor $f: \Ac\to\Bc$ between perfect dg-categories is a fibration in $\Mrt$ iff it is a fibration in $\Qeq$.
 
 (c)  A dg-functor $f: \Ac\to\Bc$  between perfect dg-categories is a Morita equivalence
  iff it is a quasi-equivalence. For perfect dg-categories being a quasi-equivalence is equivalent to 
 $H^0(f): H^0(\Ac)\to H^0(\Bc)$ being an equivalence of triangulated categories.
 \end{prop}
 
 \noindent {\sl Proof:} (a) is the definition of left Bousfield localization, see \cite{hirschhorn}, Def. 3.3.1.
 Part (b) is \cite{hirschhorn}, Prop. 3.3.16(1). Part (c)  holds because for a perfect $\Ac$ the Yoneda
 functor $y: \Ac\to\Perf_\Ac$ is a quasi-equivalence. \qed
 
 \vskip .2cm
 
 Each of the two model structures on $\dgCat$ gives the corresponding concept of homotopy limits
  and colimits which we will denote by $\holim^\Mrt$ etc. 
   Homotopy (co)fiber products are particular cases of these.

 \paragraph {C.  Global sections.} 
 We define the category of ``global sections'' of the flober $\Fen$  
 as the homotopy fiber product in $\Mrt$:
  \be\label{eq:H^0-flober}
  \begin{gathered}
 \HH^0(\CC, \Fen) \,\,:=\,\, D^b(X_-)\times^h_{D^b(X_0)} D^b(X_+) \,\,=
 \\
 =\,\,
 \holim^{\Mrt}\bigl\{ D^b(X_-) \buildrel Lp_-^*\over\lra D^b(X_0)
 \buildrel Lp_+^*\over\lla D^b(X_+) \bigr\}. 
 \end{gathered}
 \ee
   
 \begin{thm}\label{thm:H^0-flober} 
 $\HH^{0}(\CC, \Fen)$ is identified with
$\Perf_Z$, where $Z$ is as in \eqref{eq:general-flop}. 
\end{thm}

\noindent{\sl Proof:} 
  At the level of triangulated categories,  each of the $Lp_\pm^*$ embeds $D^b(X_\pm)$ into $D^b(X_0)$
 as a full triangulated subcategory. This is because the composition $R(p_\pm)_*\circ Lp_\pm^*$
 is identified with the identity of $D^b(X_\pm)$. Indeed, according to Proposition
 2.7 in  \cite{BB},  $Rp_{\pm *} \Oc_{X_0} = \Oc_{X_{\pm}}$. 
 
  We note further that $D^b(X_\pm)\simeq \Perf_{X_\pm}$
 because $X_\pm$ are assumed smooth. 
 
 \vskip .2cm
 
\noindent  Let us now consider $D^b(X_0)$ as a dg-category, and $\Perf_{X_\pm}$ as dg-subcategories there.
 Let us further consider each  $\Perf_{X_\pm}$  as a {\em homotopy strictly full dg-subcategory} in
 $D^b(X_0)$. That is, we assume that
  it includes, with each object,
all  objects quasi-somorphic to it. Doing this does not change 
the quasi-equivalence classes of our dg-categories. However, it assures that
the embeddings $\alpha_\pm: \Perf_{X_\pm}\to D^b(X_0)$ are {\em fibrations}  in $\Qeq$,
as follows from the explicit description of fibrations in \cite{tabuada}, Prop. 1.13 which we recall.

\begin{prop}
A dg-functor $g: \Ac\to\Bc$ is a fibration in $\Qeq$, if and only if it is surjective on Hom-complexes and the following 
 lifting condition
holds:
\begin{itemize}
\item [(F)] Let $A\in\Ob(\Ac)$ and $u\in\Hom^0_\Bc(g(A), B))$ be a closed degree $0$ morphism in $\Bc$ which
becomes an isomorphism in the category $H^0(\Bc)$. Then there is $A'\in\Ob(\Ac)$ such that $g(A')$ is 
equal (not just isomorphic) to $B$ and there is a closed $u'\in\Hom^0_\Ac(A, A')$ with $g(u')=u$ which
becomes an isomorphism in $H^0(\Ac)$.\qed
\end{itemize}
\end{prop}

Because all the three categories in question are perfect, the embeddings $\alpha_\pm$ are
also fibrations in $\Mrt$.  
 So  we have replaced the  original $Lp_\pm^*$ by fibrations and therefore
\[
\HH^0(\CC, \Fen) \,\,=\,\, \Perf_{X_+} \cap \Perf_{X_-}
\]
(intersection of homotopy strictly full subcategories in $D^b(X_0)$). So our statement reduces to:

\begin{prop}\label{prop:perf-intersection}
$ \Perf_{X_+} \cap \Perf_{X_-}$ is equal to $\Perf_{Z}$, embedded via $L(f_+ p_+ = f_- p_-)^*$
and extended to a homotopy strictly full subcategory. 
\end{prop}

\noindent{\sl Proof of the proposition:} 
  It is clear that the $\Perf_{Z}$  is contained in the intersection, so only the opposite inclusion needs to be proved. 
  Since our statement is about inclusion of classes of objects, we will work in the context of triangulated,
  not dg, categories.

  Suppose  $A_0\in D^b(X_0)$ lies in the  intersection, so there are
  two objects $A_-\in \Perf _{X_-}$ and
   $A_+\in \Perf_{X_+}$ such that $A_0\simeq Lp^*_-(A_-)\simeq Lp^*_+(A_+)$.
  Without loss of generality we can assume that $A_0 =  Lp^*_-(A_-)$. 
  Define
   \[
  A_Z \,\, =\,\, Rf_{-*}  \, R p_{-*} (A_0) \,\, =\,\, Rf_{-*} (A_-) \,\,\in\,\,  D ^b (Z). 
  \]
  The isomorphism $\phi : Lp^*_-(A_-)\simeq Lp^*_+(A_+)$ implies that $ A_Z\simeq R (f_{+})_*A_+$.
  Now, because of the possibly singular nature of $Z$ and $X_\pm$ the functors $Lf_\pm^*$ take values in
  the  left unbounded derived categories $D^-(X_\pm)$, and we can form the exact triangles in these
  categories  induced by the adjunction morphisms:
\[
\begin{gathered}
Lf_-^*(A_Z)\lra  A_- \lra C_-,
 \\
Lf_+^*(A_Z)\lra   A_+ \lra C_+. 
\end{gathered}
\]
Note that $C_\pm \in {\cal C}_\pm$  are objects of the $D^{-}$ versions of the null-categories. 
Further, the morphisms $\phi$  and $\Id$ (the latter expressing the commutativity of the diagram 
\eqref {eq:general-flop}), fit into a commutative square which, by axioms of triangulated categories,
extends to a morphism of exact triangles: 
\[
\xymatrix{
Lp_-^*Lf_-^*(A_Z)
\ar[d]_{\Id} 
 \ar[r] &   Lp_-^* (A_-) \ar[r]
 \ar[d]^{\phi}
 &    Lp_-^*(C_-)
 \ar@{..>}[d]
\\
Lp_+^*Lf_+^*(A_Z)\ar[r] &   Lp_+^* (A_+) \ar[r]& Lp_+^*(C_+). 
}
\]
This morphism is an isomorphism since $\phi$ and $\Id$ are. 
 That is,  $Lp_-^*(C_-)\simeq Lp_+^*(C_+)$.  It follows that 
  \[
  T_{-,+}(C_- )\,\,=  \,\, R(p_{+})_* Lp^*_-(C_-) \,\,\simeq \,\,C_+, 
  \]
  and similarly $T_{+,-}(C_+)\simeq C_-$.
This means that 
  the flop-flop functor $F= T_{+,-}\circ T_{-,+}$ takes  $C_-$ into itself.
  
  On the other hand, it is known that 
   $F$ acts on the null-category ${\cal C}_-$ by the shift by 2 (see \cite{BB} Th. 5.7).
   This implies that $C_-=0=C_+$. Therefore, 
   \[
   A_-\simeq Lf_-^*(A_Z), \quad A_+\simeq Lf_+^*(A_Z) , \quad A_0\simeq Lp^*_-Lf_-^*(A_Z),
   \]
   i.e. $A_0$  is isomorphic to an object in the image of $D ^b (Z)$ in $\Perf _{X_0}$.

It remains to show that $A_Z$ lies in $\Perf_Z$. Since $Lf^*_-(A_Z)=A_-$, we have by adjunction:
\[
R{\cal H}om_Z(A_Z, {\cal O}_z)\,\,=\,\, R(f_{-})_* \,\,  R{\cal H}om_{X_-} (A_-, {\cal O}_x ),
\]
for any closed point $x\in X_-$ and its image $z\in Z$. 
(Here we used the obvious fact that $R(f_-)_* \Oc_x = \Oc_z$.)
Further, since the object in the right hand side is clearly in $D^b(Z)$ and  the morphism $f_-$ is surjective,
it follows that $R{\cal H}om_Z(A_Z, {\cal O}_z)\in D^b(Z)$ for any closed point $z\in Z$.  This implies that $A_Z$ is perfect. 

Proposition \ref{prop:perf-intersection} and Theorem \ref {thm:H^0-flober} are proved. \qed

\begin{rem}
The only essential geometric feature of the situation used in the proof of Theorem \ref {thm:H^0-flober}
is that $T_{+,-}\circ T_{-,+}$ acts on the null-category $\Cc_-$ by a shift.  So the argument can possibly
generalize to other situations. 
\end{rem}

\paragraph{D.  $\HH^1_c$ and $\HH^1$.} 
We define
\[
\HH^1_c(\CC, \Fen) \,\,=\,\, \Kc\,\, :=\,\, \Ker \bigl\{ D^b(X_0) \buildrel (Rp_{-*}, Rp_{+*})\over\lra
D^b(X_+) \oplus D^b(X_-)\bigr\},
\]
i.e., as the full subcategory in $D^b(X_0)$ consisting of objects $A$ such that $Rp_{-*}(A)$
and $Rp_{+*}(A)$ are quasi-isomorphic to $0$. This category has been used in \cite{BB}
to construct a spherical pair.

We further define
\[
\HH^1(\CC,\Fen) \,\,=\,\, D^b(X_0) \bigl/ \la Lp_+^* (D^b(X_+)), Lp_-^*(D^b(X_-))\ra,
\]
the quotient by the minimal thick subcategory generated by the 
pullbacks of the $D^b(X_\pm)$. 

\begin{rem}
Although $\HH^0(\CC, \Fc)$ is defined as the homotopy  limit of the diagram of categories
\[
\bigl\{ D^b(X_-) \buildrel Lp_-^*\over\lra D^b(X_0)
 \buildrel Lp_+^*\over\lla D^b(X_+) \bigr\},
 \]
the category  $\HH^1(\CC, \Fc)$ is not defined as the homotopy colimit of the same diagram: that homotopy colimit is
 $D^b(X_0)$. Our definition, which is a direct categorification of the cokernel of the map of vector spaces of the form
 $E_- \oplus E_+ \to E_0$, can be perhaps seen as some ``derived functor of the homotopy limit''.  
\end{rem}

\begin{ex}
 For the Atiyah flop (Example \ref{ex:a-flop-loc}), the category $\Kc$ is equivalent to $D^b(\Vect_\k)$
 and is generated by the sheaf $\Oc_{\PP^1\times\PP^1} (-1,-1)$, where $\PP^1\times\PP^1$ is the 
 exceptional fiber in $X_0$. The category $\HH^1(\CC,\Fen)$ is in this case also identified with
 $D^b(\Vect_\k)$. 
 
 In general,  $\Kc = \HH^1_c(\CC,\Fen)$, if considered as a dg-category, is not smooth. 
 The relation between it and $\HH^1(\CC,\Fen)$ probably fits into the ``categorical Poincar\'e duality'',
 see Remark \ref{rem:CPD}. 

\end{ex}

 \paragraph {E. $\HH^2$ with compact support} 
We define $\HH^2_c(\CC, \Fen)$ as the homotopy pushout in $\Mrt$: 
\be\label{eq:H-2-c}
\begin{gathered}
\HH^2_c(\CC, \Fen) \,\,:=\,\, D^b(X_-)\bigsqcup^h_{D^b(X_0)} D^b(X_+) \,\,=
\\
=\,\,
\hocolim^\Mrt\bigl\{ 
D^b(X_-) \buildrel R(p_-)_*\over\lla D^b(X_0)
 \buildrel R(p_+)_*\over\lra D^b(X_+) \bigr\}.
 \end{gathered}
\ee

\begin{thm}\label{thm:H^2-c} 
 $\HH^2_c(\CC, \Fen)$  
is identified with $D^b(Z)$, where $Z$ is as in \eqref{eq:general-flop}. 
\end{thm}

\begin{rem}\label{rem:CPD} 
Suppose that $Z$ and $X_\pm$ are projective. 
The relations $\HH^0(\CC,\Fen)$ and $\HH^2_c(\CC,\Fen)$, i.e., between
$\Perf(Z)$ and $D^b(Z)$, is then that of duality. More precisely, let us consider them as dg-categories.
To\"en \cite{toen} has introduced a duality operation on dg-categories called the {\em Morita duality}:
\[
\Ac \,\,\mapsto  \Ac^\vee \,\,=\,\, \ul{\RHom} (\Ac, \Perf_\k),
\]
where $\ul\RHom$ is an approproately defined internal Hom functor in $\dgCat$. 
Here $\Perf_\k$ is the dg-category of cochain complexes over $\k$ with bounded, finite-dimensional
cohomology (a particular case of the notation $\Perf_\Ac$ from \S \ref{sec:1d-flob}B).

It follows from \cite{BvdB}, Thm. A.1,  that for any projective variety $Z$ we have 
$\Perf_Z^\vee\simeq D^b(Z)$ and so 
\[
\HH^2_c(\CC, \Fen) \,\,\simeq \,\, \HH^0(\CC, \Fen)^\vee. 
\] 
Note that a perverse schober $\Fen$ can be seen as a categorical analog of  a perverse sheaf $\Fc$
 with an identification $\Fc \buildrel\sim\over\to\Fc^*$, where $\Fc^*$ is the Verdier dual
perverse sheaf. This corresponds to the $\Hom$-pairing  on any category. So the above
identification can be seen as an instance of a {\em categorical Poincar\'e duality}. 

\end{rem}

  \noindent{\sl Proof of Theorem \ref{thm:H^2-c}:} 
  We will use the concept of the Drinfeld quotient $\Cc/\Ac$ of a dg-category $\Cc$ by a full dg-subcategory $\Ac$,
  see \cite{Dri}. This construction is invariant under quasi-equivalences in the following sense. Suppose
  \[
  \xymatrix{\Ac\ar[d]_\alpha\ar[r]&\Cc\ar[d]^\beta
  \\
  \Bc\ar[r]&\Dc 
  }
  \]
  be a commutative diagram of dg-functors with horizontal arrows being fully faithful  embeddings and $\alpha$ and $\beta$ being
  quasi-equivalences. Then we have a quasi-equivalence of Drinfeld quotients $\Cc/\Ac\to\Dc/\Bc$
  compatible with $\alpha$ and $\beta$. 
  
  \vskip .2cm
  
   To relate this with our geometric situation, we prove the following theorem of independent interest.

  \begin{thm}\label{thm:quotient}
Let $f: X\to Y$ be a projective morphism of quasi-projective schemes over $\k$ such that each
geometric fiber of $f$ has dimension $\leq 1$ and $Rf_*(\Oc_X)\simeq \Oc_Y$. Let
$\Cc_f\subset D^b(X)$ be the kernel of $Rf_*$.
 Let us endow $D^b(X), D^b(Y)$ and $\Cc_f$ with natural structures of dg-categories.
 Then $Rf_*$ lifts to a quasi-equaivalence between $D^b(Y)$ and the Drinfeld quotient $D^b(X)/\Cc_f$.
 \end{thm}
 
 The proof will be given in \S F below. 
 
\begin{rem}  Theorem  \ref{thm:quotient} has been stated in \cite{BO2} (in the equivalent terms of Verdier rather than Drinfeld quotient) with no restriction on the dimension of fibers of $f$. The authors of \cite{BO2} later realized that they do not have the proof of this statement. To the best of our knowledge, the problem remains unsolved for this general case. 
\end{rem}

We  apply Theorem  \ref{thm:quotient} to the two rows in the  diagram of dg-categories
\[
\xymatrix{
\Lc_+\ar[d]_{h} \ar[r]&  D^b(X_0)\ar[d]^{Rp_{-*}} \ar[r]^{ Rp_{+*}}  & D^b(X_+)\ar[d]^{Rf_{+*}}
\\
\Cc_- \ar[r] &D^b(X_-) \ar[r]^{ Rf_{-*}} & D^b(Z),
}
\]
Here $\Cc_-$ is the null-category \eqref{eq:nullcat}. We conclude that $D^b(Z) \simeq D^b(X_-)/\Cc_-$
(Drinfeld quotient). Similarly,
 $\Lc_+$ is the kernel of $Rp_{+*}$, and we conclude that 
  $D^b(X_+)=D^b(X_0)/\Lc_+$. The dg-functor $h$ is induced by $Rp_{-*}$. 
  
  \vskip .2cm

  We now note the following. 
  
  \begin{lem}\label{lem:drin-pushout}
  Let $g: \Ac\to\Bc$ be any dg-functor. Then the homotopy pushout
   $\Bc\bigsqcup_\Ac^h 0$ (both in $\Qeq$ and $\Mrt$)
  is quasi-equivalent to the Drinfeld quotient $\Bc/\Imm(g)$, where $\Imm(g)$ is the full dg-subcategory on objects from
  the image of $g$. 
  \end{lem}
  
  \noindent{\sl Proof:}  In any model category, the homotopy pushout of a diagram of cofibrant objects can be obtained
  by replacing one of the functors by a cofibration and then taking the usual push out. 
  So we  first replace $\Ac$ and $\Bc$ by cofibrant (in $\Qeq$ or $\Mrt$) dg-categories, and, after changing the
  notation, assume that $\Ac$ and $\Bc$ are already cofibrant. Second, we replace the dg-functor $\Ac\to 0$
  by $\Ac\to\Ac/\Ac$, which is a cofibration (because the natural functor  to any Drinfeld quotient is a cofibration). 
  Now, the usual pushout $\Bc \bigsqcup_\Ac (\Ac/\Ac)$ in $\dgCat$ is quasi-equivalent to $\Bc/\Imm(g)$. \qed

  \begin{lem}\label{lem:two-pushots}
 Let $u: \Sc_1\to\Tc_1$ be a dg-functor of perfect dg-categories such that $H^0(u): H^0(\Sc_1)\to H^0(\Tc_1)$
 is fully faithful and let $\Tc_1\to\Tc_2$ be any dg-functor. Let $\Sc_2\subset \Tc_2$ be the strictly full dg-subcategory on objects
 from the image of $\Sc_1$. Then the homotopy pushout $\Tc_2\bigsqcup^h_{\Tc_1} (\Tc_1/\Sc_1)$
 is identified with $\Tc_2/\Sc_2$. 
\end{lem}

\noindent {\sl Proof:}  Using   associativity of pushouts and Lemma \ref  {lem:drin-pushout},   we rewrite:  
\[
\Tc_2 \bigsqcup^h_{\Tc_1} (\Tc_1/\Sc_1)\,\, 
 \simeq \,\, \Tc_2 \bigsqcup^h_{\Tc_1}
 \biggl( \Tc_1\bigsqcup^h_{\Sc_1} 0\biggr)
 \,\,\simeq \,\,
 \biggl(\Tc_2\bigsqcup^h_{\Tc_1} \Tc_1\biggr)\bigsqcup_{\Sc_1} 0  
 \,\,  \simeq \,\,   \Tc_2\bigsqcup^h_{\Sc_1} 0 \,\,  \buildrel \ref  {lem:drin-pushout} \over \simeq \,\,
\Tc_2/\Sc_2. 
\]
\qed

\vskip.2cm

To apply Lemma  \ref {lem:two-pushots} , we take $\Tc_1=D^b(X_0)$, $\Sc_1=\Lc_+$ and $\Tc_2= D^b(X_-)$. 
This gives 
\[
D^b(X_-)/
\Sc_2 \,\,\simeq \,\, D^b(X_-) \bigsqcup^h_{D^b(X_0)} ( D^b(X_0)/\Lc_+),
\]
 where $\Sc_2$ is the strictly full subcategory on objects from the image of $\Lc_+$, i.e.,
the strictly full subcategory in $C_-$ on the objects from the image of the functor $h$. 
By  Theorem   \ref{thm:quotient}  we have $D^b(Z) = D^b(X_-)/\Cc_-$ and $D^b(X_+) = D^b(X_0)/\Lc_+$. 
Thus to finish the proof of 
Theorem \ref{thm:H^2-c},  it suffices to establish   the following fact.

\begin{lem}\label{lem:h-surjective} $\Sc_2=\Cc_-$, i.e.,
the functor $h$ is essentially surjective on objects. 
\end{lem}

\vskip .2cm

 \noindent{\sl Proof:} Let $\Ac_- = \Cc_-\cap \Coh(X_-)$, an abelian category.
 It was shown in \cite{BB} that:
 
 \begin{enumerate} 
 
  \item[(1)] There is an essentially surjective (spherical) functor $\Psi_-: D^b(\Ac_-)\to \Cc_-$
which extends the embedding of $\Ac_-$.

\item[(2)] There is an exact  functor $\wt p_-: D^b(\Ac_-)\to D^b(X_0)/\Kc$ such that
 $\Psi_- = R(p_-)_*\circ \wt p_-$. Here   $\Kc$
is the kernel in \S D above.

\item[(3)] We also have  $R(p_+)_*\circ \wt p_-=0$. 
\end{enumerate}
Our statement follows from these properties. \qed
 
 \vskip .2cm
 
 This finishes the proof of Theorem 
 \ref{thm:H^2-c}.

\paragraph{F.  Proof of Theorem \ref{thm:quotient}. }
 Note that the homotopy category of the Drinfeld quotient by a full dg-subcategory whose homotopy category is a thick subcategory of the ambient homotopy category is the Verdier quotient of the corresponding homotopy categories. For this reason, it is enough to prove the equivalence of $D^b(Y)$ regarded as triangulated category,  with the Verdier quotient 
$D^b(X)/\Cc_f$. 

\vskip .2cm

We have the obvious functor $\Phi :D^b(X)/{\cal C}_f\to D^b(Y)$.
We will prove that $\Phi$ is an equivalence. This statement consists of three steps (I)-(III) below.

\vskip .2cm

\noindent \ul  {\sl (I) $\Phi$ is injective on $\Hom$-spaces.} 
 Take a pair of objects $A,B\in D^b(X)$ and a morphism $\alpha : A\to B$. 
Suppose that $\Phi (\alpha )=Rf_*(\alpha )=0$. 
As in the proof of Proposition \ref{prop:perf-intersection}, we consider the functor $Lf^*$ taking values in
$D^-(X)$ and form
 an exact triangle with the first arrow the counit $c$ of the adjunction:
\begin{equation}\label{ffac-triangle}
Lf^*Rf_* A \buildrel c\over \lra  A\lra C, \quad C\in D^-(X). 
\end{equation}
By applying $Rf_*$ to this triangle, we see that $C$ belongs to ${\cal C}_f^-$, the $D^-$ version of the null-category for $f$.
The composition of $c$  with $\alpha$ gives the morphism $\gamma: Lf^*Rf_*A \to B$ that corresponds by adjunction to $Rf_*\alpha$, hence $\gamma$ is zero. Then the long exact sequence on Hom-spaces from the triangle
\eqref{ffac-triangle}  to $B$ shows that $\alpha$ factors via a morphism $C\to B$. 
If $C$ were in ${\cal C}_f$ (the bounded version of the null category),
 then this would imply that $\alpha$ is the zero morphism in the quotient category.
 
 \vskip .2cm

As $C$ is, in general, an object of  $\Cc_f^-$ only, we proceed to 
 find an appropriate bounded replacement for $C$. To this end, we consider a t-structure on $D^-(X)$ 
 which restricts to a bounded t-structure on $D^b(X)$ and such that functor $Rf_*$ is t-exact for this t-structure. 
For instance, we can choose any of the Bridgeland's t-structures on the unbounded category $D(X)$ \cite{bridg-flop}, which are numbered by $p\in \mathbb{Z}$, by fixing $p$,   and restrict it to $D^-(X)$. 

Consider the exact triangle of truncations for $C$ in this t-structure:
\[
\tau_{<k}C\to C \to \tau_{\ge k}C, \quad k\in\ZZ. 
\]
Since $B$ is in the bounded derived category, we can find $k \ll 0$ such that ${\rm Hom}(\tau_{<k}C, B)=0$. Then, applying functor ${\rm Hom}( -,B)$ to the above triangle of truncations gives a lift of our morphism 
$C\to B$ to a morphism $\tau_{\ge k}C\to B$. Together with the composite
 $A\to C \to \tau_{\ge k}C$, this gives a factorization of $\alpha$ via $\tau_{\ge k}C$. 
 Since $Rf_*$ is t-exact and $C$ is in ${\cal C}^-_f$, it follows that $\tau_{\ge k}C\in {\cal C}_f$.
  Hence $\alpha =0$ in the quotient category.
This proves the claim (I).  

\vskip .2cm

\noindent \ul  {\sl (II) $\Phi$ is surjective on $\Hom$-spaces.} 
Let $\beta \in {\rm Hom} (Rf_*A, Rf_*B)$. It induces a morphism ${\tilde \beta}: Lf^*Rf_*A \to B$, 
such that $Rf_*{\wt \beta }=\beta$.
As before, the source of $\wt \beta$ lies, a priori, in $D^-(X)$, and we
consider a triangle of truncations:
$$
\tau_{<k}(Lf^*Rf_*A)\to Lf^*Rf_*A \to \tau_{\ge k}(Lf^*Rf_*A)
$$
with respect to the t-structure as in the proof of (I). 
Again, since $A$ and $B$ are bounded, we can choose $k$ negative enough so that
 \[
 {\Hom}(\tau_{<k}(Lf^*Rf_*A), A)\,\,=\,\, 0\,\, =\,\, {\Hom}(\tau_{<k}(Lf^*Rf_*A), B),
 \]
  which implies that the canonical morphism $Lf^*Rf_*A\to A$ has a lift to a morphism $\delta :\tau_{\ge k}(Lf^*Rf_*A) \to A$ and ${\tilde \beta}$ has a lifting to a morphism $\beta ': \tau_{\ge k}(Lf^*Rf_*A) \to B$. 

By considering the octahedron related to the commutative triangle 

 \[
 \xymatrix{
  Lf^*Rf_*A \ar[r]
 \ar[dr] & \tau_{\ge k}Lf^*Rf_*A 
 \ar[d]^{\delta}
 \\
 & A, 
 }
 \]
we conclude that the cone of the morphism $\delta$ lies in ${\cal C}_f$. Hence we got a pair of morphisms (a roof):
$$
A \buildrel \delta \over\lla \tau_{\ge k}(Lf^*Rf_*A) \buildrel \beta '\over\lra B, 
$$
which can be interpreted as a 
morphism $A \dasharrow B$ in the quotient category. Moreover, one can easily see that it is mapped by $\Phi$ into $\beta$. 
This  proves the claim (II). Together, (I) and (II) mean
that $\Phi$ is is a fully faithful functor.

\vskip .2cm

\noindent \ul  {\sl (III) $\Phi$ is essentially surjective (on objects).} 
 For any coherent sheaf $\Fc$ on $Y$, consider the zeroth cohomology sheaf with respect to our chosen t-structure
 as above:
$$
A={\cal H}^0(Lf^*(\Fc))
$$
Since $Rf_*Lf^*(\Fc)=\Fc$ and $Rf_*$ is t-exact, it follows that $Rf_*(A)=\Fc$, i.e. $\Fc$ lies in the image of the functor.
Since the functor is fully faithful, the d\'evissage of any object $B\in D^b(Y)$ into its t-cohomology sheaves 
proves that  $B$ lies in the essential image of $\Phi$. This ends the proof of Theorem  \ref{thm:quotient}.

\vfill\eject

\section{The web of flops. Role of arrangements. }\label{sec:webflop}

  \noindent {\bf A. Perverse sheaves on  arrangements. }
  Proposition \ref{prop:dirac} is the simplest instance of a more general result
  \cite{KS} which we now recall. 
  
   Let   $\Hc$ be an arrangement of hyperplanes in $\RR^n$.
   The complexification $\Hc_\CC$ defines a natural stratification of $\CC^n$, and we denote by
   $\Perv(\CC^n, \Hc)$ the corresponding category of perverse sheaves. In \cite{KS} it was described
   in terms of the {\em chamber structure} on $\RR^n$ given by the arrangement $\Hc$.
   More precisely, $\Hc$ subdivides $\RR^n$ into locally closed real {\em cells} (of various
    dimensions).
   Open cells are called {\em chambers}. We have a poset $(\Cc = \Cc_\Hc, \leq)$ formed by cells and inclusions
    of their closures.

 \begin{thm}\label{thm:perv-arr}
   $\Perv(\CC^n, \Hc)$ is equivalent to the category of diagrams formed by  
   vector spaces $E_C, C\in\Cc$, and linear maps
   \[
 \xymatrix{
 E_C \ar@<.4ex>[r]^{\gamma_{CC'}}&E_{C'}  \ar@<.4ex>[l]^{\delta_{C'C}}, 
 } \quad \quad C\leq C',
 \]
 satisfying the following relations:  
 \begin{itemize}
 
  \item[(0)] Transitivity:  if $C''\leq C\leq C'$, then $\gamma_{CC''}=\gamma_{C'C''}\gamma_{CC'}$ and
 $\delta_{C''C}=\delta_{C'C}\delta_{C''C'}$.
 
 \item[(1)] Idempotency: 
 $ \gamma_{CC'}\delta_{C'C}=\Id$ for any $C\leq C'$. 
  This, together with (0), allows us to unambiguously define
 the map $\varphi_{CC'}: E_C\to E_{C'}$ for any $C, C'$ as the composition
  $\varphi_{CC'}=\gamma_{DC'}\delta_{CD}$ where $D$ is any cell such that $D\leq C, C'$.

  \item[(2)] Collinear transitivity: if three cells $C_1, C_2, C_3$ are {\em collinear}, i.e., there are 
  points $c_i\in C_i$ such that $c_2$ lies in the closed straight line interval $[c_1, c_3]$, then
  $\varphi_{C_1C_3}=\varphi_{C_2C_3}\varphi_{C_1C_2}$. 
  
  \item[(3)] Invertibility: if $C, C'$ are two cells of the same dimension $r$ lying in a linear subspace of
  dimension $r$ and separated by a cell $D\leq C, C'$ of dimension $(r-1)$ (``wall''), then
  $\varphi_{CC'}$ is an isomorphism. \qed
  \end{itemize}
  \end{thm}
  
 Proposition  \ref{prop:dirac} is obtained when $n=1$ and $\Hc$ consist of one ``hyperplane'' $\{0\}\subset\RR$.

 \vskip .3cm
 
 \noindent {\bf B. Web of flops: the movable cone picture.}  
 General construction of $3$-dimensional flops, of which the Atiyah flop \eqref{eq:a-flop} is the simplest
 instance, can be considered as an algebro-geometric analog of the procedure of surgery on a knot
 in $3$-dimensional topology. That is, we remove a curve $C_+$ from a variety $X_+$ and fill $X_+ - C_+$
 in a new way, by a  curve $C_-$,  to get $X_-$. The curves $C_+$  (in smooth $X$)
 for which this is possible
 ({\em flopping curves}) have been classified. If $C_+$ is irreducible, then $C_+\simeq \PP^1$
 and the normal bundle $N_{C_+/X_+}$ can be one of the three types:  $\Oc(-1)^{\oplus 2}$ (the Atiyah flop),
 $\Oc\oplus \Oc(-2)$ (the so-called {\em pagoda flops}), and $\Oc(1)\oplus\Oc(-3)$, see 
 \cite{pinkham}.

 \vskip .2cm
 
 More generally, we can start with a reducible curve $C_+=\bigcup C_i$ in a 3-fold $X=X_+$ such
 that we can make a flop along some component $C_i$, getting a new 3-fold $X_1$, then  it may be possible
 to flop $X_1$ along the strict preimage of some other component $C_j$, getting $X_2$ and so on. We get in this
 way a ``web of flops'', a system of 3-folds $X_\nu$  and flops connecting them. This system of flops can have
 loops: we may be able to obtain the same $X_\nu$ by two or more different sequences of flops. 
 
 \vskip .2cm

  According to Y. Kawamata \cite{kaw}, the  structure of iterated
  flops is governed by the {\em chamber structure of the movable cone} $M_X$.
 By definition, $M_X\subset\on{Pic}(X)\otimes \RR$ is
 the cone generated by
  line bundles $\Lc$ on $X$ such that the locus of base points of the 
  linear system $|\Lc|$ has codimension $\geq 2$. 
    The open part of this cone is subdivided into {\em chambers}
 (certain open cones).  For $\Lc$ inside each chamber, 
 the variety $X_\Lc$ is the same;  in particular,   the
sign of the  degree  (i.e., the property of being ample or not)
of  the restriction of $\Lc$ to any fixed component of $C_+$  is the same. 
 
 Considering the closures of the chambers, their intersections etc.,
 one gets a decomposition of the movable cone into cells of all dimensions. In particular, 
  when we cross  a {\em wall}  (a cell of real codimension 1)
 between neighboring chambers, the $X_\Lc$ undergoes a
  flop.  Further,  cells of codimension $\geq 2$ can be considered as relations, syzygies etc. among the flops.

 \vskip .3cm
 
 \noindent{\bf C. The elephant picture. The role of arrangements. The schober HMMP. } 
 It is a fundamental fact \cite{BO}\cite{bridg-flop}
  that any  3-dimensional flop results in a  derived equivalence:
 $D^b(X_+)\buildrel\sim\over \to D^b(X_-)$. The {\em Homological Minimal Model Program}
 (HMMP) studies, in particular,   such derived equivalences and relations among them. 
 If we have a web of flops $(X_\nu)$, then
 all $X_\nu$ have equivalent derived categories, and loops give self-equivalences. 
 We get a local system
 of categories on the oriented graph whose vertices are the $X_\nu$ and edges are the flops. 
 In fact, it
 was proved by Donovan-Wemyss \cite{DW} that this extends to a local system of categories on $\CC^n-\Hc_\CC$,
 the complement of a certain hyperplane arrangement on $\CC^n$ with equations
 of the hyperplanes having real coefficients.
     Individual $X_\nu$ correspond thereby to chambers of the real arrangement. 
 
 \vskip .2cm
  More precisely, 
 contracting a (possibly reducible) flopping curve 
 $C_+\subset X_+$ produces a $3$-fold $Z$ with an isolated singular point $z$. 
  According to M. Reid's
 ``elephant'' picture \cite{reid},   hyperlplane sections of $Z$ through $z$ are partial desingularizations of 
  a Kleinian (ADE) singularity $\AAA^2/\Gamma$. 
  So $Z$, near $z$,  can be seen in terms of the total space of a 1-parameter deformation
  of $\AAA^2/\Gamma$, which can be 
   obtained by base change from the universal   deformation  of $\AAA^2/\Gamma$
  as constructed  by  Brieskorn \cite{brieskorn-german}. Brieskorn's deformation is defined over the
  Cartan subalgebra $\hen$ of the corresponding ADE Lie algebra, with singular fibers
  forming the {\em root arrangement} of hyperplanes in $\hen$, and the Donovan-Wemyss arrangement
  is a certain subarrangement of this,  corresponding to a partial, not full, desingularization of 
  $\AAA^2/\Gamma$. 
  
  \vskip .2cm
  
   Therefore,  Theorem
 \ref{thm:perv-arr} appears very suggestive from the point of view of 3-dimensional birational geometry
 of flops. 
 It leads to the following. 
 
 \vskip .2cm
 
  \begin{prpsl}\label{prpsl:HMMP}
    Local systems of triangulated categories appearing in HMMP
 typically defined on open strata of some  stratified spaces, admit natural extensions to
 perverse schobers on the entire spaces. 
 \end{prpsl}

 \vskip .3cm
 
 \noindent{\bf D. The general concept of  $\Hc$-schobers.} 
  Note that  the conditions 
 of Theorem
 \ref{thm:perv-arr}  do not
 involve addition in the relations on the  $\gamma$'s and $\delta$'s and so can be formulated 
 for a diagram of triangulated categories, not
 vector spaces. 
 
 More precisely, let $\Hc$ be an arrangement of hyperplanes in $\RR^n$. We consider the poset
 $(\Cc = \Cc_\Hc, \leq)$ of cells of $\Hc$ as a category. We denote by $C_{\min}\in\Cc$
 the intersection of all the hyperplanes of $\Hc$; one can assume that $C_{\min}=\{0\}$
 without changing the combinatorics, so we make this assumption in the sequel. 
 
 \begin{Defi} 
  By a {\em triangulated 2-functor} on $\Cc$
 we will mean a 2-functor $\Fc$ from $\Cc$ to the 2-category of triangulated categories and exact functors.
 That is, $\Fc$ consists of the data:
 \begin{itemize}
 \item  For each $C\in\Cc$, a triangulated category $\Ec_C$.
 
 \item For each pair  $C\leq C'$ in $\Cc$, an exact functor $\gamma_{CC'}: \Ec_C\to\Ec_{C'}$.
 
 \item  For each length 3 chain $C\leq C'\leq C''$ in $\Cc$, an isomorphism of functors
 $\kappa_{C, C',C''}: \gamma_{C'C''}\gamma_{CC'}\to\gamma_{CC''}$, and these
 isomorphisms  satisfy the compatibilty
 condition for each length 4 chain of cells. 
 \end{itemize}
 \end{Defi} 
 
 Let now $\Fc = (\Ec_C, \gamma_{CC'})$ be a triangulated 2-functor on $\Cc$. Assume, in addition,
 that each $\gamma_{CC'}$ has a right adjoint $\delta_{C'C}$. Then $(\Ec_C, \delta_{C' C})$
 form a triangulated 2-functor on $\Cc^\op$, the opposite poset of $\Cc$. We now impose the following
 additional condition:
 
 \begin{itemize}
 \item[(1)] (Idempotency) Each natural transformation $\Id_{\Ec_C}\to \gamma_{CC'}\delta_{C' C}$
 (the unit of the adjunction), is an isomorphism of functors.
 \end{itemize}
 This implies, first of all, that each $\delta_{CC'}$ is an embedding of a left admissible
 \cite{bondal-kapranov} full triangulated subcategory, so we can 
  think of  all the  $\Ec_C$ as such subcategories in $\Ec_0$, the category correspond to the minimal
  cell $\{0\}=C_\min$. 
  
  Let us define the  {\em flopping functors} 
  \be\label{eq:flop-functors}
  \varphi_{CC'} = \gamma_{0C'} \delta_{C0}:\,\, \Ec_C\lra \Ec_{C'}
  \ee
  for any two cells $C, C'\in\Cc$. The idempotency condition provides the  canonical identification
  $\varphi_{CC'}  \, \simeq \, \gamma_{DC'}\delta_{CD}$ where $D$ is any cell such that $D\leq C, C'$.
  In particular,
  \[
  \varphi_{CC'} = \gamma_{CC'}, \,\, C\leq C'; \quad \varphi_{CC'} = \delta_{CC'}, \,\, C\geq C'.
  \]
  Further, given three arbitrary faces $C, C', C''$, the counit $\delta_{C' 0} \gamma_{0C'}\to\Id_{\Ec_{C'}}$
  induces a natural transformation (not necessarily an isomorphism)
  \[
  \kappa_{CC'C''}: \varphi_{C'C''}\varphi_{CC'}\lra \varphi_{CC''}.
  \]
  For any four cells $C, C', C'', C'''$ the evident square of $\kappa$'s commutes:
  \be\label{eq:kappa-comm}
  \begin{gathered}
  \kappa_{CC'C'''} \circ (\kappa_{C C'' C'''} \cdot \varphi_{CC'}) = \kappa_{C C'' C'''}
  \circ (\varphi_{C''C'''}\cdot \kappa_{C C'' C'''})
  \\
  \text{as transformations }  \, \varphi_{C'' C'''}\varphi_{C'C''}\varphi_{CC'}\lra\varphi_{C C'''}.
  \end{gathered}
  \ee

  Let $[\Cc]$ be the category whose objects are all elements of $\Cc$ and there is a
  unique morphism between any two objects. 
  We can express \eqref{eq:kappa-comm} by saying that $(\Ec_C,\varphi_{CC'}, \kappa_{CC'C''})$
  form a lax 2-functor $[\Cc]\to  \Cc at$.

 \begin{Defi}\label{def:H-flober}
 Let $\Hc$ be an arrangement of hyperplanes in $\RR^n$. An {\em $\Hc$-schober} is
 a triangulated 2-functor $\Fc = (\Ec, \gamma_{CC'}, \kappa_{CC'C''})$ on $\Cc_\Hc$
 with each $\gamma_{CC'}$ admitting  a right adjoint $\delta_{C'C}$, satisfying the
 idempotency condition (1) above, as well as the following two conditions:
   \begin{enumerate}

  \item[(3)] (Collinear transitivity) If three cells $C, C', C''$ are collinear, then the
  natural transformation $\kappa_{C C' C''}$ is an isomorphism of functors. 
  
    \item[(4)] (Invertibility)  if $C, C'$ are two cells of the same dimension $r$ lying in a linear subspace of
  dimension $r$ and separated by a cell $D\leq C, C'$ of dimension $(r-1)$ (``wall''), then
  $\varphi_{CC'}$ is an equivalence of categories.
  
 \end{enumerate} 
\end{Defi} 
 
 We will use $\Hc$-schobers as  a possible avatar of
 a perverse schober on $\CC^n$ smooth with respect to $\Hc_\CC$. 
 Given an $\Hc$-schober  $\Fen$, forming the rational Grothendieck groups $K(\Ec)$ gives
 a diagram which gives a perverse sheaf $K(\Fen)\in\Perv(\CC^n, \Hc)$. 
 One way to interpret Proposal 
 \ref{prpsl:HMMP} is to  construct $\Hc$-schobers formed by derived categories
 of varieties appearing in a web of flops. We will refer to them as {\em flobers}.

 \vfill\eject
 
 \section{The Grothendieck resolution}\label{sec:groth}
 
 \noindent{\bf A. Definition of the resolution and the $\gen$-web of flops.} 
 Let $\k$ be an algebraically closed field of characteristic $0$. 
Let $\gen$ be a reductive Lie algebra over $\k$,
  corresponding to a reductive algebraic group $G$, 
 and $F\simeq G/B$ be its flag variety. It can be seen as parametrizing all the Borel subalgebras $\ben\subset\gen$, so there
 is a tautological bundle $\pi: \ul\ben\to F$. The {\em Grothendieck resolution}, see, e.g. \cite{bez} \cite{BR} \cite{chriss-ginz}
 is the total space $\wt\gen$ of this bundle, i.e.,
 \[
 \wt\gen \,\,=\,\, \bigl\{ (x,\ben)\in \gen\times F \bigl| x\in\ben\bigr\}. 
 \]
 It comes with a natural projection $g: \wt\gen\to \gen$, a proper map.

  Let $\gen_\rss\subset\gen$ be the open subset of
   regular semisimple elements (``matrices with distinct eigenvalues''). Over $\gen_\rss$, the projection $g$ is
  an unramified Galois covering
 with Galois group $W$ (the Weyl group of $\gen$).

  In fact, $g$ factors into the composition  (Stein factorization)
of a finite morphism $\varpi$ and  a map $f$ with connected fibers (which, in our case, happens to
be birational, since these fibers generically consist of one point):
 \be\label{gen-f}
 \wt\gen \buildrel f\over \lra Z:=\, \gen \times_{\hen/W} \hen \buildrel \varpi\over \lra \gen.
 \ee
 Here $\hen$ is the Cartan subalgebra in $\gen$, and we use the identification $\gen/\on{Ad}(G) = \hen/W$
 to construct the ``characteristic polynomial map''
 \[
 \chi: \gen \lra \gen/\on{Ad}(G) = \hen/W. 
 \]
 The variety $Z$ is singular while $\wt\gen$ is smooth.
 
 In particular,  we can view points of $Z$ as pairs
 $(x, K)$, where $x\in\gen$ and $K$ is a connected component of  $g^{-1}(x)$, the variety of all Borel subalgebras containing $x$. 
 
 \begin{ex}\label{ex:sl2Groth}
 Let $\gen=\sen\len_2$, then $F=\PP^1$, and $\ul\ben\simeq \Oc(-1)^{\oplus 2}$. The variety $Z$ is the 3-dimensional
 quadratic cone (the 2-sheeted covering of $\gen\simeq\AAA^3$ ramified 
 along the quadratic cone in $\AAA^3$).
 So $f$ is the flopping contraction for the Atiyah flop (local model, Example \ref{ex:a-flop-loc}). 
  \end{ex}
  
  \begin{rems}
  (a) The morphism $\wt\gen\to\hen$ implicit in \eqref{gen-f} comes from the following well known but important phenomenon
  \cite{chriss-ginz}. The vector bundle $\ul\ben^{\on{ab}} = \ul\ben/[\ul\ben, \ul\ben]$ of abelianized Borel subalgebras on $F$, is trivial,
  and its space of global sections $\hb = H^0(F,  \ul\ben^{\on{ab}})$ can be seen as the
  {\em universal Cartan subalgebra} of $\gen$. More precisely, for any particular Cartan $\hen\subset\gen$,
  a choice of a Borel $\ben\supset\hen$ identifies $\hen$ with $\hb$; a different choice of $\ben$ through $\hen$
  changes this identification by an element of $W$. So we have a canonical map $\wt\gen\to\hb$, and a more intrinsic
  definition of $Z$ is $Z=\gen\times_{\hb/W}\hb$. 
  
  In the sequel we  fix a distinguished Cartan $\hen\subset\gen $ together with a Borel $\ben_+\supset\hen$ so writing
  $Z$ as in \eqref{gen-f} becomes unambiguous.   This
    choice  fixes $\Delta\supset\Delta_+\supset\ds$, the systems of
   roots, positive roots and simple roots of $\gen$, which we will use whenever needed.

  \vskip .2cm
  
  (b) Since $Z$ is affine and $f$ is proper, $Z = \Spec \, H^0(\wt\gen, \Oc)$ is the affinization of $\wt\gen$.  
  \end{rems}

  The Grothendieck resolution gives rise to a very explicit ``$\gen$-web of flops" constructed as follows.
 Note that the variety $Z$ is acted upon by $W$, so we define the $w$-flopped contraction as
 the base change
 \[
 f_w: X_w := w^*\wt\gen \lra Z, \quad w\in W.
 \]
Although all the $X_w$ are isomorphic as algebraic varieties, they are connected by
 nontrivial birational isomorphisms coming from the birational identifications $f_w$  with $Z$. 
  This is similar to defining a  3d flop by base change along an involution (which is $z\mapsto z^t$ in
  Example \ref{ex:a-flop-loc}). 
  
  \begin{rem}
 According to Brieskorn \cite{brieskorn}, his simultaneous resolution of an ADE singularity is obtained
 from the Grothendieck resolution for  the corresponding $\gen$,  by a base change: restriction to
 the {\em Brieskorn slice}, a transversal slice to the subregular nilpotent orbit. Combining this
 with the elephant picture of 3-dimensional flops, see \S \ref{sec:webflop}C, we can say
 that the Grothendieck resolution is ``the mother of all  flops'' (at least of $3$-dimensional ones). 
  \end{rem}

  \vskip .3cm

  \noindent {\bf B. Singular nature of $Z$.} By Example
  \ref {ex:sl2Groth}, 
  the variety $Z\subset \gen\times\hen$ can be considered as the $\gen$-analog of the
  3-dimensional quadratic cone. This analogy makes it natural to introduce the following varieties
  which we will use later.
  
  \vskip .2cm
   First, let $\pr: \hen\to\hen/W$ be the projection. 
  Inside $\hen$ we have the {\em root hyperplanes} $\hen _\alpha = \alpha^\perp$, $\alpha\in\Delta_+$.
   We denote by 
  \[
  \nabla \,\,=\,\, \pr\bigl( \bigcup_{\alpha\in\Delta_+} \hen_\alpha \bigr)\,\,\subset \,\, \hen/W
  \]
  the {\em discriminantal hypersurface in} $\hen/W$, i.e., the locus of  
  points represented by non-free orbits of  $W$. 
  We further denote by 
 $D=\chi^{-1}(\nabla)\subset \gen$ the preimage of $\nabla$ and call it the
 {\em discriminantal hypersurface in} $\gen$. 
 Thus $\gen - D = \gen_{\rss}$ is the open set of {\em regular semisimple elements}
 (``matrices with distinct eigenvalues'').

   By construction, the hypersurface $D$ is 
    the branch locus of the finite morphism $\varpi: Z\to \gen$. 
  Ramification points of $Z$, i.e., points lying over   $D$, are  generically smooth. 
 
  For $\gen =\sen\len_2$, the variety $D$
 becomes a 2-dimensional quadratic cone in $\sen\len_2=\AAA^3$ (the cone of
 nilpotent elements).
 \vskip .2cm
 
 Singularities of $Z$ are governed by the closed subvariety $\gen_\nr\subset\gen$
 formed by  non-regular elements. We recall several equivalent characterizations of this concept.
 
 \begin{itemize}
 \item[(R)] An element  $x\in\gen$ is called {\em regular},
 if the centralizer $\gen^x$ of $x$ is abelian, and {\em non-regular}  otherwise, see   \cite{chriss-ginz}.
  The sets of
 regular and non-regular elements of $\gen$ are denoted $\gen_\reg$ and $\gen_\nr$. 
 
 \item[(R')] 
 Alternatively, $\gen_\nr$
 consists of the {\em critical points} of $\chi$, i.e., of $x\in \gen$ such that the differential $d_x\chi$
 is not of full rank, see \cite{kostant}, Th. 0.1. 
 
 \item[(R'')]  Any element $x\in\gen$ has the (unique) {\em Jordan} (or {\em Cartan})
 {\em decomposition}
 $x=y+z$ where $y$ is semisimple, $z$ is nilpotent and $[y,z]=0$. In this case $\gen^y$,
 the centralizer of $y$, is reductive, and it is classical (\cite{kostant}, Prop. 0.4) that $x$
 is regular iff $z$ is a principal nilpotent element in $\gen^y$, i.e., lies in the open orbit
 in the nilpotent cone of $\gen^y$. 
 
 \item[(R''')] $x$ is regular if and only if   $F^x = \pi (g^{-1}(x)) \subset F$, the variety 
 of Borel  subalgebras  
 through $x$, is $0$-dimensional, i.e., consists of finitely many points. 
 \end{itemize}

 \begin{prop}\label{prop:nr-sing}
 (a1)  $D_\sing$, the singular locus of $D$, is the union of two irreducible components:
 \[
 D_\sing \,\,=\,\, \gen_\nr \,\cup \, \chi^{-1}(\nabla_\sing). 
 \]
\vskip .2cm 
 
 (a2) The variety 
  $\gen_\nr$  has codimension $3$ in $\gen$ and $2$ in $D$. The variety $\chi^{-1}(\nabla_\sing)$
  is empty for $\gen=\sel_2$, and in all other cases has codimension $2$ in $\gen$
  and codimension $1$ in $D$. 
  
  \vskip .2cm
  
  (a3) 
 The intersection of $D$ with the transversal
 slice at a generic point of $\gen_\nr$ has a 2-dimensional quadratic cone singularity.
 
 \vskip .2cm
 
 (b) The variety $Z_\sing$ coincides with the preimage $\varpi^{-1}(\gen_\nr)$.
 It has codimension $3$ in $Z$. The intersection of $Z$ 
 with the transversal
 slice at a generic point of $\gen_\nr$ has a 3-dimensional quadratic cone singularity.
 
 \end{prop}
 
 \noindent{\sl Proof:} (a1) We first prove that $D_\sing$ is contained in 
 $ \gen_\nr \cup  \chi^{-1}(\nabla_\sing)$, that is, a point $x$ lying in $D\cap \gen_{\reg}$
 and mapped by $\chi$ into a smooth point of $\nabla$, is a smooth point of $D$.
 Indeed, by (R') the map $\chi$ is a smooth map near $x$, so near $x$ the hypersurface
 $D=\chi^{-1}(\nabla)$ is smooth. 
 
 To see the reverse inclusion, $ \gen_\nr \cup  \chi^{-1}(\nabla_\sing)\subset D_\sing$,
 we first note that, by the same argument using (R'), the
   the intersection
$\gen_\reg \cap \chi^{-1}(\nabla_\sing)$ is contained in $D_\sing$. 
So it remains to prove that $\gen_\nr\subset D_\sing$. This will follow from the
part (a3)  below. 

\vskip .2cm

(a2) The statement about codimension of $\gen_\nr$ is well known. The variety $\nabla_\sing$
is the image, under $\pr:\hen\to \hen/W$, of the locus of points lying on more than one
root hyperplane. So it is empty for $\gen=\sel_2$ and has codimension $2$ in $\hen$
in all other cases. Our statement follows from this by applying $\chi^{-1}$.
 
 \vskip .2cm
 
 (a3) Using  the Jordan decomposition in (R''), we see that the generic non-regular case
 is   when $\gen^y\simeq \sen\len_2 \oplus (\text{abelian})$ (``simple coincidence of eigenvalues") and $z=0$. 
 Denote by $\gen_\nr^\circ$ the open part of $\gen_\nr$ formed by such semisimple $x=y$.
 For $y\in\gen_\nr^\circ$,  a transversal slice to  $\gen_\nr$ at $y$ can be taken to consist of
 $y+t$ where $t$ lies in the $\sen\len_2$-part of  $\gen^y$. Being identified with $\sen\len_2$, it is 
 3-dimensional, and its
 intersection with $D$ is the cone of nilpotent elements in $\sen\len_2$. 
 
 \vskip .2cm
 
 (b) First, let us prove that $\varpi^{-1}(\gen_\reg)$ consists of smooth points of $Z$. Indeed,
 let $x\in\gen_\reg$. Then  by (R') there is a neighborhood $U$ of $x$ in $\gen$ such that
 $\chi: U\to \hen/W$ is a smooth morphism. Therefore $\varpi^{-1}(U)$ is the pullback, under this
 smooth morphism, of $\hen$ mapping to $\hen/W$. Since $\hen$ is smooth, so is the pullback. 
 
 Next, let us prove that any point of $Z$ lying over any  $x\in \gen_\nr$ is singular. For this it is enough to
 assume that $x\in\gen_\nr^\circ$ (``generic'' case) and to prove, in this case, the last statement of (b). 
 Assuming this, we see that over such $x$ 
  the ramification of $\varpi$ is simple quadratic, as any preimage $z\in Z$
 of $x$ has the stabilizer in $W$ identified with $\ZZ/2$. Together with (a3), this means that near
 such $z$, the variety $Z$ looks like the product of a smooth manifold and a double covering of $\AAA^3$
 ramified along a quadratic cone. \qed

 \vskip .2cm
 
 Denote $S=Z_\sing$. We now analyze $S$ in terms of the projection $\kappa: Z\to\hen$.
  
 \begin{prop}\label{prop:components}
 Irreducible components of $S$ are in bijection with positive roots $\alpha\in\Delta_+$. The component
 $S_\alpha$ corresponding to $\alpha$, projects surjectively to the root hyperplane $\hen_\alpha\subset \hen$. 
 \end{prop}
 
 \noindent {\sl Proof:} Inside $S=\varpi^{-1}(\gen_\nr)$ we have the open dense subset
 $S^\circ = \varpi^{-1}(\gen_\nr^\circ)$, where $\gen_\nr^\circ$ is, as in the proof of
 Proposition \ref{prop:nr-sing}(a), the variety of semi-simple elements $y$ with 
 $\gen^y\simeq \sen\len_2\oplus (\text{abelian})$. 
 An element $y\in \hen$ lies in $\gen_\nr^\circ$ iff it belongs to $\hen^\circ_\alpha = \hen_\alpha -\bigcup_{\beta
 \neq\alpha} \hen_\beta$ for some $\alpha\in\Delta_+$. This means that $\gen_\nr^\circ$ is
 the (not necessarily disjoint) union of ``conjugation sweeps'' $\on{Ad}(G)\cdot \hen^\circ_\alpha$.
 Accordingly, consider the action of $\on{Ad}(G)$ on $\gen\times\hen$ through the first factor,
 and let 
 \[
 \on{diag}(\hen^\circ_\alpha)\,\subset\,  \hen\times\hen \, \subset \, \gen\times\hen
 \]
 be the image of $\hen_\alpha^\circ$ under the diagonal embedding. We conclude that
  $S^\circ$ is the (disjoint!) 
 union of the $S^\circ_\alpha = \on{Ad}(G)\cdot \on{diag}(\hen^\circ_\alpha)$.
 Denoting by $S_\alpha$ the closure of $S^\circ_\alpha$, we get the statement. \qed
 
 \vskip .3cm

  \begin{figure}[H]
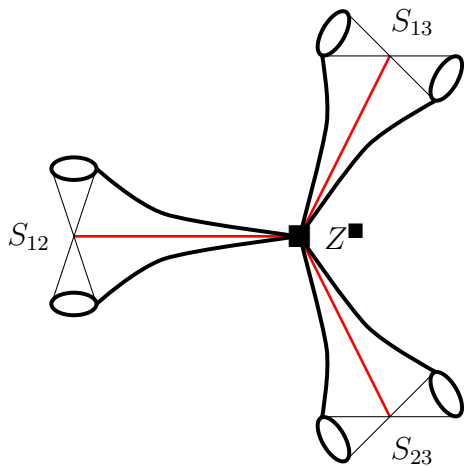

  \centering
 
 \btp[scale=.3]
  
\node (0) at (0,0){$\blacksquare$};
 
 \node at (2,0){$Z^\blacksquare$};    
    \draw[red, line width = 1] (0,0) -- (-10,0); 
    \draw [red, line width=1] (0,0) -- (4,8);
    \draw [red, line width =1] (0,0) -- (4, -8);

  \draw (1,8) -- (7,8);
  \draw (2,10) -- (6,6);
  
  \draw (-11,-3) -- (-9,3);
  \draw (-11,3) -- (-9,-3); 
  
  \draw (2, -10) -- (6, -6);
  \draw (1, -8) -- (7, -8);

 \draw [line width = 1.5]  (-10,3) ellipse (1 and .5);
 \draw [line width = 1.5]  (-10,-3) ellipse (1 and .5);
 
 \draw [line width = 1.5]  (1.5, 9) circle [x radius = 1.1, y radius = .5, rotate =60]; 
 \draw[line width = 1.5]   (6.5, 7) circle [x radius = 1.1, y radius = .5, rotate =60]; 
 
 \draw[line width = 1.5]  (1.5, -9) circle [x radius = 1.1, y radius = .5, rotate =-60]; 
 \draw [line width = 1.5] (6.5, -7) circle [x radius = 1.1, y radius = .5, rotate =-60]; 
 
 \draw[line width = 1.5]  plot [smooth] coordinates  { (0,0) (1.2, 5) (1,8)} ; 
  \draw[line width = 1.5]  plot [smooth] coordinates  { (0,0) (3,4) (6,6)} ; 
  
 \draw[line width = 1.5]  plot [smooth] coordinates  { (0,0) (-6,1) (-9,3)} ;   
 \draw[line width = 1.5]  plot [smooth] coordinates  { (0,0) (-6,-1) (-9,-3)} ;   
 
  \draw[line width = 1.5]  plot [smooth] coordinates  { (0,0) (3, -4) (6, -6)} ;  
  \draw[line width = 1.5]  plot [smooth] coordinates  { (0,0) (1.2, -5) (1,-8)} ;  
  
  \node at (-12,0){$S_{12}$};   
 \node at (5,9.5) {$S_{13}$};  
  \node at (5,-9.5) {$S_{23}$};

 \etp
  \caption{Structure of $Z$ up to codimension 4. Red lines represent components $S_\alpha$ of $Z_\sing$,
  with labels given for $\gen=\sel_3$.}
  
  \end{figure}

  \begin{rem}
  Several important properties of the  variety $Z$ have been established in  \cite{hotta-kash}
  using the 
  morphism $f: \wt\gen\to Z$, which is a small resolution of  singularities.    \end{rem}

 \noindent {\bf C. The $\gen$-web of flops up to codimension $4$.}
 Let 
 $
 \gen_\nr^\blacksquare\,\,=\,\,\gen_\nr - \gen_\nr^\circ
 $
 be the closed subset in $\gen_\nr$ formed by non-regular elements more complicated 
 (i.e, other) than semisimple elements $y$ with $\gen^y\simeq\sen\len_2\oplus (\text{abelian})$. 
  Using the Jordan decomposition,
 it is easy to see that $\gen_\nr^\blacksquare$ has codimension $1$ in $\gen_\nr$ and
 codimension $4$ in $\gen$. Denote
 \[
 \gen^\Box \,\,=\,\,\gen - \gen_\nr^\blacksquare
 \]
 be the complement of this codimension $4$ subvariety. Let also 
 \[
 Z^\blacksquare = \varpi^{-1}(\gen^\blacksquare), \quad Z^\Box = Z-Z^\blacksquare = \varpi^{-1}(\gen^\Box)
 , \quad 
 \wt\gen^\Box =  f^{-1}(Z^\Box).
 \]
  The open subvariety $Z^\Box$ is invariant under the $W$-action on $Z$,
 so we have the restriction of the $\gen$-web of flops to it:
 \[
 f_w^\Box:  X_w^\Box = w^* \wt\gen^\Box \lra Z^\Box. 
 \]
 This web of flops can be described directly in terms of the standard Atiyah flops. Indeed, the singular locus
 of $Z^\Box$ is the disjoint union of  smooth varieties  $S^\circ_\alpha$, $\alpha\in\Delta_+$,
 and transversely to each $S_\alpha$,
 the variety $Z^\Box$ looks like a 3-dimensional quadratic cone. 
 Therefore, a neighborhood of $S_\alpha\subset Z^\Box$ 
 has two small desingularizations corresponding to the two desingularizations of the 3d quadratic cone. 
 One of them is given by  the restriction of $f^\Box: \wt\gen^\Box\to Z^\Box$. Let us call it the
 {\em positive desingularization} of $Z^\Box$ along $S_\alpha^\circ$, and the other one will be called
 the {\em negative} one. Thus we have $2^{|\Delta_+|}$ possible small desingularizations of $Z^\Box$. 
 Among these, we have $|W|$ desingularizations $X_w^\Box$. They are described as follows.
 
 \begin{prop}\label{prop:X_w-desingul}
 Let $w\in W$ and $\alpha\in\Delta_+$. Then $f_w^\Box: X_w^\Box\to Z^\Box$ is, along $S_\alpha^\circ$:
 \[
 \begin{cases}
 \text{a positive desingularization, if } & w(\alpha)\in\Delta_+;
 \\
 \text{a negative desingularization, if } & w(\alpha)\in \Delta_-  . 
 
 \end{cases}
 \]
 \end{prop}
 
 \noindent{\sl Proof:} 
 We analyze the sets naturally labeling the combinatorial objects we want to compare. Let $\Irr(S)$ be the set of
 irreducible components of $S$. It has a $W$-action coming from the $W$-action on $Z$. 
 From the proof of Proposition \ref{prop:components} we see that, as a $W$-set, $\Irr(S)$ is identified with
 the set of root hyperplanes $\hen_\alpha$, the set  $ \Delta/\pm$ of pairs
 $[\alpha] := \{\alpha, -\alpha\}$ of opposite
roots.  In the remainder of this proof we will use the more precise  notation $S_{[\alpha]}$ for components of $S$. 
 
 \vskip .2cm
 
 Let further $\wt\Irr(S)$ be the set of pairs consisting of a component of $S$ and a choice of one of two
 of small desingularizations of $Z$ near that component. This set is also equipped with a $W$-action and we have
   a $W$-equivariant surjective map $\wt\Irr(S)\to \Irr(S)$ whose fibers have cardinality $2$. 
   
   \begin{lem}\label{lem:desing}
   We have a commutative diagram with vertical arrows being isomorphisms of $W$-sets:
   \[
   \xymatrix{
   \wt\Irr(S)\ar[r]\ar[d]& \Irr(S)
   \ar[d]
   \\
   \Delta\ar[r]&\Delta/\pm. 
   }
   \]
   \end{lem}
 
 The lemma implies Proposition \ref{prop:X_w-desingul}. Indeed,   $\Delta_+\subset\Delta$
 is a fundamental domain for the action of $\{\pm 1\}$ on $\Delta$. So if we call a positive desingularization 
 of $S_{[\alpha]}$ the one corresponding to the positive root $\alpha$ inside $[\alpha]$, then the action of
 $w\in W$ will send it to a positive or negative desingulariation according to whether $w(\alpha)$ is
 a positive or negative root. 
 
 \vskip .2cm
 
 \noindent{\sl Proof of Lemma \ref{lem:desing}:} 
  We first consider 
   the two small desingularizations $\wt Q_\pm\to Q$ of a 3-dimensional quadratic
 cone $Q\subset \AAA^4$. They are
   naturally labelled by the two families of planes 
 (2-dimensional linear subspaces)
 lying on $Q$, i.e., by the two families of generators of the 2-dimensional quadric 
 $\PP(Q) \simeq\PP^1\times\PP^1$. Call them the $(+)$- and $(-)$-families. Planes
 from the $(\pm)$-family lift bijectively into $\wt Q_\pm$ while planes from the $(\mp)$-family
 lift to $\wt Q_\pm$ as blowups. 
 
 \vskip .1cm
 
 We now look at the  two families of 2-planes in the transverse slices to the $S_{[\alpha]}$ and show  that
 they are naturally labelled by the roots (positive or negative)  $\alpha\in\Delta$ themselves inside the
 pairs $[\alpha]$ of opposites. For this it is sufficient to take  the slice at a point $(x,K) \in S_{[\alpha]}$ 
 such that  $x$ is a generic point of 
  the root hyperplane $\hen_\alpha\subset\hen\subset \gen$ in the chosen Cartan and 
  $K$ is the component of the variety of Borels through $x$
 such that  the standard Borel $\ben_+$ belongs to $K$. 
  This slice is provided
 by the variety $Z({\sen\len_2[\alpha]})$ corresponding to the
 $\sen\len_2$-subalgebra $\sen\len_2[\alpha]\subset\gen$ with the set of weights $[\alpha]$. 
 The two families of $2$-planes in $Z({\sen\len_2[\alpha]})$ are naturally labelled by the two Borels
 of $\sel_2[\alpha]$ through the chosen Cartan, i.e., by the two roots in $[\alpha]$. Indeed, such a Borel
 is a 2-dimensional subspace and its lift to $Z({\sen\len_2[\alpha]})$ is a $2$-plane representing
 the corresponding family. \qed

  \vfill\eject
  
  \section{The $\gen$-web of flops and partial blowdowns.}\label{sec:blowdowns}

  \noindent{\bf A. The  movable cone picture and partial blowdowns.} 
    \begin{wrapfigure}{r}{0.3\textwidth}
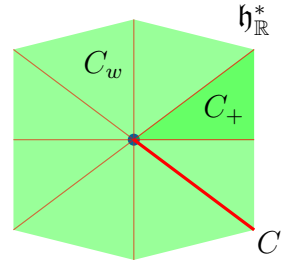


  \btp[scale=.4]
  
  \node (0) at (0,0){};
  \fill[color=blue] (0) circle (0.2); 
  
  \draw[  draw=red, fill=green, opacity = 0.4] (0,4) -- (0,0) -- (4,3);
   \draw[  draw=red, fill=green, opacity = 0.6] (4,3) -- (0,0) -- (4,0);
  \draw[ draw=red, fill=green, opacity = 0.4] (4,0) -- (0,0) -- (4,-3);
   \draw[ draw=red, fill=green, opacity = 0.4] (4,-3) -- (0,0) -- (0, -4);
  \draw[ draw=red, fill=green, opacity = 0.4] (0,-4) -- (0,0) -- (-4,-3);
  \draw[ draw=red, fill=green, opacity = 0.4] (-4,-3) -- (0,0) -- (-4,0);
 \draw[draw=red, fill=green, opacity = 0.4] (-4,0) -- (0,0) -- (-4,3); 
  \draw[ draw=red, fill=green, opacity = 0.4] (-4,3) -- (0,0) -- (0,4);
  
  \draw[red, line width = 1.2] (0,0) -- (4,-3){};
  
  \node at (4,4) {$\hen^*_\RR$};
  
  \node at (3,1){$C_+$}; 
  \node at (-1,2.5){$C_w$}; 
  \node at (4.5,-3.4){$C$}; 
  
  \etp
  \caption{Cells of the coroot arrangement.}

    \end{wrapfigure}
 We now explain the analog, for the Grothendieck resolution,
 of  the relation between  flops and chambers of  the 
   movable cone. This relation is particularly transparent in our 
    ``local''  situation, as all divisors are movable.
    We keep the notation of the previous section.

 The space $\hen_\RR^*$  of real
  weights carries the   {\em coroot arrangement} $\Hc=\Hc_\Delta$ of hyperplanes 
  $H_\alpha$,  $\alpha\in\Delta_+$. 
   Explicitly, $H_\alpha$ is the orthogonal complement of the coroot $\alpha^\vee$
   corresponding to $\alpha$.   
   Therefore we have the decomposition of $\hen^*_\RR$ into cells of $\Hc$. The chambers
   (open cells) $C_w$  are labelled by $w\in W$ via $C_w = w^{-1}(C_+)$, where $C_+$ is the
   {\em dominant chamber}, on which all $\alpha^\vee, \alpha\in\Delta_+$ are positive. 
   
   \vskip .2cm
   
   This decomposition is  the $\gen$-analog of the  chamber structure of the movable
   cone for 3d flops. That is,  denote $X=\wt\gen$. 
 We have 
 \[
 \Pic(X) = \Pic(G/B) = \hen^*_\ZZ,
 \]
  the lattice of  integer weights, the isomorphism given
 by pullback along $\pi$. For $\lambda\in \hen^*_\ZZ$
 we denote by $\Lc(\lambda)=\pi^*\Oc(\lambda)$ the corresponding line bundle. 
 Then all $\Lc(\lambda)$
 are movable. We denote
 \[
 Z_\lambda \,\,=\,\,\Proj_Z \biggl( \bigoplus_{n=0}^\infty f_*\,\,  \Lc(n\lambda) \biggr). 
 \]
 This is the local analog of ``the image of the rational map to a projective space defined
 by the linear system $|\Lc(\lambda)|$''. 
 
 Thus, if $\lambda \in C_+$, then $Z_\lambda=X$, since $\Lc(\lambda)$
 is ample. If $\lambda=0$, then $Z_\lambda = Z$ is the affinization of $X$. 
 
 \begin{prop}\label{prop:bdowns}
 (a) $Z_\lambda$ depends only on the cell $C$ of $\Hc$ containing $\lambda$.
 So we will use the notation $Z_C$ for it. 
 
 \vskip .2cm
 
 (b) For any two cells $C\leq C'$ we have a regular birational morphism $q_{C',C}: Z_{C'}\to Z_C$.
 These morphisms are transitive whenever $C\leq C'\leq C''$. 
 
 \vskip .2cm
 
 (c) For any $w\in W$ we have $Z_{w^{-1}(C)} = w^* Z_C$ as a variety
 over $Z$.

 \end{prop}

 In particular, $Z_{C_w} = X_w$ is the flopped Grothendieck resolution. The diagram formed by the
 $Z_C$ and $q_{C'C}$ is the $\gen$-analog of the diagram \eqref{eq:a-flop} defining the Atiyah flop. 
 
 \paragraph{B. $W$-equivariance of the $Z_\lambda$.} To prepare for the proof of Proposition
 \ref{prop:bdowns}, we start with the following statement.
 
 \begin{prop}\label{prop:Z-lambda-eq}
 For any $w\in W$ and $\lambda\in\hen^*_\ZZ$ we have $Z_{w^{-1}(\lambda)} = w^* Z_\lambda$ as
 a variety over $Z$. 
  \end{prop}
  
  \noindent {\sl Proof:} 
   The smooth locus of $Z$ is, 
  by Proposition \ref{prop:nr-sing}, equal to $Z_\reg = \varpi^{-1}(\gen_\reg)$. 
   It is
  a $W$-invariant open subset whose complement $Z_\sing$ has, by Proposition \ref{prop:nr-sing},
  codimension $3$ in $Z$. Let $ \wt\gen_\reg= f^{-1}(Z_\reg)\subset\wt\gen$. 
  The Grothendieck resolution map $f: \wt\gen_\reg \to Z_\reg$
  is an isomorphism.  
   The complement $\wt\gen- \wt\gen_\reg$ has, by the analysis of \S \ref{sec:groth}C, codimension 2 in
  $\wt\gen$.
   Over the generic part of $S=Z_\sing$ given by the disjoint union of the $S_\alpha^\circ$,
  it is a $\PP^1$-bundle (coming from the  small desingularizations of the transverse 3d quadratic cones). 
  This means that each vector bundle $\Mc$ on $\wt\gen_\reg$ extends canonically to a torsion free
  sheaf $j_*\Mc$ on $\wt\gen$, where $j: \wt\gen_\reg \hookrightarrow \wt\gen$  is the embedding. 
  In particular, if $\Mc$ was already the restriction of a vector bundle $\Nc$ from $\wt\gen$ to $f^{-1}(Z_\smth)$,
  the $j_*\Mc = \Nc$.

  Because of the identification with $Z_\reg$, the group $W$ acts by automorphisms of the open
  subvariety $ \wt\gen_\reg \subset \wt\gen$.
   Proposition \ref {prop:Z-lambda-eq} now  reduces to the following  fact.

   \begin{lem}\label{lem:w-lambda}
   Let $\lambda\in\hen^*_\ZZ$ and $w\in W$. Then
   \[
   w^*\bigl( \Lc(\lambda)|_{\wt\gen_\reg}\bigr) \,\, \simeq\,\,  \Lc(w^{-1}\lambda)|_{\wt\gen_\reg}. 
   \]
   \end{lem}
   
   Indeed, the lemma implies that for any $n\geq 0$, the coherent sheaf
   $f_*\Lc(w^{-1}(n\lambda))$ on $Z$ is identified with $w^* (f_*\Lc(n\lambda))$,
   and so $Z_{w^{-1}\lambda}\simeq w^* Z_\lambda$ as a scheme over $Z$.

   \vskip .2cm
   
   \noindent{\sl Proof of Lemma \ref{lem:w-lambda}:} This is well known. 
     We are grateful to
   R. Bezrukavnikov for the explanation regarding this fact and for suggesting the proof below. 
   
   \vskip .2cm
   
   Let $\Lc$ and $\Rc$ be  the  line bundles in the left and right hand side of the isomorphism
   claimed in the proposition. They are naturally $G$-equivariant bundles, so it suffices to
    establish the equality
   $[\Lc]=[\Rc]$ of their classes in the $G$-equivariant Picard group $\Pic^G(\gen_\reg)$.
    Because $\wt\gen-\wt\gen_\reg$
   has codimension $2$, the restriction map $\Pic^G(\wt\gen)\to\Pic^G(\wt\gen_\reg)$
   is an isomorphism. On the other hand, the projection $\wt\gen\to F$ identifies $\Pic^G(\wt\gen)$ with
   $\Pic^G(F)= \hen^*_\ZZ$. 
   
   Inside $\wt\gen_\reg$,  consider the locally closed subvariety $\wt O= (\varpi f)^{-1}(O)$,
   where $O$ is the adjoint orbit of a regular (semisimple) element $x\in\hen$. Then 
   $\wt O = \bigsqcup_{y\in W}  O_y$ is the disjoint union of   isomorphic lifts of $O$ labelled by $y\in W$. 
   For each $y$ we have that $\Pic^G(O_y)=\hen^*_\ZZ$ and so the restriction
   $\Pic^G(\wt\gen_\reg)\to\Pic^G(O_y)$ is an isomorphism. It is therefore enough to show that
   for some (or, what is equivalent, for any) $y$
   the classes, in $\Pic^G(O_y)$, of the  restrictions of $\Lc$ and $\Rc$, are equal. 
   
   To see this, we note that the action of $w$ on $\wt\gen_\reg$  permutes the $O_y$, and
   the map from $\Pic^G(O_1)\simeq \hen^*_\ZZ$ to $\Pic^G(O_w)\simeq \hen^*_\ZZ$
   induced by $w$, coincides with the action of $w$ on $\hen^*_\ZZ$. \qed

 \paragraph {\bf C. Blowdowns via parabolic Grothendieck resolutions.} 
 For the proof  of Proposition 
 \ref{prop:bdowns}, and for future use, we give an independent construction of the varieties $Z_C$
 satisfying parts (b) and (c), and then identify $Z_C$ with $Z_\lambda$ for any $\lambda\in C$. 
 We start with the following standard fact. 
 
 \begin{prop}\label{prop:cells-parab}
 (a) The poset $\Cc$ of cells in $\Hc$ is anti-isomorphic to the the poset of parabolic subalgebras in $\gen$
 containing $\hen$. Explicitly, the parabolic subalgebra $\pen_C$ corresponding to $C$ is spanned by $\hen$
 and  the root vectors $e_\alpha$ for all $\alpha\in\Delta$ such that $\alpha|_C \geq 0$. 
 
 \vskip .2cm
 
 (b) Under the above identification, cells $C$ lying in the closure of the dominant chamber $C_+$
 correspond to parabolic subalgebras containing the standard Borel $\ben_+$. \qed
 \end{prop}
 
 For $C\in\Cc$ let $P_C\subset G$ be the parabolic subgroup corresponding to $\pen_C$
 and $F_C=G/P_C$ be the variety of conjugates of $\pen_C$. We then have the {\em parabolic Grothendieck
 resolution}
 \[
 \wt \gen_C  \,\,=\,\, \bigl\{ (x,\pen)\in \gen\times F_C \bigl| x\in\pen\bigr\} \buildrel \pi_C\over\lra F_C, \quad
 \pi_C^{-1}([\pen])=\pen. 
 \]
  Here $[\pen]\in F_C$ is the point corresponding to a parabolic $\pen$ conjugate to $\pen_C$.
  For any such $\pen$ let $\nen_\pen$ be its nilpotent radical,  and $\len_\pen = \pen/\nen_\pen$ 
   the Levi quotient. 
 Thus $\len_\pen$ is a reductive Lie algebra, and we can form its own Grothendieck resolution
  $\wt\len_\pen\to\len_\pen$ and its affinization
  $Z(\len_\pen)=\Spec \,\, H^0(\wt\len_\pen, \Oc)$.
  Denote $Z(\pen)$ the fiber product, which induces a finite covering of $\pen$,
 \[
 Z(\pen) \,\,:=\,\, \pen\times _{\len_\pen} Z(\len_\pen) \lra \pen.
 \]
 Such covering can be formed for any parabolic subalgebra $\pen$ conjugate to $\pen_C$, and we denote by
 $Z_C$  the universal family of these coverings: 
 \be\label{eq:Z=universal}
 Z_C\,\, := \,\, Z(\ul\pen)\,=\, G\times_{P_C} Z(\pen_C)  \lra F_C. 
 \ee
 Thus $Z_0=Z$ while $Z_{C_+}=\wt \gen$. 
 
 \vskip .2cm
 
 If $C\leq C'$, then $\pen_{C'}\subset\pen_C$, so $P_{C'}\subset P_C$, and we have the projection
 $\rho_{CC'}: F_{C'}\to F_C$. Considering $\wt\gen_C$ and $\wt\gen_{C'}$ as universal bundles of
 parabolic subgroups over $F_C$ and $F_{C'}$ respectively, we have the $G$-equivariant map
 $\alpha_{C'C}: \wt\gen_{C'}\to \rho_{CC'}^* \wt\gen_{C}$, where $\rho_{CC'}^*$ means the
 induced (pullback) bundle. The map $\alpha_{C'C}$ is normalized uniquely by the requirement that
 over the base point $[\pen_{C'}]$  of $F_{C'}$, it is equal to the embedding $\pen_{C'}\to\pen_C$. The map 
 $\alpha_{C'C}$ gives rise to a morphism $q_{C'C}: Z_{C'}\to Z_C$ satisfying part (b) of  Proposition 
 \ref{prop:bdowns}. 
 Part (c)  ($W$-equivariance of the $Z_C$) also follows directly from the construction. 
 
 Proposition 
 \ref{prop:bdowns} now reduces to the following fact. 
 
  \begin{prop}\label{prop:Zl=ZI}
 If $\lambda\in\hen^*_\ZZ$ lies in a cell $C\in \Cc$, then $Z_\lambda$ is identified with $Z_C$
 as a variety over $Z$. 
 \end{prop}
 
 \paragraph {\bf D. Proof of Propositions   \ref{prop:Zl=ZI}.} 
 Because of $W$-equivariance of the $Z_C$ and $Z_\lambda$, 
   it is enough to consider the case when $C\leq C_+$
 lies in the closure of the dominant chamber. 
 As is well known, such faces
 are labelled by subsets $I\subset \ds$. 
 The face $C_I \leq  C_+$ is given by the conditions 
 \[
 (\alpha^\vee, c) = 0, \alpha\in I, \quad  (\alpha^\vee, c)>0, \alpha\in \ds -  I. 
 \]
 We denote the  corresponding  standard parabolic subalgebra $\pen_I = \pen_{C_I}\supset \ben_+$,
 the parabolic Grothendieck resoluton associated to $C_I$ by $\wt\gen_I$ and set $Z_I=Z_{C_I}$. 
   Let 
 $\rho_I: F\to F_I$ be  the canonical projection. Intrinsically, 
 any Borel $\ben$ is contained in a unique $\pen$ conjugate to $\pen_I$. 
This gives rise
 to the morphism $g_I: \wt\gen\to \wt\gen_I$.
 
 \begin{prop}
 The variety $Z_I$ fits into a Stein factorization 
  \[
 \xymatrix{
 \wt \gen \ar[r]^{f_I}
 \ar[dr]_{g_I} & Z_I 
 \ar[d]^{\varpi_I} 
 \\
 & \wt \gen_I
 }
 \]
 with $\varpi_I$ being finite and $f_I$ having connected fibers (and, in our case, being
 birational). 
 \end{prop}
 
 \noindent{\sl Proof:}    The proposition means that   $\k$-points of $Z_I$ are identified with triples $((x,\pen),K')$
 where $(x,\pen)\in \wt\gen_I$  and $K'$ is a connected component of  $g_I^{-1}(x,\pen)$, 
 the variety of Borels $\ben$ such that $x\in\ben\subset\pen$. Now, such Borels are in bijection
 with Borel subalgebras in $\len_\pen$ containing $\ol x$, the image of $x$ in $\len_\pen$. So  $K'$ corresponds 
 (in a bijective
 fashion) to a connected
 component $K''$ of the varieties of Borels in $\len_\pen$ through $\ol x$. Now, the set of such components $K''$
 is,  by definition,
 the fiber of $Z(\pen)\to\pen$ over $x\in\pen$. \qed
 
 \vskip .2cm

  We now prove Proposition \ref{prop:Zl=ZI}. Consider the commutative diagram
 \[
 \xymatrix{
 Z& \ar[l]_f \wt\gen  = \ul\ben
 \ar[d]^{g_I} 
 \ar[dl]_{f_I}
  \ar[r]^\pi & F = \{\ben\}
 \ar[d]^{\phi_I}
 \\
 Z_I \ar[u]^{q_I}
 \ar[r]_{\varpi_I}& \wt\gen_I = \ul\pen \ar[r]_{\pi_I}& F_I = \{\pen\}.
 }
 \]
If $\lambda\in C_I$, then the line bundle $\Oc(\lambda)$ on $F$ is the pullback $\Oc(\lambda)
= \phi_I^* \Oc_{F_I}(\lambda)$ of a very ample line bundle $\Oc_{F_I}(\lambda)$ on $F_I$.  
Let $\Lc_I(\lambda) = \pi_I^*\Oc_{F_I}(\lambda)$, a line bundle on $\wt\gen_I$. Note that
$\varpi_I^* \Lc_I(\lambda)$ is very ample on each fiber $q^{-1}_I(z)$, $z\in Z_I$. Indeed, any  $q_I^{-1}(z)$
is identified with a closed subvariety of $F_I$, the variety of all parabolics of type $I$,
and the restriction of $\varpi_I^* \Lc_I(\lambda)$ is identified with the restriction of $\Oc_{F_I}(\lambda)$
from $F_I$. Now, we identify the summand $f_*\,  \Lc(n\lambda)$ in the definition of $Z_\lambda$:
\[
f_* \,\Lc(n\lambda) \,\, \simeq \,\,q_{I*}\, f_{I*}\, g_I^* \, \Lc_I(n\lambda) \,\, \simeq \,\, 
q_{I*} \varpi_I^* \Lc_I(n\lambda),
\]
  where the second isomorphism comes from the fact that $f_I$ is proper with connected fibers.
  We also use that $f_{I*} f_I^* = \Id$. 
  Now, since  $\varpi_I^* \Lc_I(\lambda)$ is very ample on fibers of $q_I$, we have that
  \[
  Z_\lambda \,\,=\,\, \Proj_Z \biggl( \bigoplus_{n=0}^\infty  f_*\,\,  \Lc(n\lambda) \biggr)
  \,\,\simeq \,\, 
  \Proj_Z \biggl( \bigoplus_{n=0}^\infty q_{I*} \,\,  (\varpi_I^* \Lc_I (\lambda))^{\otimes n}  \biggr)
  \,\, \simeq \,\, Z_I. 
  \]
  Proposition  \ref{prop:Zl=ZI}  and therefore Proposition  \ref{prop:bdowns} are proved.

 \paragraph { E. Blowdowns up to codimenion 4.} We now continue the analysis of
    \S \ref{sec:groth}C from which we retain the notation.
     For a cell $C$ of $\Hc$ let $Z_C^\Box = q_C^{-1}(Z^\Box)$. 
    Recall that the singular locus of $Z^\Box$ is the disjoint union of smooth components $S^\circ_\alpha$,
    $\alpha\in \Delta_+$, and  the transverse slice to $S^\circ_\alpha$ in $Z^\Box$ is a 3d quadratic
    cone with two possible small desingularizations, positive and negative. The morphism $q_C: Z_C^\Box\to Z^\Box$
    is therefore a partial desingularization of $Z^\Box$. More precisely:
    
    \begin{prop}
    Let $\alpha\in\Delta_+$ and $C$ be a cell of $\Hc$. Then, in the variety $Z_C$:
    \begin{itemize}
    \item The components $S_\alpha$ with $\alpha^\vee|_C=0$ remain singular (nothing happens).
    
    \item The components $S_\alpha$ with  $\alpha^\vee|_C>0$ are resolved in the positive way.
    
     \item The components $S_\alpha$ with  $\alpha^\vee|_C<0$ are resolved in the negative way. 
    \end{itemize}
    \end{prop}
    
    \noindent{\sl Proof:} 
 As with   Proposition \ref{prop:X_w-desingul}, this follows by looking at a point
 $\wt x = (x, K)\in S_\alpha$ 
 with $x\in\hen_\alpha$ generic and $K$ being the component of Borels through $x$ such that  $\ben_+\in K$.
 Since the 
  the transversal slice to $S_\alpha$ at $\wt x$ is identified with $Z({\sen\len_2[\alpha]})$,
  our statement reduces to the case of $\sen\len_2$. \qed

  \vfill\eject

  \section{Fiber products and varieties of simplices}

  \paragraph{ A. Fiber products and horizontal components.}
  The system of partial blowdowns, constructed in \S \ref{sec:blowdowns}, is the $\gen$-analog
  of the diagram \eqref{eq:a-flop}, describing the Atiyah flop. However, the derived
  equivalence associated to the flop, is constructed using the diagram
  \eqref{eq:a-fib-prod} featuring the fiber product of the two small desingularizarions. 
  It is  now our goal  to construct such a diagram in the $\gen$-situation. 
  
  \vskip .2cm
  
  As before, we denote by $\Hc$ the coroot arrangement in $\hen_\RR^*$.
  For any cell $C$ of $\Hc$ we define the {\em $C$-incidence variety}
  \[
  IX_C \,\,=\,\, \varprojlim_{C'\geq C} Z_{C'} \,\,\subset \,\, \prod_{C'\geq C} Z_{C'}. 
  \]
 This variety can be reducible. At the same time, all $Z_{C'}$ project to $Z=Z_0$
 in a compatible, birational way, in fact biregularly over
 \[
 Z_{\on{rss}} \,\,=\,\, \varpi^{-1}\gen_{\on{rss}} \,\,\subset Z.  
 \]
  Therefore, $IX_C$ has the {\em horizontal component}  $X_C$
 defined
 as the closure of the image of  $Z_\rss$
  under the system
 of the rational maps inverse to the  $q_{C'}: Z_{C'}\to Z$. 
 Recall that the chambers (open cells) $C_w$ of $\Hc$ are labelled by $w\in W$,
 and $q_{C_w}=f_w: Z_{C_w}=X_w \to Z$  is the flopped Grothendieck resolution.
 
 \begin{prop}
 $X_C$ is identified with the closure of the image of
 \[
 (f_w^{-1}): Z_\rss \lra \prod_{w: \,\, C_w\geq C} X_w.
 \]
 \end{prop}
 
  \noindent{\sl Proof:} We show that the projection $\prod_{C'\geq C} Z_{C'} \to \prod_{C_w\geq C}
  Z_{C_w}$ is an isomorphism of $X_C$ onto its image, which we temporarily denote by $J$.
   Indeed, the inverse map is given
  by sending a point $(x_w)_{C_w\geq C} \in J$ to the system $(x_{C'})_{C'\geq C}\in X_C$
  where $x_{C'}$ is defined as $q_{C_w, C'}(x_w)$ for any $C_w\geq C'$.  Here we note
   that for any $C'\geq C$,  the value $q_{C_w, C'}(x_w)$ is, for any $(x_w)\in J$, independent of the choice of
  $C_w\geq C'$, because it is so for $(x_w)$ lying in the image of the map $(f_w)^{-1}$ on $Z_\rss$. 
  
  \qed
  
 The following is clear by construction. 
 
 \begin{prop}
  The $X_C$ form an representation of the poset $\Cc_\Hc$ of faces of $\Hc$, 
  by algebraic varieties.
 More precisely,  for any $C\leq C'$ we have a regular, proper, birational map
 $p_{C, C'}: X_{C}\to X_{C'}$, and these maps are transitive: for any $C\leq C'\leq C''$
 we have $p_{C,C''}=p_{C',C''}\circ p_{C, C'}$.  \qed
 \end{prop}

 \paragraph{B. Role of the varieties $X_C$.} 
  The natural next step would be to pass from the diagram $(X_C)$ to an
 $\Hc$-schober (Definition \ref{def:H-flober}), a
diagram   $\Fen$
  of triangulated categories  $(\Ec_C, \gamma_{CC'}, \delta_{C'C})$. 
   Intuitively, $\Ec_C$ should be some coherent
 derived category of $X_C$, while $\gamma_{CC'}$  and $\delta_{C'C}$
 should be the pushforward and
 pullback along $p_{CC'}$.  The diagram $\Fen$ should represent a perverse schober on
 $\hen^*_\CC$ smooth with respect to $\Hc_\CC$. 
 
  However,  the varieties $X_C$ can be singular, 
 so a straightforward  implementation of this idea is  difficult.  
 For instance, in order to have both inverse and direct images well defined in the singular case,
  we would be forced to work with
 (left) unbounded derived categories which will destroy  the Grothendieck groups
 and meaningful decategorifications.

 One possible approach is to construct natural desingularizations
 $\wt X_C$ of the  singular $X_C$ and to define $\Ec_C = D^b(\wt X_C)$. 
 We make first steps in this direction in \S \ref{sec:schubert}.

 However, Theorem \ref{thm:H^0-flober} suggests the following conjecture describing the analog
  of $\HH^0(\hen^*_\CC,\Fen)$. This
  conjecture  can be formulated
 without  reference to a choice of $\wt X_C$.

 \begin{conj}
 The intersection of the pullbacks of the categories $D^b(X_w)$ inside $D^b(X_0)$ is equivalent to $\Perf(Z)$. 
 \end{conj}

  \paragraph{C.  Reduction of $X_C$ to  $X_0$ for a Levi.} 
  In order to indicate the dependence of our varieties $X_C, Z_C$ as well as the arrangement $\Hc$,
  the Cartan subalgebra $\hen$ etc. 
  on $\gen$ we will use, when necessary, the notations $X_C(\gen), Z_C(\gen), \Hc(\gen), \hen(\gen)$, etc. 
  
  \vskip .2cm
  
  Let $C$ be a cell of $\Hc(\gen)$, and $\pen_C$ be the corresponding parabolic subalgebra
  containing $\hen$, see Proposition \ref{prop:cells-parab}. Let $\len_C$ be the Levi quotient
  of $\pen_C$ and $\men_C$ be the semi-simplification (quotient by the center) of $\len_C$.
  We realize $\men_C$ as the subalgebra in $\gen$ generated by the root generators
  $e_\alpha$ with $\alpha|_C=0$. 
  
  At the same time, let $\RR\cdot C\subset\hen^*_\RR(\gen)$ be the real subspace spanned by $C$. Then
  we have identifications
  \[
  \hen^*_\RR(\len_C) \simeq  \hen^*, \quad \hen_\RR^*(\men_C) \simeq \hen^*_\RR(\gen)\bigl/\RR\cdot C. 
  \]
  and
  \[
  \Hc(\len_C) \simeq \Hc(\gen)^{\geq C}, \quad \Hc(\men_C) \simeq \Hc(\gen)/C.
  \]
  Here $\Hc(\gen)^{\geq C}$ is the sub-arrangement of $\Hc(\gen)$ consisting of the hyperplanes
  which contain $C$, while $\Hc(\gen)/C$ is the {\em quotient arrangement} in the quotient space 
  $\hen^*_\RR(\gen)\bigl/\RR\cdot C$ whose hyperplanes are the images of the hyperplanes from 
  $\Hc(\gen)^{\geq C}$, see, e.g., \cite{KS}. Accordingly, the cells of $\Hc(\gen)/C$ are the images,
  under the quotient map, of the cells $C'$ of $\Hc$ satisfying $C'\geq C$. We will denote such images
  by $C'/C$. 
  
  \vskip .2cm
  
  In particular, for any $C'\geq C$ we have the
   varieties $Z_{C'}(\len_C)$ and $Z_{C'/C}(\men_C)$
  associated to the reductive groups $\len_C$ and $\men_C$. We have the composite projection
  $\pen_C\to\len_C\to\men_C$. 
  
  \begin{prop}\label{prop:Z_c'-pullback}
  For any $C'\geq C$ we have an identification
  \[
  Z_{C'}(\gen) \,\,\simeq \,\, G\times_{P_C} \bigl( \pen_C \times_{\men_C} Z_{C'/C}(\men_C)\bigr). 
  \]
  These identifications take, for any $C''\geq C'\geq C$, the projection $Z_{C''/C}(\men_C)\to Z_{C'/C}(\men_C)$,
  to the projection $Z_{C''}\to Z_{C'}$. 
  \end{prop}
  
  The proposition means that $Z_{C'}(\gen)$ is obtained  from $Z_{C'/C}(\men_C)$
  by first performing a pullback onto $\pen_C$ and then forming the universal family over $F_C$. 
  
  \vskip .2cm
  
  \noindent{\sl Proof:} This is a direct consequence of the definition of $Z_{C'}$ as the universal family over
  $F_{C'}$  given by   \eqref{eq:Z=universal}. \qed

  \begin{cor}\label{cor:x-wall-sm}
  (a) For any cell $C$ of $\Hc=\Hc(\gen)$ we have an identification
  \[
  X_C(\gen) \,\,\simeq \,\, G\times_{P_C}\bigl(\pen_C\times_{\men_C} X_0(\men_C)\bigr). 
  \]
  Here $X_0(\men_C)$ is the principal component of the
   fiber product of all the flopped Grothendieck resolutions of $\men_C$
  over $Z(\men_C)$. 
  
  \vskip .2cm
  
  (b) In particular,   if $C$ is a wall (codimension 1 cell), then $X_C(\gen)$ is smooth. It coincides with the fiber
  product $X_{w'}(\gen) \times_{Z_C(\gen) }X_{w''}(\gen) $, where $w', w''\in W$ are such that $C$  separates
  two chambers $C_{w'}$ and $C_{w''}$.
  
  \end{cor}
  
  The corollary implies that each $X_C(\gen)$ locally looks like the product of $X_0(\men_C)$ and
  a smooth manifold, because it is obtained from $X_0(\men_C)$ by first forming a pullback under a smooth map 
  and then forming a universal family over a smooth base. 
  
  \vskip .2cm
   
   \noindent{\sl Proof:} (a) follows directly from Proposition \ref {prop:Z_c'-pullback}. Part (b) follows
   from (a) since $\men_C\simeq\sen\len_2$ in this case, and the case of $\sen\len_2$
   corresponds to the Atiyah flop where the statement is well known, see \S \ref{sec:flops-and-schobers}B. \qed

  \paragraph{D. The variety of $\gen$-simplices.} 
 By an (ordered) $\gen$-{\em simplex} we mean a pair $\hen\subset\ben\subset\gen$ formed by
 a Cartan subalgebra $\hen$ and a Borel $\ben$ containing it. 
 Denote by $h,b$ the dimensions of the $\hen$'s and $\ben$'s. Then we have the quasi-projective variety
  $T^\circ_\gen\subset \Fl(h,b,\gen)$
 formed by ordered $\gen$-simplices $\hen\subset\ben$. It is nothing but the quotient $G/H$ where
 $H\subset G$ is the maximal torus.

 \begin{ex}
 Let $\gen=\sen\len_n$. A Cartan subalgebra $\hen\subset\gen$ is the same as a choice of
 $n$ points $\{x_i\}_{i\in I}$, $|I|=n$ in $\PP^{n-1}$  which are in (linearly) general position, 
 i.e., form an {\em unordered simplex}. Explicitly, $\hen$ consists of linear operators which
are diagonal in a basis   of $\k^n$ formed by vectors $e_i$  that represent  $x_i$.
 A choice of a Borel $\ben\supset\hen$ amounts to ordering the $x_i$, i..e, writing them
 as $(x_1, \cdots, x_n)$ (an ordered simplex). The variety $T^\circ_\gen$ is in this case
 identified with the open set $(\PP^{n-1})^n_\gene\subset (\PP^{n-1})^n$
 formed by $n$-tuples of points in general position.

 \end{ex}
 
  We recall  that $G$-orbits on $F\times F$ are parametrized by elements of $W$ (the Bruhat
  decomposition). In other words, 
    for any two Borel subalgebras $\ben, \ben'$ we have their {\em relative position}
  $d(\ben, \ben')\in W$. Thus $d(\ben,\ben')=1$ iff $\ben=\ben'$ 
   while $d(\ben, \ben')=w_0$
  (the maximal element in $W$), iff $\ben$ and $\ben'$ are in general position (i.e., 
  $(\ben, \ben')$ lies in the open
  $G$-orbit in $F\times F$). 
  In particular, given $\hen\subset\ben$, for any $w\in W$ we have a unique Borel
  $\ben_w = \ben_w(\hen) \supset \hen$ such that $d(\ben_w, \ben)=w$. 
  
  \vskip .2cm
  
  Let $F^W$ be the product of $W$ copies of $F$; its points will be written as systems
  $(\ben_w)_{w\in W}$ of Borels. We then have a regular embedding
  \[
  u: T^\circ_\gen \lra F^W, \quad (\hen\subset\ben) \mapsto (\ben_w(\hen))_{w\in W}. 
  \]
  If we view $T^\circ_\gen$ as $G/H$ and $W$ as the quotient of the normalizer $N(H)$ by $H$, then
  the $w$-component of $u$ is the second component  of  the $w$-action map
  \[
  \on{Ad}(w) : G/H\lra G/H, \quad (\hen, \ben) \mapsto (\hen, \ben_w(\hen)). 
  \]

  \begin{Defi}  
  The {\em variety of $\gen$-simplices} $T_\gen \subset F^W$ is defined as the
  closure of the image of $u$. 
  \end{Defi}

  In the case  $\gen=\sen\len_n$ there is a more direct description of $T_\gen$, cf.
  \cite{babson1}. Instead of $F^W=F^{n!}$, consider the product of Grassmannians
  $\Gr(|I|, \k^n)$ labelled by all subsets $I\subset \{1,\cdots, n\}$. We then have the regular embedding
  \be\label{eq:map-grass-sim}
  \begin{gathered}
  v: (\PP^{n-1})^n_\gene \lra \prod_{I\subset \{1, \cdots, n\}} \Gr(|I|, \k^n),
  \\
  (x_1, \cdots, x_n) \,\,\mapsto \,\, \bigl( \on{Span}(x_i)_{i\in I}\bigr)_{I\subset \{1, \cdots, n\}}, 
  \end{gathered}
  \ee
which associates to $(x_1, \cdots, x_n)$ the system of the subspaces spanned
by all possible subsets of $\{x_1, \cdots, x_n\}$. 

\begin{prop}\label{prop:T-sln} 
$T_{\sel_n}$ is identified with the closure of the image of $v$. 
\end{prop}

\noindent {\sl Proof:} Denote $\Gr_I = \Gr(|I|, \k^n)$ the factor in the target of $v$ corresponding to $I$. 
The flag variety $F$ for $\sel_n$ is, classically,  embedded in the product of Grassmannians $\Gr(p, \k^n)$,
$p=1,\cdots, n-1$ (as the incidence variety). Let us think of the factor $\Gr(p, \k^n)$ in this
embedding as $\Gr_{\{1,\cdots, p\}}$. 

Further, for each $w\in W$ we consider the $w$-th factor $F_w=F$  in $F^W$ and embed in into the
product of Grassmannians $\Gr_{w(\{1, \cdots, p\})}$, $p=1,\cdots, n-1$.  By taking
the direct product of these embeddings over all $w\in W$, we 
  embed
$F^W$ into the following product of Grassmannians:
\[
\eta: F^W \lra \prod_{p=1}^{n-1} \prod_{I\subset \{1, \cdots, n\}, \atop |I|=p} \Gr_I^{W_{\{1,\cdots, p\}}^I}
\]
  where $W_{\{1,\cdots, p\}}^I\subset W$
is the set of permutations $w\in W=S_n$ such that $w(\{1,\cdots, p\}) = I$. 
We now notice:

\begin{lem}
The composition $\eta\circ u$ takes values in the product, over all $p, I$, of the small
diagonals  $\Gr_I \subset \Gr_I^{W_{\{1,\cdots, p\}}^I}$. The resulting morphism
from $T_{\sel_n}^\circ = (\PP^{n-1})^n_\gene$ to $\prod_I \Gr_I$
is identified with $v$. 
\end{lem}

\noindent{\sl Proof of the lemma:} Let $x=(x_1,\cdots, x_n)\in (\PP^{n-1})^n_\gene$
correspond to $\hen\subset\ben\subset \sel_n$. Then, for $w\in S_n$, the Borel subalgebra
$\ben_w(\hen)$ corresponds to the flag 
\[
{\Span}(x_{w(1)})\subset \Span (x_{w(1)}, x_{w(2)})\subset \cdots \Span(x_{w(1)}, \cdots,
x_{w(n-1)})\subset\k^n. 
\]
Our statement follows from this immediately. \qed

Proposition \ref{prop:T-sln}  is proved.

\paragraph{E. $X_0$ and the variety of simplices.}
We now realize the biggest horizontal component $X_0\subset \prod_{w\in W} X_w$
as a ``fibration in cones" over the variety  of simplices $T_\gen$. 

We denote by 
\be\label{eq:h-gene}
\ul\hen_\gene\,\,=\,\,\bigl\{  \bigl( (\hen\subset\ben), x\bigr) \in  T_\gen^\circ\times \gen
 \bigl| \,x\in\hen\bigr\} \,\,\buildrel\alpha\over\hookrightarrow  \,\, T_\gen\times\gen
 \ee
 the universal bundle of Cartan
subalgebras over $T_\gen^\circ$. 

\begin{prop}\label{prop:X0-cone}
$X_0$ is identified with the closure of $\ul\hen_\gene$ in $T_\gen\times\gen$. In particular,
it is 
  a closed subvariety in $T_\gen\times\gen$, conic (dilation invariant) in the $\gen$-direction.
\end{prop}

\noindent{\sl Proof:}
We note that we have a projection $Z_\rss\to T^\circ_\gen$ whose fiber over $(\hen\subset\ben)$
is $\hen_\reg$, the open subset of regular elements in $\hen$. 
So we write $Z_\rss=\ul\hen_\reg$ which exhibits it as an open dense subset in $\ul\hen_\gene$.

 Next, each $X_w=w^*\wt\gen$ is identified with $\wt\gen\subset F\times\gen$, but with projection
to $Z$ twisted by $w$. So  $X_0$ can be conisdered as  a subvariety in $\wt g^W = \prod_{w\in W} \wt\gen$,
and, as such, is embedded into $F^W\times\gen^W$. By definition, $X_0$ is the closure 
(in $\wt\gen^W \subset F^W\times\hen^W$) of the image of
$Z_\rss\subset Z$ under the system of maps inverse to the projections $X_w\to Z$
(biregular over $Z_\rss$). We now notice that the composite map
$Z_\rss\to F^W\times \gen^W$ lands in $T_\gen\times \gen$ (with $\gen\subset\gen^W$
being the small diagonal). This means that $X_0$ is identified with the closure of
$\ul\hen_\reg$ in $T_\gen\times\gen$ which is the same as the closure of the vector bundle
$\ul\hen_\gene$, in which it is open and dense. \qed

\vskip .2cm

\paragraph{F. The variety of reductions and the Cartanization of $X_0$.} 
Proposition \ref{prop:X0-cone} gives an affine morphism $\rho: X_0\to T_\gen$ whose
fiber over a generic point $\bb=(\hen\subset\ben)\in T^\circ_\gen$ is the Cartan subalgebra
$\hen$. The fiber $\rho^{-1}(\bb)$ over an arbitrary point $\bb\in T_\gen$ can be seen as the union of
the limit positions of such Cartans for all 
 1-parameter curves ($\k[[t]]$-points) $\bb(t)$ in $T_\gen$ such that $\bb(0)=\bb$ and 
 $\bb(t)\in T^\circ_\gen$ for $t\neq 0$. In particular, if we write $\bb\in T_\gen\subset F^W$
  as a system of Borels $(\ben_w)_{w\in W}$,
 then $\rho^{-1}(\bb)$ is contained in the intersection $\bigcap_{w\in W}\ben_w$. 
 We now produce  partial desingularization of $T_\gen$ and $X_0$. 
 
 \begin{Defi} \cite{iliev1, iliev2}
 Let $R^\circ_\gen \subset\Gr(h,\gen)$ be the variety of all  Cartan subalgebras.
 The {\em variety of
 reductions} for $\gen$ is defined to be the closure
  $R_\gen = \ol{ R_\gen^\circ} \subset\Gr(h,\gen)$.  
 \end{Defi}
 
   By construction,
 $R_\gen$ carries the universal bundle  $\ul \aen\subset R_\gen\times\gen$
 of rank $h$ whose fibers are ``limit positions of Cartan subalgebras''. They are 
 abelian subalgebras in $\gen$.  We then have the embedding
 \be\label{eq:u-hat}
 \wh u: T^\circ_\gen \lra F^W \times R_\gen, \quad (\hen\subset\ben) \mapsto \bigl(
 (\ben_w(\hen))_{w\in W}, \hen\bigr)
 \ee
 and we define the {\em Cartanized variety of $\gen$-simplices} 
  $ \wh T_\gen$ to be the closure of the image of $\wh u$. Thus we have a regular birational
  map $\tau: \wh T_\gen\to T_\gen$ which is the composition of the closed embedding
  $\wh T_\gen \subset T_\gen \times R_\gen$ and the projection of the product to
  $T_\gen$. 
  
  \vskip .2cm
  
  Similarly, we consider the embedding 
  \[
  \ul\hen_\gene \,  \,\,\buildrel\wh \alpha\over\hookrightarrow  \,\, T_\gen  
    \times \gen\times R_\gen, 
  \quad ((\hen\subset\ben), x) \mapsto \bigl( \alpha ((\hen\subset\ben), x), \hen\bigr),
 \]
 where $\alpha$ is as in \eqref {eq:h-gene}. Define the variety $\wh X_0$ as the closure of
 the image of $\wh\alpha$.   
 \begin{prop}\label{prop:hat-x-0}
 $\wh X_0$ lies (after permuting the second and third factors in the target of $\wh\alpha$)
 in $\wh T_\gen\times\gen\subset T_\gen\times R_\gen \times\gen$.
 Further,  
 we have a commutative diagram
 \[
 \xymatrix{
\ul\aen \ar[d]& \ar[l] \wh X_0 \ar[d]_{\wh\rho} \ar[r]^\sigma& X_0\ar[d]^\rho
 \\
 R_\gen&\ar[l] \wh T_\gen\ar[r]_\tau& T_\gen
 }
 \]
 where:
 \begin{enumerate}
 \item[(a)]  $\sigma$ and $\tau$ are    regular,  birational and proper.
 
 \item[(b)] The left square is Cartesian, in patricular, $\wh\rho$ exhibits $\wh X_0$ as the total
 space of a rank $h$ vector bundle on $\wh T_g$.
 
 \item[(c)] For $\bb\in T_\gen$ the preimage $\rho^{-1}(\bb)\subset\gen$ is the union, over
 all $\wh \bb\in\tau^{-1}(\bb)$, of the vector subspaces $\wh \rho^{-1}(\wh \bb)$. 
  \end{enumerate}
  \qed
 \end{prop}
 
  \noindent{\sl Proof:} The first statement (about $\wh X_0$) follows from the construction of $\wh T_\gen$
  as the closure of the lift of $T^\circ_\gen$ into $T_\gen\times R_\gen$.   The construction and commutativity
  of the diagram is just a restatement of the steps above. Let us prove the claimed properties of the diagram.
  
  \vskip .2cm
  
  (a) $\sigma$ is  regular  and proper since it is obtained from a proper morphism
  (projection)  $T_\gen\times \gen\times R_\gen\to T_\gen\times \gen$ by taking the closure of a lift of $\ul\hen_{\on{gen}}$
  to the source and mapping in to the closure of the image of this lift in the target. It is birational since it is
  an isomorphism on $\ul\hen_{\on{gen}}$.  The argument for $\tau$ is similar.
  
  \vskip .2cm
  
  (b)This  follows because we can also define $\wh X_0$ as the closure of  $\ul\hen_{\on{gen}}$ in the total space of
   the pullback of the  vector bundle
  $\ul \aen$ to $T_\gen\times R_\gen$.
  
  \vskip .2cm
  
  (c) By definition of $X_0$ as the closure,
   $\rho^{-1}(\bb)$ is the union of all the limit positions of the Cartans associated to $\bb'\in T^\circ_\gen$
  approaching $\bb$ in $1$-parameter families. These limit positions will represent points of $R_\gen$
  and thus points  of $\wh T_\gen$, more precisely of  $\wh\bb\in\tau^{-1}(\bb)$. 
  The individual elements of  these limit positions will become vectors in the
  fibers of the tautological bundle  $\aen$ over $R_\gen$, i.e., vectors in the vector subspaces
  $\wh\rho^{-1}(\wh \bb)$.  \qed

 \vfill\eject
 
 \section{The flop diagram for $\sel_3$ and Schubert's variety of complete triangles}\label{sec:schubert}
 
  \paragraph{A. Complete triangles vs. Cartanization.} From now on we restrict the discussion to $\gen=\sel_3$.
  The variety $T=T_{\sel_3}$ of $\sel_3$-simplices is the classical {\em variety of triangles} \cite{schubert, semple, magyar}.
  It has dimension $6$ and is embedded into $(\PP^2)^3\times (\check \PP^2)^3$,
  see \eqref{eq:map-grass-sim}.
 Here $\check \PP^2 = \Gr(2,\k^3)$  is the projective plane of lines in $\PP^2$.   
  Below are some of  the well known properties of $T$. 
  
  \begin{prop}
  (a) $T$ coincides with the incidence variety
  \[
 \bigl\{ (x_1, x_2, x_3, l_{12}, l_{13}, l_{23}) \in (\PP^2)^3 \times (\check \PP^2)^3 \bigl | \,\,\,
x_i, x_j \in l_{ij}\bigr\}. 
  \]
  (b) The singular locus $T^\sing$ consists of $(x_i, l_{ij})$
  such that $x_1=x_2=x_3$ and $l_{12}=l_{13}=l_{23}$. It is thus identified with $F$, the flag variety for $\sel_3$. 
  
  \vskip .2cm
  
  (c) Near $T^\sing = F$, the variety $T$ locally looks like the product of a $3$-dimensional affine space
  and a $3$-dimensional quadratic cone. 
  \end{prop}
  
  This has been proved in \cite{semple}, Th. I,
  where explicit charts were constructed. For a more modern treatment of (a) and (b), see \cite{magyar}, Th. 1(b)
  and Lemma 7.

  \vskip .2cm
  
  Because of the $3$d quadratic cone nature of the singularities of $T$, it has two small desingularizations,
  connected by a flop:
  \be
  T^\Sch\buildrel \tau_\Sch \over\lra T \buildrel \tau_\FM\over\longleftarrow T^\FM.
  \ee
  Thus for any $\bb\in T^\sing = F$, the preimages $\tau_\Sch^{-1}(\bb)$ and $\tau_\FM^{-1}(\bb)$
  are isomorphic to $\PP^1$. 
  The variety $T^\Sch$ is known as  {\em Schubert's variety of complete triangles}  and $T^\FM$
  is the {\em Fulton-MacPherson blowup} of $(\PP^2)^3$. They are both  acted upon by 
  $SL_3$. Intrinsically, they are distinguished by the fact 
   that the action  of the maximal torus $(\GG_m)^2 \subset SL_3$ on  $T^\Sch$ has isolated fixed points,
   while the  fixed points in  $T^\FM$ can be non-isolated, see \cite{magyar, roberts}. 
   
   \vskip .2cm
   
   We recall \cite{magyar, roberts, semple}  the definition of $T^\Sch$ going back to Schubert \cite{schubert}.
   Consider the $6$-dimensional space $S^2(\k^{3*})$ of quadratic forms on $\k^3$. Every triple
   $(x_1, x_2, x_3)\in (\PP^2)^3_\gene$ defines a 3-dimensional subspace
   \[
   N_{x_1, x_2,  x_3} \,\,=\,\, \bigl\{ q\in S^2(\k^{3*}) \bigl| \, q(x_1) = q(x_2) = q(x_3) =0\bigr\}
   \]
  (the net of quadrics through $x_1, x_2, x_3$). 
   We then consider the embedding
  \[
  \begin{gathered}
  v_\Sch:  (\PP^2)^3_\gene \lra (\PP^2)^3 \times (\check \PP^2)^3\times \Gr(3, S^2(\k^{3*})), 
  \\
  (x_1, x_2, x_3) \mapsto 
  \bigl(  (x_1, x_2, x_3), (\Span(x_1, x_2), \Span(x_1, x_3), \Span(x_2, x_3)), N_{x_1, x_2, x_3}\bigr). 
  \end{gathered}
  \]
  By definition, $T^\Sch$ is the closure of the image of $v_\Sch$. The fact that $T^\Sch$
 thus defined,  is smooth, was proved directly  in \cite{semple}, Th. II. 
  
  \begin{prop}\label{prop:schb-cart}
  The Cartanized variety of triangles $\wh T$ is isomorphic to  Schubert's variety $T^\Sch$ (as a variety over $T$, so that $\tau$ in Proposition \ref{prop:hat-x-0}
 corresponds to $\tau_\Sch$). 
 In particular, $\wh T$ is smooth, while $T^\Sch$ carries a  rank 2 bundle of abelian Lie subalgebras
  in $\sel_3$. 
  \end{prop}
  
  \noindent {\sl Proof:} This is a consequence of the interpretation, given in \cite{iliev1}, 
  of the variety $R = R_{\sel_3}$ of reductions for $\sel_3$. A Cartan subalgebra $\hen\subset\sel_3$ can be
  seen as  corresponding to some $x=(x_1, x_2, x_3)\in  (\PP^2)^3_\gene$ defined uniquely
  up to permutation; we write $\hen=\hen_x$. Thus, tautologically,
  $R$ is the closure of the image of the embedding
  \[
   (\PP^2)^3_\gene/S_3 \lra \Gr(2, \sel_3), \quad x\mapsto \hen_x. 
  \]
  It was shown in \cite{iliev1} Prop. 4.1 that $R$ is equal to the subvariety in $ \Gr(2, \sel_3)$
  formed by all $2$-dimensional abelian Lie subalgebras, i.e., to the intersection
  \[
  \Gr(2, \sel_3) \cap \PP\bigl(  \Ker\bigl\{ \Lambda^2 \sel_3 \buildrel [-,-]\over\lra \sel_3\bigr\}\bigr)
  \]
  in the Pl\"ucker embedding of  $\Gr(2, \sel_3)$.

  At the same time, any  $x$ as above gives 
  a 3-dimensional subspace $N_x\subset S^2(\k^{3*})$, so we have the embedding
  \[
  \nu:   (\PP^2)^3_\gene/S_3 \lra \Gr(3,S^2(\k^{3*})) \subset \PP\bigl(\Lambda^3(S^2(\k^{3*}))\bigr), 
  \quad x\mapsto N_x. 
  \]
 It was proved in  \cite{iliev1}, Th. 4.2, that the closure of the image of $\nu$ is identified with $R$. 
  This is based on the isomorphism of $\sel_3$-modules
  \[
  \Ker\bigl\{ \Lambda^2 \sel_3 \buildrel [-,-]\over\lra \sel_3\bigr\} \,\,\simeq \,\, \Lambda^3  S^2 (\k^{3*}). 
  \]
   describing the ambient spaces of the two Pl\"ucker embeddings, see \cite{iliev1}, Lemma 4.3. 
   
   \vskip .2cm

  We can therefore compare the embeddings $\wh u$ from \eqref{eq:u-hat} defining $\wh T$ and the embedding
  $v_\Sch$ defining $T^\Sch$. We see that they both factor through
  the same map into $ (\PP^2)^3 \times (\check \PP^2)^3\times R$ with $R$ embedded
  in $\Gr(2, \sel_3)$ in the first case and in $\Gr(3, S^2(\k^{3*}))$ in the second case. \qed
  
  \begin{cor}
  The Cartanized variety $\wh X_0$ for $\gen=\sel_3$ is smooth.
  \end{cor}
  
 \noindent{\sl Proof:}  Indeed, $\wh X_0$ is the total space of a rank 2 vector bundle
  (of abelian Lie subalgebras in $\sel_3$)
  over the smooth variety $\wh T=T^\Sch$. 
  
  \qed
  
  Using Proposition \ref{prop:hat-x-0}, we can now get a more detailed picture of the original 
  variety $X_0\subset T\times\sel_3$ in
   terms of the projection $\rho: X_0\to T$. 
  
  \begin{cor}\label{cor: lim-cartans} 
  (a) If  $\bb\in T$ is a smooth point, then $\rho^{-1}(\bb)$ is a 2-dimensional abelian Lie
  subalgebra in $\sel_3$. 
  
  \vskip .2cm
  
  (b) If $\bb\in T^\sing = F$ corresponds to a Borel subalgebra $\ben\subset\sel_3$, then 
  $\rho^{-1}(\bb)= [\ben, \ben]$ is 3-dimensional. 
  
  \vskip .2cm
  
  (c) In the situation (b), the preimage $\tau^{-1}(\bb)\subset \wh T$ is identified with the variety
  of all 2-dimensional abelian Lie subalgebras in $\nen= [\ben,\ben]$. Such subalgebras  are precisely the
  2-dimensional vector subspaces in $\ben$ containing the center $[\nen,\nen]$. In particular,
  there are $\PP^1$ of them. 
  
  \vskip .2cm
  
  (d) The variety $X_0$ is singular along $T^\sing\times\{0\}\subset T\times\sel_3$. 
   \qed
  \end{cor}
  
  \noindent{\sl Proof:}  (a)  follows from the identification of $\tau$
  with $\tau_\Sch$ and the known fact that $\tau_\Sch$ is bijective outside of $T^\sing$. 
  
    \vskip .2cm
  
  Part (b)  is a consequence of 
   the following elementary lemma which explains
  the appearance of $\nen$.
  
  \begin{lem}
  Let $p_i(t), i=1,\cdots, 3$, be three $\k[[t]]$-points of $\PP^2$ with the following properties:
  \begin{enumerate}
  \item[(1)] For $t\neq 0$, i.e., as $\k((t))$-points, the $p_i(t)$ are in general position.
  
  \item[(2)] At $t=0$, all three $p_i(t)$ evaluate to the same point $p\in\PP^2(\k)$,
  and the three lines $(p_i(t), p_j(t))$ evaluate to the same line $l\subset \PP^2$ (containing $p$).
   \end{enumerate}
   Let $\ben\subset\gen = \sel_3$ be the Borel subalgebra fixing the flag $(p,l)$ and $\nen=[\ben,\ben]$.
   Let further $x(t)\in\sel_3(k[[t]])$ be a family of matrices such that for $t\neq 0$ the matrix
   $x(t)$ has three points $p_i(t)$ as  eigendirections. Then $x(0)\in\nen$. Moreover,
   each element of $\nen$ can be obtained as $x(0)$ for an appropriate $p_i(t)$, $i=1,\cdots, 3$ and
   $x(t)$ as above.
  \end{lem}
  
  \noindent{\sl Proof of the lemma:} Clearly $x(0)$ fixes $p,l$ so $x(0)\in\ben$. Let $\xi$ be the gobal vector field on $\PP^2$
  given by $x(0)$. The two independent
  eigenvalues of $x(0)\in\sel_3$ can be read off the two eigenvalues of the transformation induced by $\xi$ on the tangent
  space $T_p\PP^2$. Now, the first eigenvalue (on the line $T_p l$) vanishes because the points $p_1(t)$ and $p_2(t)$
  specialize for $t\to 0$, to $p$ and the line $(p_1(t), p_2(t))$ specializes to $l$. The second eigenvalue
  (on the quotient $T_p\PP^2/T_pl$) vanishes because the third point $p_3(t)$ is such
  that the lines $(p_i(t), p_3(t))$, $i=1,2$ specialize, as $t\to 0$,  to the same line $l$.
  This shows that $x(0)\in\nen$. The fact that any element of $\nen$ can be obtained like this is obvious
  and left to the reader. \qed

  \vskip .2cm
  
  Part (c) follows from (b), from 
   Proposition \ref{prop:hat-x-0}(c) and the fact (\cite{iliev1}
  Prop. 4.1) that $R$ coincides with the variety of all 2-dimensional
  abelian Lie subalgebras in $\sel_3$. 
  
    \vskip .2cm
  
  To see (d), consider the Zariski tangent space
  to $X_0$ at $(\bb, 0)\in T^\sing\times\{0\}$. 
  It contains, as transverse direct summands, first,  the tangent space to $T$ at $\bb$
  (of dimension $7$), and, second,  the fiber of the projection $\rho$ over $\bb$, i.e., $\ben$
  (of dimension $3$).  So it is $10$-dimensional, while $X_0$ is $8$-dimensional. 
  This finishes the proof of Corollary \ref{cor: lim-cartans}.

 \paragraph{B. The flop diagram for $\sel_3$.} 
 The coroot arrangement $\Hc$ for $\sel_3$ consists of three lines in the real plane
 $\hen^*_\RR\simeq\RR^2$. It
  has $13$ cells: six chambers $C_w$,
 six walls (open half-lines) and the point $0$. Accordingly, the diagram $(X_C, p_{CC'})_{C\leq C'}$
 consists of 13 varieties: 
 
 \begin{itemize}
 \item [(2)]  Six flopped Grothendieck resolutions $X_w=X_{C_w}$, $w\in S_3$. They are smooth.
 
 \item[(1)] Six binary fiber products $X_C$ corresponding to the half-lines $C$. They are
 smooth, see Corollary \ref{cor:x-wall-sm}. 
 
 \item[(0)] The ``central''  variety $X_0$, singular.  
  \end{itemize}
  
  \begin{Defi}\label{def:sl3-flop-d}
  (a) The {\em flop diagram for $\sel_3$} is the diagram 
  $(Y_C, l_{CC'})_{C\leq C'}$   of proper birational regular maps between smooth varieties,
  obtained from $(X_C, p_{CC'})$
  by replacing $X_0$ with
   its Cartanization
  $\wh X_0\buildrel\sigma\over\to X_0$ and leaving the other varieties unchanged. That is, we put 
  \[
  Y_C = \begin{cases}
  X_C, & \text{if } C\neq 0,
  \\
  \wh X_0,& \text{if } C=0;
 \end{cases}
  \quad\quad 
  l_{CC'} = \begin{cases}
  p_{CC'}, & \text{if } C\neq 0,
  \\
  p_{0C'}\circ \sigma, & \text{if }  C=0<C'. 
  \end{cases}
  \]
  (b) The {\em Schubert transform} is the diagram of birational morphsms
  \[
  Y_{-C} \buildrel l_{0, -C}\over\lla  Y_0 \buildrel l_{0,C}\over\lra Y_C
  \]
   associated to any $1$-dimensional face (ray) $C$ of $\Hc$.  
    \end{Defi} 
    The remainder of this section prepares the ground for the 
   study of  the diagram formed by the categories $D^b(Y_C)$ in 
     the next section \ref {sec:sl3-flob}.

   \paragraph{C. Central fibers. The partial triangle picture.}
We first  summarize the main features of the varieties $Y_C$ and $X_C$ for arbitrary $C$.

\begin{prop}
   (a) We have a diagram
   \[
   \xymatrix{
    F_C & \ar[l]_{ \rho_C} Y_C\ar[r]^{ \tau_C} & \gen
    \\
   &  F_C\ar[u] \ar[r]& \{0\} \ar[u],
   }
   \]
   in which the square is Cartesian.  
   
   (b) The map $\rho_C$ represents $Y_C$ as the total space of a vector bundle on $F_C$ which
   we denote $L_C$. Every fiber of $L_C$ is a Lie subalgebra in $\gen$.

  (c)  The map $(\rho_C, \tau_C):  Y_C\to  F_C\times\gen$ embeds $Y_C$
   as a closed subvariety in $ F_C\times\gen$. In other words, it realizes $ L_C$ as a vector subbundle of
   the trivial bundle  $F_C\times\gen$.

    \qed
   
     \end{prop}

   The variety $F_C =  \tau_C^{-1}(0)$ will be called the {\em central fiber}
   of $Y_C$. We now summarize the structure of the vector bundles $L_C$. 
   
   Denote by $F= F(1,2,\k^3)\subset \PP^2 \times \check \PP^2$ the flag variety of $SL_3$.

   \begin{prop}\label{prop:F_C}
    (a) If $\dim(C)=2$, then $F_C\simeq F$ has dimension $3$, the bundle $L_C$ has rank $5$ and
   every fiber of $ L_C$ is a Borel subalgebra in $\gen$. 
   
   \vskip .2cm
   
   (b) If $\dim(C)=1$, then $F_C$ is isomorphic, as a $SL_3$-manifold,  to $F\times_{\PP^2} F$ 
   or $F\times_{\check \PP^2} F$ and
   has dimension $4$. The bundle $L_C$ has rank $4$. 
   
   \begin{itemize}
   
 \item [(b1)]   Let  $ q= (p,l,l')\in F\times_{\PP^2} F$
   be a point with $l\neq l'$.  Then, in the above identification,
    the fiber of $L_C$ over $q$ is identified with the intersection of
   two Borels  corresponding to the flags $(p,l)$ and $(p, l')$, i.e. with the space of global vector fields on
   $\PP^2$ which preserve $p$, $l$ and $l'$. 
   At $q=(p,l,l)$ the fiber $L_q$ is the (well defined)
   limit position of the above intersections of  Borels corresponding to neighboring points $q'$.  
   Explicitly, $L_q$ consists of global vector fields on $\PP^2$ which preserve $p$ and $l$ and
   whose restricttion to $l$ has vanishing linear part at $p$. 
   
   \item[(b2)] Let  $q= (p, p', l) \in F\times_{\check \PP^2} F$ is a point with $p \neq p'$.  Then, in the above identification,
   the fiber of $L_C$ over $q$ is identified with the intersection of the Borels corresponding to the flags $(p,l)$ and
   $(p', l)$, i.e., with the space of global vector fields on $\PP^2$ preserving $p$, $p'$ and $l$. At $q= (p,p,l)$
   he fiber $L_q$ is the (well defined)
   limit position of the above intersections of  Borels corresponding to neighboring points $q'$.  
   Explicitly, $L_q$ consists of global vector fields  on $\PP^2$ whose linear part of $p$ (an endomorphism of $T_p\PP^2$
   preserving the subspace $T_p l$) induces zero endomorphism of the quotient $T_p\PP^2/T_pl$. 
       \end{itemize} 
   
    \vskip .2cm
   
   (c) If $C=\{0\}$, then  $F_C = \wh T$ is the Schubert space of complete triangles and has dimension $6$.
   The bundle $L_C$ has rank $2$ and is formed by   Cartan subalgebras in $\gen$ and their limit positions. 
   \end{prop}
   \qed

           \begin{wrapfigure}{r}{0.3\textwidth}
  
    \begin{framed}\raggedleft

  \btp[scale=.8]
     \centering
   \draw[line width=0.7] (-2,0) -- (2,0);
   \node at (-1,0) {$\bullet$}; 
  \node at (1,0) {$\bullet$};   
  \node at (0,2) {$\bullet$}; 
    \draw[line width=0.7] (-1.5, -1) -- (0.5, 3);
     \draw[line width=0.7] (1.5, -1) -- (-0.5, 3); 
  \node at (-1.5, 0.5) {$p_1$};  
   \node at (1.5, 0.5) {$p_2$};  
     \node at (0.5, 2) {$p_3$};   
     \node at (2.5,0) {$l_{12}$};   
\node at (-0.8, 3.4) {$l_{23}$};     
\node at (0.8, 3.4) {$l_{13}$};    

  \etp
  
  \caption{The symbolic triangle.}\label{fig:triangle}
    \end{framed}
     
    \end{wrapfigure}

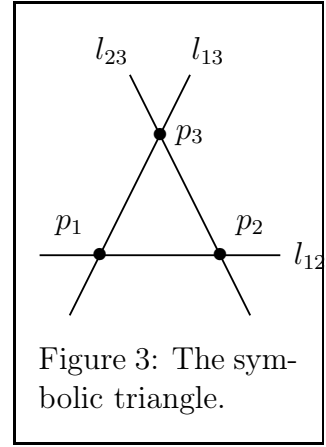
  \begin{figure}[H]
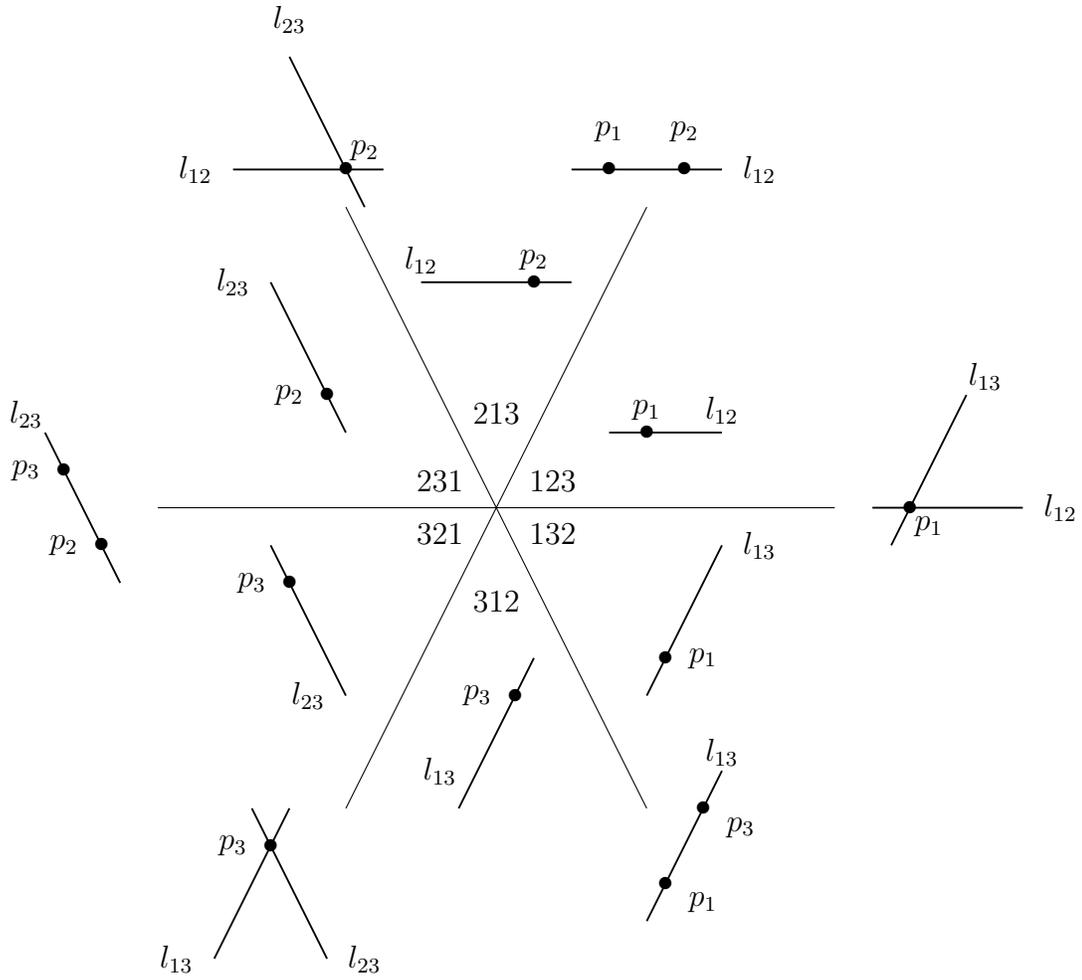

  \centering
 
 \btp[scale=.5]
 
 \draw (-9,0) -- (9,0); 
 \draw (-4, -8) -- (4, 8); 
  \draw (-4,  8) -- (4, -8); 
 
 \node at (1.5,0.7 ){$123$}; 
 \node at (0,2.5) {$213$}; 
  \node at (-1.5,0.7 ){$231$}; 
 \node at (-1.5, -0.7 ){$321$}; 
  \node at (0, -2.5) {$312$}; 
 \node at (1.5, -0.7 ){$132$};  
 
 \draw[line width=0.7] (3,2) -- (6,2);
 \node at (4,2){$\bullet$}; 
 \node at (4,2.6) {$p_1$}; 
 \node at (6, 2.6) {$l_{12}$}; 
 
 \draw[line width=0.7] (-2,6) -- (2,6); 
 \node at (1,6){$\bullet$}; 
 \node at (1, 6.6) {$p_2$}; 
\node at (-2, 6.6){$l_{12}$}; 

 \draw[line width=0.7] (-4,2) -- (-6,6); 
\node at (-4.5, 3) {$\bullet$};  
\node at (-5.5 ,3) {$p_2$}; 
\node at (-7,6) {$l_{23}$}; 

 \draw[line width=0.7] (-6, -1) -- (-4, -5); 
\node at (-5.5, -2) {$\bullet$}; 
\node at (-6.5, -2) {$p_3$}; 
\node at (-5, -5) {$l_{23}$}; 

 \draw[line width=0.7] (1, -4) -- (-1, -8); 
\node at (0.5, -5) {$\bullet$}; 
\node at (-0.5, -5) {$p_3$}; 
\node at (-1.5, -7) {$l_{13}$}; 

 \draw[line width=0.7] (4, -5) -- (6, -1); 
\node at (4.5, -4) {$\bullet$}; 
\node at (5.5, -4) {$p_1$}; 
\node at (7, -1){$l_{13}$}; 

\draw [line width=0.7]  (-12, 2) -- (-10, -2); 
\node at (-11.5, 1) {$ \bullet$}; 
\node at (-10.5, -1){$\bullet$};  
\node at (-12.5, 1){$p_3$}; 
 \node at (-11.5, -1) {$p_2$};  
\node at (-12.5, 2.5) {$l_{23}$};

\draw [line width=0.7] (-3,9) -- (-7, 9) ; 
\node at (-4,9) {$\bullet$}; 
\draw [line width=0.7] (-3.5, 8) -- (-5.5, 12) ; 
 \node at (-3.5, 9.5) {$p_2$}; 
 \node at (-8, 9) {$l_{12}$}; 
 \node at (-5.5, 13) {$l_{23}$}; 
 
 \draw [line width=0.7] (2,9) -- (6, 9) ; 
 \node at (3,9) {$\bullet$}; 
 \node at (5,9) {$\bullet$}; 
 \node at (3, 10) {$p_1$}; 
 \node at (5,10) {$p_2$}; 
 \node at  (7,9) {$l_{12}$}; 
 
  \draw [line width=0.7] (10,0) -- (14,0) ; 
 \node at (11,0) {$\bullet$}; 
   \draw [line width=0.7] (10.5, -1) -- (12.5, 3) ; 
 \node at (11.5, -0.5) {$p_1$}; 
\node at (15,0) {$l_{12}$};  
 \node at (13,3.5){$l_{13}$}; 
 
   \draw [line width=0.7] (4, -11) -- (6, -7) ; 
\node at (4.5, -10) {$\bullet$}; 
\node at (5.5, -8) {$\bullet$};  
 \node at (5.5, -10.5){$p_1$}; 
\node at (6.5, -8.5) {$p_3$};  
\node at (6, -6.5) {$l_{13}$};  
 
 \node at (-6, -9) {$\bullet$}; 
    \draw [line width=0.7] (-6.5, -8) -- (-4.5, -12) ; 
    \draw [line width=0.7] (-5.5, -8) -- (-7.5, -12) ;   
\node at (-7, -9) {$p_3$};     
 \node at (-8.5, -12) {$l_{13}$}; 
 \node at (-3.5, -12) {$l_{23}$}; 
 
 \etp
  \caption{The partial triangle notation for faces $C$ and varieties $Y_C, F_C$. }
  \label{fig:partial}
  
  \end{figure}

    It is convenient to depict the central fibers $F_C$ and the varieties $Y_C$ symbolically by 
    ``partial triangles'', i.e.,   by certain parts of the picture (Fig. \ref{fig:triangle})
     consisting of three points $p_1, p_2, p_3$ and three 
    lines $l_{12}, l_{13}, l_{23}$ joining them. 
                  This is depicted in  Fig. \ref{fig:partial}.  The six chambers labelled by permutations of $1,2,3$,
    correspond to six flags in the triangle, the six rays correspond to parts formed by either two
    vertices  and an edge  through them, or by two edges and a vertex common to them.
    The face  $\{0\}$ correspond to the full triangle of Fig. \ref{fig:triangle}.

   \vfill\eject

  If $C\neq \{0\}$, then the central fiber $F_C$ is simply the incidence variety associated to
  the corresponding picture. For instance, if $C$ is the right  horizontal ray in Fig. \ref{fig:partial},
  then $F_C$  is formed by $(p_1, l_{12}, l_{23})\in\PP^2\times \check \PP^2\times\check \PP^2$
  such that $p_1\in l_{12}$ and $p_1\in l_{13}$. 
  
  If $C=\{0\}$, then the incidence variety associated to the picture is $T$, the space
  of triangles, and the central fiber $F_0$ is its desingularization $\wh T$. 
  
  \paragraph{D. Fiber products: 1-ray and 2-ray  varieties.} 
  Let us call any variety $Y_C=X_C$ associated to $C$ which is a ray (1-dimensional face) of $\Hc$,
  a {\em 1-ray variety}, and the corresponding central fiber $F_C = F\times_{\PP^2} F$ or
  $F\times_{\check \PP^2} F$ a {\em 1-ray central fiber}.  
  
  Note that  $F_C$ is just  an elementary
  Schubert correspondence  (i.e., the correspondence associated to a simple reflection in the Weyl group)
  in $F\times F$,  while $Y_C\subset \wt \gen\times\wt\gen$ is the lift of such a correspondence
  to the Grothendieck resolution. These are precisely the correspondences used by Bezrukavnikov-Riche
  \cite{BR} to construct an action of the braid group $\Br_\gen$ on $D^b(\wt\gen)$: they give the action of
  the generators of $\Br_\gen$.  
  
  In order to prove the collinear transitivity conditions for the 
$\sel_3$-flober in \S \ref{sec:sl3-flob},
   we will  need
  certain fiber products of the correspondences $Y_C$. 
  These products have been studied in \cite{riche}, \S 2. At the same time, related varieties
  (the central fibers of the fiber products) appear in  the classical study of triangle varieties \cite{roberts}. 
 We present a treatment which emphasizes this connection.


 \vskip .2cm
 
\begin{figure}[H]
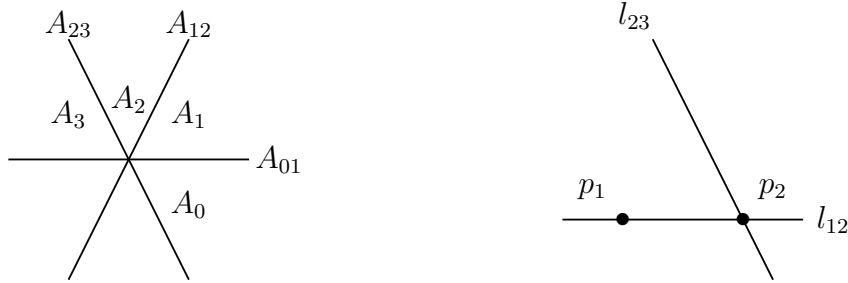
 
  \centering

  \btp[scale=.4]
     \centering
   \draw[line width=0.7] (-4,0) -- (4,0);
   \draw[line width=0.7] (-2,4) -- (2,-4);
  \draw[line width=0.7] (-2,-4) -- (2,4);
  \node at (2,1.5) {$A_1$}; 
 \node at (0,2) {$A_2$};   
 \node at (-2, 1.5) {$A_3$};  
   \node at (2,-1.5) {$A_0$}; 
 \node at (2, 4.5) {$A_{12}$}; 
 \node at (-2, 4.5) {$A_{23}$}; 
 \node at (5,0) {$A_{01}$}; 
    \etp
    \quad\quad\quad\quad
      \quad\quad\quad\quad
  \btp[scale=.8]
     \centering
   \draw[line width=0.7] (-2,0) -- (2,0);
   \node at (-1,0) {$\bullet$}; 
  \node at (1,0) {$\bullet$};

     \draw[line width=0.7] (1.5, -1) -- (-0.5, 3); 
  \node at (-1.5, 0.5) {$p_1$};  
   \node at (1.5, 0.5) {$p_2$};  
 
     \node at (2.5,0) {$l_{12}$};   
\node at (-0.8, 3.4) {$l_{23}$};     
    
 \etp
     \caption{Four consecutive chambers and the partial triangle code for a 2-ray variety.}\label{fig:4-cham}
 
    \end{figure}

 Let $A_0, \cdots, A_3$ be four consecutive chambers of $\Hc$ and $A_{01}, A_{12}, A_{23}$ be the rays
 between them, see Fig. \ref {fig:4-cham}. 
 We denote by 
 \[
 Y_{(2)} \,=\, Y_{ 12, 23} \, := \,     Y_{A_{12}} \times_{Y_{A_2}} Y_{A_{23}}
 \]
 and call the {\em 2-ray variety}  the fiber product of two 1-ray varieties over $Y_{A_2} \simeq  \wt\gen$. 
 (The chamber $A_0$ is not used in this definition but will be used in defining the $3$-ray varieties later.)
 
  We now extend the partial triangle notation to 2-ray varieties. If  we assume that the
  positions of  $A_0, \cdots, A_3$  in Fig. \ref{fig:4-cham} correspond to Fig. \ref{fig:partial}
  (i.e., $A_1$ corresponds to $123$, $A_2$ corresponds to $213$ etc.)
 then $Y_{(2)}$ is encoded by the partial triangle  in  Fig. \ref{fig:4-cham} on the right. 
 We define the {\em 2-ray central fiber} as the fiber product of the central fibers of the factors of $Y_{(2)}$, i.e., as
  the 
 incidence correspondence associated to the partial triangle code: 
 \[
 \begin{gathered}
 F_{(2)} \,=\,  F_{ 12, 23} \, := \,     F_{A_{12}} \times_{F_{A_2}} F_{A_{23}}\,\,=
 \\
 =\,\,
 \bigl\{ (p_1, p_2, l_{12}, l_{23}) \in \PP^2\times\PP^2\times\check \PP^2\times\check \PP^2
 \bigl| \,\, p_1, p_2\in l_{12}, \,\, p_2\in l_{23}\bigr\}. 
 \end{gathered} 
 \]
 This is nothing but the Demazure resolution of the Schubert correspondence in $F\times F$
 associated  to the reduced decomposition $(23) (12)$ in  $W=S_3$. 
 By construction we have a diagram
 \[
 \xymatrix{
 Y_{(2)}\ar[dr]_\rho \ar[r]^{\hskip -0.5cm i} & F_{(2)} \times\gen
 \ar[d]^{\pi}
 \\
 & F_{(2)}
  }
 \]
 where $i$ is a closed embedding and $\pi$ is the projection. Each fiber of $\rho$ is a Lie subalgebra in $\gen$.

 In the following it will be convenient to identify elements $\xi\in\gen$ with global
 regular vector fields on $\PP^2$.

 \begin{prop}\label{prop:Y2} 
 (a) $F_{(2)}$ is a smooth variety of dimension 5, which is a $\PP^1$-bundle over $F_{A_{12}}$,
 as well as a $\PP^1$-bundle over $F_{A_{23}}$. 
 
 \vskip .2cm
 
 (b) The map $\rho$ realizes $Y_{(2)}$ is the total space of a rank $3$ vector bundle $L_{(2)} = L_{12, 23}$ over 
 $F_{(2)} = F_{12, 23}$,
 whose fibers are Lie subalgebras in $\gen$. In particular, $Y_{(2)}$ is smooth.
 
 \vskip .2cm
 
 (c) The projection $\sigma: Y_{(2)}\to Y_{A_{12}}$ exhibits $Y_{(2)}$ as the blowup of $Y_{A_{12}}$
 (the total space of a rank $4$ bundle $L_{A_{12}}$ on $F_{A_{12}}$)
 along the total space of a rank $2$ subbundle $E$. 
 Explicitly, the fiber of $E$ at $(p_1, p_2,  l)$ consists of global vector fields $\xi$ on $\PP^2$
 which belong to $(L_{A_{12}})_{(p_1, p_2, l)}$ and whose linear part at $p_2$ is
 the scalar operator in $T_{p_2} \PP^2$. In the case $p_1=p_2$ this linear part must be $0$,
 since the correspinding elements of $\gen=\sel_3$ must be nilpotent. 
 
 Similarly for the projection to $Y_{A_{23}}$. 
 
 \vskip .2cm
 
 (d) The square 
 \[
 \xymatrix{
 Y_{(2)} \ar[d] \ar[r] & Y_{A_{12}}\ar[d]
 \\
 Y_{A_{23}} \ar[r]& Y_{A_2}
 }
 \]
 is Tor-independent, i.e., $Y_{(2)}$ is identified with the derived fiber product 
 $Y_{A_{12}} \times^h_{Y_{A_2}} Y_{A_{23}}$. 
 
 \end{prop}
 
 \noindent {\sl Proof:} (a) is obvious.  Let us prove (b).
 By definition of the fiber product, the
 fiber of $\rho$  over $q = (p_1, p_2, l_{12}, l_{23}) $ is the intersection
 \[
 (L_{(2)})_q \,\, := \,\, (L_{A_{12}})_{(p_1, p_2, l_{12})} \,\cap \,  (L_{A_{23}})_{(p_2, l_{12}, l_{23})} 
 \,\,\subset \,\,\gen
 \]
 of two Lie subalgebras in $\gen$. If $q$ is a generic point (i.e., $p_1\neq p_2$ and $l_{12}\neq l_{23}$), 
 this intersection is simply the subalgebra in $\gen$ which preserves the  points
 $p_1, p_2$ and the lines $l_{12}$ and $l_{23}$, and this subalgebra has dimension 3. 
 We now claim that for any other  $q\in F_{(2)}$ the above intersection has the same
 dimension 3.  By $PGL_3$-equivariance, it suffices to check  one point in each
 $PGL_3$-orbit on $F_{(2)}$. We consider only the ``most degenerate" case
 $p_1=p_2 = p$ and $l_{12}=l_{23} = l$. In this case, by Proposition \ref{prop:F_C}(b), the subalgebra 
 $(L_{(2)})_q$ consists of global vector fields $\xi$ on $\PP^2$ which preserve $p$ and $l$ and
 whose linear part at $p$ has trivial restrictions to the tangent space  $T_pl$ and to the quotient
 $T_p(\PP^2)/T_pl$. Global $\xi$ preserving $p$ and $l$ form a Borel subalgebra $\ben\subset\gen$,
 and the linear parts in question correspond to the diagonal elements of a triangular matrix.
 So the 
  condition on the linear parts 
  simply means that $\xi$ lies in the commutant $\nen = [\ben, \ben]$. This commutant  has dimension $3$.
   \vskip .2cm
   
   (c) We look at the fibers of $\sigma$. Suppose we have a point $(q',\xi)\in Y_{A_{12}}$,
   so $q'= (p_1, p_2, l_{12}) \in F_{A_{12}}$ and $\xi\in (L_{A_{12}})_q$. Then $\sigma^{-1}(q',\xi)$
   consists of $(q',\xi')$ such that 
    $(q= p_1, p_2, l_{12}, l_{23})$ projects to $q'$ and  $\xi'=\xi  \in (L_{(2)})_q$. In other words,
    the freedom consists in choosing the new line $l_{23}$ which must be invariant under  $\xi$.

    Now, let  $\Lambda\in\End( T_{p_2}\PP^2)$ be 
     the linear part of $\xi$ at $p_2$. If $\Lambda$
      is not a scalar operator, there is precisely one possible
    $l_{23}$: the line tangent to the other eigen-direction of $\Lambda$
    (or to the unique eigen-direction, if $\Lambda$ is not semisimple).
     If, however, $\Lambda$ is a scalar operator, then any $l_{23}$ through $p_2$
    is good and we have a $\PP^1$ worth of choices. So $\sigma$ has fibers $\PP^1$ over the
    total space of the subbundle $E$ and is an isomorphism on the complement. Since
    both the source and target of $\sigma$ are smooth varieties of the same dimension $8$,
    the identification with the blowup follows.    
    
    \vskip .2cm
    
    Finally, part (d) follows from the next general proposition. It has been formulated in \cite{Kuz2} with the reference to \cite{Kuz1}, but the latter text does not explicitly contains the proof of the fact. We are thankful to A. Kuznetsov for explanations and reproduce the proof here for the sake of completeness. 
    
    \begin{prop}\label{prop:tor-ind}
    Let 
    \[
    \xymatrix{
    X_T = X\times_S T
    \ar[d]_{\psi'}  \ar[r]^{\hskip 1cm \varphi'} & X \ar[d]^{\psi}
    \\ 
    T\ar[r]_{\varphi} & S
    }
    \]
    be a Cartesian square  in which $X,T,S$ are smooth irreducible varieties.
    Suppose further that each component of $X_T$ has the same dimension which is equal to
    \[
    \dim(X_T) \,\,=\,\,\dim(X) + \dim(T) -\dim(S). 
    \]
    Then the square is Tor-independent. 
    \end{prop}
    
    \noindent{\sl Proof of Proposition \ref {prop:tor-ind}:}  The statement is obvious, if $\varphi$ is
    a smooth   moprhism.  Let us point another particular case.
    
    \begin{lem}\cite{Kuz1}
    Suppose, in addition, that $\varphi$ is a closed embedding (of a smooth subvariety into a smooth
    variety). Then the square is Tor-independent.
    \end{lem}
    
    \noindent{\sl Proof of the lemma:} The statement is local. Locally $T\subset S$ is given by
    $n$ independent equations $f_1 = \cdots = f_n=0$. So $\Oc_T$ has a Koszul resolution
    \[
    \Oc_T \sim \bigl(\Oc_S[\xi_1, \cdots, \xi_n], d\bigr), \quad \deg(\xi_i)=-1, \,\, d(\xi_i)=f_i. 
    \]
   Accordingly, $\Oc_T\otimes^L_{\Oc_S} \Oc_X$ is represented by the Koszul complex
    \[
  \Kc^\bullet  \,\,:= \,\,   \Oc_X \sim \bigl(\Oc_S[\xi_1, \cdots, \xi_n], d\bigr), \quad  d(\xi_i)=\wt f_i :=
  f_i\circ\psi \in\Oc(X). 
    \]
    But the condition on the dimension of $X_T$ implies that $\wt f_1,\cdots, \wt f_n$ form
    a regular sequence and so 
     $\Kc^\bullet\sim \Oc_{X_T}$. \qed   
     
     \vskip .2cm
     
     \noindent{\sl End of proof of Proposition \ref {prop:tor-ind}:} 
     We factor $\varphi$ as the composition
     \[
     T\buildrel\gamma\over\lra S\times T \buildrel \pi_S\over\lra S, \quad 
     \gamma(t) = (\varphi(t), t), \,\,\pi_S(s,t) = s
     \]
     of a regular embedding and a smooth projection. So we decompose our square 
     as the concatenation of two other Cartesian squares
     \[
     \xymatrix{
     X_T\ar[d]  \ar[r] & X\times T\ar[d]  \ar[r]^{  \hskip 0.5cm \pi_X} & X\ar[d] 
     \\
     T\ar[r]_{\hskip -0.5cm \gamma} & S\times T \ar[r]_{\hskip 0.5cm \pi_S}& S
     }
     \]
     of which the right one is Tor-independent because $\pi_S$ is smooth an the left one
     is Tor-independent by the lemma. \qed

\begin{wrapfigure}[11]{r}{0.3\textwidth}
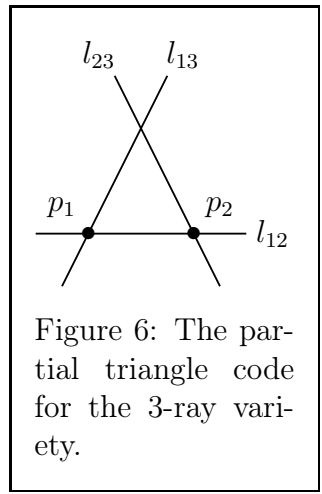

  
    \begin{framed}\raggedleft

  \btp[scale=.7]
     \centering
   \draw[line width=0.7] (-2,0) -- (2,0);
   \node at (-1,0) {$\bullet$}; 
  \node at (1,0) {$\bullet$};   
    \draw[line width=0.7] (-1.5, -1) -- (0.5, 3);
     \draw[line width=0.7] (1.5, -1) -- (-0.5, 3); 
  \node at (-1.5, 0.5) {$p_1$};  
   \node at (1.5, 0.5) {$p_2$};  
     \node at (2.5,0) {$l_{12}$};   
\node at (-0.8, 3.4) {$l_{23}$};     
\node at (0.8, 3.4) {$l_{13}$};    

  \etp
  
  \caption{The partial triangle code  for the $3$-ray variety.}\label{fig:3-ray}
    \end{framed}
    \end{wrapfigure}

    \paragraph{E. $3$-ray varieties.} 
           Let $A_0, \cdots, A_3$ be four consecutive chambers as in Fig. \ref{fig:4-cham}. 
    We then have two $2$-ray varieties
    \[
    Y_{01, 12} \,=\, Y_{A_{01}}\times_{Y_{A_1}} Y_{A_{12}}, \quad Y_{12, 23} = Y_{A_{12}}\times_{Y_{A_2}} Y_{A_{23}},.
    \]
     both isomorphic to $Y_{(2)}$. 
    The {\em $3$-ray variety} is defined as their fiber product over $Y_{A_{12}}$:
    \[
    Y_{(3)} \,=\, Y_{01, 12, 23}\ : =\, Y_{01, 12} \times_{Y_{A_{12}}} Y_{12, 23}\,=\, Y_{A_{01}}\times_{Y_{A_1}} Y_{A_{12}} 
    \times_{Y_{A_2}} Y_{A_{23}},
    \]
  so that we have the Cartesian square
 \be\label{eq:3-ray-sq}
  \xymatrix{
  Y_{(3)} \ar[r]^{\zeta_2} \ar[d]_{\zeta_1}  & Y_{01, 12}\ar[d]^{\sigma_1}
  \\ 
  Y_{12, 23} \ar[r]_{\sigma_2} & Y_{A_{12}}.
   }
  \ee 
  We encode $Y_{(3)}$ 
  by the partial triangle in   Fig. \ref{fig:3-ray}. Note the absence of the point $p_3$. 
    By the  {\em $3$-ray central fiber} we will mean
    the fiber product of the central fibers of the factors of $Y_{(3)}$,  which is the same as
  the 
 incidence correspondence associated to the partial triangle code: 
  \[
    \begin{gathered}
   F_{(3)} \,=\,  F_{01, 12, 23} \,\,=\,\ F_{01, 12} \times_{F_{A_{12}}} F_{12, 23}\,=\, F_{A_{01}}\times_{F_{A_1}} F_{A_{12}} 
    \times_{F_{A_2}} F_{A_{23}}\,\,= 
   \\
   = \,\,
   \bigl\{ (p_1, p_2, l_{12}, l_{13}, l_{23}) \in (\PP^2)^2 \times (\check\PP^2)^3 \,\bigl|
    \,\, p_1\in l_{12}, l_{12}, \,\,\, p_2\in l_{12}, l_{23} \bigr\}. 
   \end{gathered} 
    \]

    This is nothing but the Demazure resolution of the Schubert correspondence in $F\times F$ associated to
    the reduced decomposition $(23)(12)(23)$ in $W=S_3$. 
    It has been used  classically, see \cite{roberts},  as a tool to analyze  Schubert's 
     variety $\wh T$ of complete triangles.

     By construction,  we have a diagram
 \[
 \xymatrix{
 Y_{(3)}\ar[dr]_\rho \ar[r]^{\hskip -0.5cm i} & F_{(3)} \times\gen
 \ar[d]^{\pi}
 \\
 & F_{(3)}
  }
 \]
 where $i$ is a closed embedding and $\pi$ is the projection. Each fiber of $\rho$ is a Lie subalgebra in $\gen$.

 
  Inside $F_{(3)}$ we have  the  locus of {\em degenerate partial \nobreak triangles}  
    \[
    F \,\,= \,\,\bigl\{ (p_1, p_2, l_{12}, l_{13}, l_{23}) \in F_{(3)} \,\bigl| \,\, p_1=p_2, \,\, l_{12}=l_{12}=l_{23}\bigr\}
    \]
    identified with  the flag variety $SL_3/B$.

    \begin{prop}\label{prop:Y3}
    (a) The variety $ F_{(3)}$ is smooth, of dimension $6$. It can be represented as a $\PP^1\times\PP^1$-bundle
    over $F_{A_{12}}$. 
    
    \vskip .2cm
    
    (b) The projection $T\to F_{(3)}$ (and therefore the composite projection $\wh T \to F_{(3)}$)
    is birational. 
    
        \vskip .2cm
    
    (c) Let  $q = (p_1, p_2, l_{12}, l_{13}, l_{23}) \in F_{(3)}$. Then the Lie subalgebra $\rho^{-1}(q)\subset\gen$
    is $2$-dimensional, unless 
$p_1=p_2$ and $l_{13}=l_{23}$. In the latter case, 
    it is $3$-dimensional. 
    
        \vskip .2cm
    
    (d) The $3$-ray variety $Y_{(3)}$ has all components of dimension 8 (in fact, we will see that it is irreducible). 
    \end{prop}
    
    \noindent {\sl Proof:} (a) and (b) are obvious.  To see (c), we denote by $L_{01, 12}$  and $L_{12, 23}$ the
    rank $3$ bundles of Lie algebras over $F_{01, 12}$ and $F_{12, 23}$ respectively,  whose total spaces
    are the $2$-ray varieties $Y_{01,12}$ and $Y_{12, 23}$, see Proposition \ref{prop:Y2} (b). They are subbundles
    in the trivial bundles with fiber $\gen$, i.e., their fibers are Lie subalgebras in $\gen$. Then, by definition of 
    the fiber product, 
    \[
    \rho^{-1}(q) = (L_{01, 12})_{(p_1, p_2, l_{13}, l_{23})} \,\cap \, (L_{12, 23})_{(p_1, p_2, l_{12}, l_{23})} \,\,\subset 
    \,\,\, \gen. 
    \]
  
   If  $p_1=p_2$ and $l_{13}=l_{23}$, then both terms in the intersection are equal.
   
   In any other case the two terms in the intersection are  two distinct $3$-dimensional subspaces contained in the 
   $4$-dimensional space $(L_{A_{12}})_{(p_1, p_2, l_{12})}$, so the intersection is of dimension  $2$. 
   
   \vskip .2cm
   
   (d)  Note that the $2$-ray variety $Y_{01,12}$, being the blow-up of $Y_{A_{12}}$ in a smooth subvariety of codimension 2, is, at least locally over $Y_{A_{12}}$, a Cartier divisor in a ${\mathbb P}^1$ bundle over $Y_{A_{12}}$. Let $W$ be the pull-back of this  ${\mathbb P}^1$ bundle from $Y_{A_{12}}$ to $Y_{12, 23}$. Then the fiber product $Y_{(3)}$ is a Cartier divisor in it. Hence all components of $Y_{(3)}$ are of the same dimension $8$. \qed

   

  \begin{cor}\label{cor:3-ray-sq} 
  The Cartesian square \eqref{eq:3-ray-sq}
is Tor-independent.
  \end{cor}
  
  \noindent {\sl Proof:} Follows from Proposition \ref{prop:tor-ind}, since $Y_{01, 12}$, $Y_{12,23}$ and $Y_{12}$
  are smooth of dimension $8$, and $Y_{(3)}$ has all components of dimension  $8$ by Proposition \ref{prop:Y3}(d). \qed

 \paragraph{F. $3$-ray variety as a fiber product of blowups.}  
  We now recall that the maps $\sigma_1$ and $\sigma_2$ in  \eqref{eq:3-ray-sq}
 realize the $2$-ray varieties as blowups of the 1-ray variety $Y_{A_{12}}$. More precisely,
 denote for short  
 \[
 F_{(1)} = F_{A_{12}}\,\, =\,\, F\times_{\check \PP^2} F\,\, = \,\, \bigl\{ (p_1, p_2, l_{12})\in \PP^2\times\PP^2\times \check\PP^2\,\bigl| \,\, p_1\in l_{12}, \, p_2\in l_{12} \bigr\}
  \]
  the central fiber of $Y_{A_{12}}$. 
 We further denote by $L  =L_{A_{12}}$ the rank $4$ bundle of Lie algebras on $F_{(1)}$
 whose total space is $Y_{A_{12}}$, see Proposition \ref{prop:F_C}(b). 
 We write simply $Y_{A_{12}} = L$.

 Then, by Proposition \ref{prop:Y2} (c), 
 \[
 Y_{01, 12} \,\,=\,\, \Bl_{E_1}(L), \quad Y_{12, 23} = \Bl_{E_2}(L), 
 \]
 where $E_1, E_2\subset L$ are two rank $2$ subbundles of Lie algebras in $L$ defined as
 follows.  As before we identify elements $\xi$ of $\gen=\sel_3$ with global vector fields on $\PP^2$.
 Then the fibers of $E_i$, $i=1,2$, are defined by:
 \[
 (E_i)_{(p_1, p_2, l_{12})} \,\,=\,\, \bigl\{ \xi\in L_{(p_1, p_2, l_{12})}\subset\gen \,\bigl|
 \,\, \text{the linear part of $\xi$ at $p_i$ is scalar}
 \bigr\}. 
 \]
 
 \begin{prop}
 Let $q= (p_1, p_2, l_{12})\in F_{(1)}$. 
 
 \vskip .2cm
 
 (a) If $p_1\neq p_2$, then the fibers 
 $(E_1)_q, (E_2)_q
 \subset L_q$
 are transversal, i.e., intersect only at $0$. 
 
 \vskip .2cm
 
 (b) If $p_1 = p_2$, then
  the fibers $(E_1)_q$ and $(E_2)_q$ are equal. 
 \end{prop} 
 
 \noindent {\sl Proof:} (a) Suppose  $p_1\neq p_2$ and $\xi\in (E_1)_q\cap (E_2)_q$.
 We need to prove that $\xi=0$. 
 Consider the dual plane $\check \PP^2$, with the  lines $P_i$ corresponding to the points $p_i$
 and let $\eta$ be the vector field on $\check \PP^2$ corresponding to $\xi$. It suffices
 to prove that $\eta=0$.
  
 The fact that $\xi$ preserves $p_i$ and has scalar linear part at $p_i$,  means
 that $\eta$ preserves each point of $P_i$, i.e., vanishes on $P_i$ identically. 
 Now, for any global vector field $\zeta$ on  any $\PP^n$, any connected component of the
 zero locus of $\zeta$ is a projective subspace (since it corresponds to the eigenspace
 of the corresponding matrix). So the fact that $\eta$ vanishes on $P_1\cup P_2$ means
 that $\eta=0$. 
 
 Part (b) is obvious. \qed
 
 It follows from this proposition that the  intersection of the two blow-up centers $E_1\cap E_2$ has two irreducible components: one of codimension 4, the zero section of $L$,  and the other one of codimension 3, which is the locus of all points in the fibers in the part (b) of the above proposition.

   \begin{prop}\label{prop:Y3-rat}
    The variety $Y_{(3)}$ is irrducible and has rational singularities. That is, for any proper birational morphism
   $k:  M\to Y_{(3)}$ with $M$ smooth, we have $Rk_*(\Oc_M) \simeq \Oc_{Y_{(3)}}$. In fact, the singular locus $Q$ of $Y_{(3)}$ is a smooth 5-dimensional variety and the singularity in the transversal direction to it is a 3-dimensional quadratic cone. The variety $Q$ is identified via the map $Y_{(3)}\to Y_{A_{12}}$  with the blow up in the zero section of the total space of the rank 2 vector bundle $E_1|_F=E_2|_F$ over the flag variety $F\subset F_{(1)}$.

   \end{prop}

   \noindent{\sl Proof:} The map $\zeta_1 :Y_{(3)}\to Y_{12, 23}$ is an isomorphism outside the preimage of the intersection of the two blow-up centers $L_{A_{12}}\cap L_{A_{23}}$. The later has one component of codimension 3 and another one of codimension 4, as it was mentioned above. Then the preimage of the first component has an irreducible component of minimal possible codimension 2 in $Y_{12, 23}$. Since the dimension of fibers for $\zeta_1 :Y_{(3)}\to Y_{12, 23}$ has maximal possible dimension 1, the irreducibility of $Y_{(3)}$ follows if there is no component of $Y_{(3)}$ over this component of codimension 2. This follows from local calculations below.   

   The variety $F_{(1)}$, the bunlde $L$ on it and the subbundles $E_1$ and $E_2$
   are all acted upon by the group $SL_3$. The action on $F_{(1)}$ has two orbits:
   the diagonal $F \subset F\times_{\check\PP^2}F$, a codimension 1 submanifold,  and the complement to it. 
   
   So to unerstand the singularities of  the fiber product of the two blowups 
   it is enough to take a $1$ dimensional transversal slice $T$
   to $F$, consider the restriction $L'= L|_T$, its blowups along the $E_i' = E_i|_T$ and look at their fiber product.    
   
   The following local model describes this situation:
   \begin{itemize}
   
 \item   $T= \AAA^1_t := \Spec \, \k[t]$, the affine line with coordinate $t$ (we will use similar notation below). 
   
 \item  $L' = \AAA^5_{x_1, y_1, x_2, y_2, t}$ is the trivial rank $4$ bundle over $T$ with coordinates in the fiber being $x_1, y_1, x_2, y_2$. 
 
 \item $E'_1, E'_2$ are two rank $2$ subbundles which are transversal outside $t=0$ and coincide for $t=0$, 
 given explicitly by
 \[
 \begin{gathered}
 E'_1 \,\,=\,\,\bigl\{ (x_1, y_1, x_2, y_2, t) \,\bigl| \, x_2=0, \, y_2=0\bigr\}, 
 \\
  E'_2 \,\,=\,\,\bigl\{ (x_1, y_1, x_2, y_2, t) \,\bigl| \, x_2=tx_1, \, y_2=ty_1\bigr\}. 
 \end{gathered}
 \]
   \end{itemize}
   We now write the blowups $\Bl_{E'_i}(L')$ explicitly in local charts. We recall general formulas.
   
  Let $S$ be a smooth algebraic variety and $x,y\in \k[S]$ be two regular functions
  defining smooth divisors meeting transversally. Then $\Bl_{x=y=0}(S)$ is the union of  two
  charts, one given in $S\times \AAA^1_z$ by the equation $y=xz$, the other given in
  $S\times \AAA^1_{z'}$ by the equation $x=yz'$. 
  
  Accodingly, $\Bl_{E'_1}(L')$ has a local chart in $\AAA^6_{x_1, y_1, x_2, y_2, t, z_1}$ given
  by $y_2=x_2z_1$ (and  another chart with the roles of $x_2, y_2$ interchanged).
  Similarly, $\Bl_{E'_2}(L')$ has a local chart in  $\AAA^6_{x_1, y_1, x_2, y_2, t, z_2}$
  given by
  $(y_2-ty_1) = (x_2-tx_1) z_2$ (and another chart with the roles of $x$'s and $y$'s interchanged). 
  So  the fiber product $\Bl_{E'_1}(L')  \times_{L'} \Bl_{E'_2}(L')$ has $4$ charts, of which we write one,
  it is given by two equations 
  \[
  \begin{gathered}
  y_2 = x_2\cdot z_1,
  \\
  (y_2-ty_1) = (x_2-tx_1)\cdot z_2
  \end{gathered}
  \]
  in the affine space $\AAA^7_{x_1, x_2, y_1, y_2, t, z_1, z_2}$. We analyze this chart, the other charts
  being similar. 
  
  Eliminating $y_2$ from the first equation, we represent this chart as a hypersurface
  $  \Sigma \subset  \AAA^6_{ x_1, x_2, y_1,  t, z_1, z_2}$ given by the equation
  \[
  t(y_1-x_1z_2) \,\,=\,\,x_2(z_1-z_2).
  \]
  In other words, $\Sigma=h^{-1}(Q)$, where $Q\subset \AAA^4_{a,b,c,d}$ is the $3$-dimensional
  quadratic cone $\{ad=bc\}$ and $h: \AAA^6\to\AAA^4$ is the regular map given by
  \[
  a=t, \quad b=x_2, \quad c=z_1-z_2, \quad d=y_1-x_1z_2.
  \]
 We verify at once that $h$ is smooth in any point lying over $Q$. Therefore   
  $\Sigma_\sing$, the  singular locus of $\Sigma$,  is a smooth subvariety given (in $\AAA^6$) by the equations
  $a=b=c=d=0$ and near any point of $\Sigma_\sing$ the variety $\Sigma$ locally looks like the product
  of $\Sigma_\sing$ and $Q$. So it has rational singularities.

  By the the similar analysis of the  other charts of the fiber product, we get the structure of singularities of the blow-up. The locus of singularities is a smooth surface which is mapped birationally on the plane $(t=0,\ x_2=0,\ y_2=0)$ and is in fact the blow-up of this plane at the origin. The singularities of the fiber product in the transversal direction of every point of this surface are three dimensional quadratic cones.
  
  Since the plane $(t=0,\ x_2=0,\ y_2=0)$ corresponds exactly to the fiber of $E_1$ over a point in the flag variety $F\subset F_{(1)}$, the description of the variety $Q$ follows.

   \vfill\eject

   \section {Properties of the $\sel_3$-flober}\label{sec:sl3-flob}
   
   \paragraph{A.  Definition and main result. } We recall the flop diagram  $(Y_C, l_{CC'})_{C\leq C'}$  for $\sel_3$, see Definition \ref{def:sl3-flop-d}.

      \begin{Defi}\label{def:sl3-flober}
   The {\em $\sel_3$-flober} is the diagram consisting  of the coherent derived categories  $\Ec_C:=D^b_\coh(Y_C)$
       of the varieties from the $\sel_3$-flop diagram (Definition \ref{def:sl3-flop-d}) 
  and the  functors
   \[
   \gamma_{CC'} = R (l_{CC'})_*: \Ec_C\lra \Ec_{C'}, \quad
   \delta_{C'C} = L(l_{CC'})^*: \Ec_{C'}\lra\Ec_C, \quad   C\leq C'.
   \]
    We denote this diagram by $\Fen = \Fen_{\sel_3}$. 
  \end{Defi}
  The rough size of the  categories $\Ec_C$ can be understood from the following. 
  
  \begin{prop}
   Let   $E_C = K(\Ec_C)$ be
    the rational Grothendieck group of $\Ec_C$ (or, what is the same, of $Y_C$).
  Then:
  \[
  \dim_\QQ (E_C) = \begin{cases}
  6 & \text{ if } \dim(C)=2;
  \\
   12& \text{ if } \dim(C)=1;
   \\
  72 & \text{ if } \dim(C)=0.  
  \end{cases}
  \]
  \end{prop}
 
 \noindent{\sl Proof:} As each $Y_C$ is the total space of a vector bundle over its central fiber $F_C$,
 we have $K(Y_C)\simeq K(F_C)$. Now, for $\dim(C)=2$ we have $F_C = F$ is the flag variety for $\sel_3$,
 so its K-group has rank $|W|=6$. For $\dim(C)=1$ we have $F_C = F\times_{\check \PP^2} F$ is
 the total space of a $\PP^1$-bundle over $F$, so its K-group has rank $12$. Finally, for $C=0$
 we have $F_C=\wh T$ is the Schubert space of complete triangles. From its explicit representation
 as an iterated blowup of $(\PP^2)^3$ recalled in \cite{roberts}, Ex. 7.7.3, 
  it follows that $K_0(\wh T)$ is free and has the same rank as
 the total Chow group $A^\bullet(\wh T)$.  The structure of the latter group is also known, see
 \cite{roberts} Th. 3.2 and references therein, as well as \cite{collino-fulton}.  It is free of  rank 72. \qed
 
 \vskip .2cm

  The rest of the paper is devoted to the proof of the following.

     \begin{thm}\label{thm:main}
   The diagram of categories $(\Ec_C, \gamma_{CC'}, \delta_{C'C})$ satisfies all the conditions
    of  Definition \ref{def:H-flober}
   of $\Hc$-schobers except, possibly,  Condition (4) (invertibility) for the case $C'=-C$, $\dim(C)=1$, $D=0$.  
 \end{thm}
 
 Studying Condition (4) in te remaining case,  i.e., investigating the effect of the Schubert transform 
  $
 Y_{-C}\buildrel l_{0, -C}\over  \longleftarrow Y_0 \buildrel l_{0,C}\over \lra Y_C,
 $
  (see Definition \ref{def:sl3-flop-d}(b)) on derived categories, is an interesting question which we plan
  to address in a future work.

 \paragraph{B. Idempotency and invertibility.} 
  We now start the proof of Theorem \ref{thm:main}. The fact that the $\gamma_{CC'}$ form a 2-functor from
 $\Cc$ (the poset of faces)
  to triangulated categories, is a consequence of the fact that the $l_{CC'}$ form
 a functor from $\Cc$ to algebraic varieties. 
   Let us prove Condition (1) (idempotency). 
 Each $l_{CC'}: Y_C\to Y_{C'}$ is a regular birational, proper morphism between smooth varieties.
  Therefore  we have $R(l_{CC'})_*(\Oc_{Y_C}) = \Oc_{Y_{C'}}$.  So by the projection formula,
  for each $\Fc\in  \Ec_{C'}= D^b_\coh(Y_{C'})$ we have a natural identification
  \[
  \gamma_{CC'} \delta_{C'C}(\Fc) \,\,\simeq \,\, R(l_{CC'})_* (l_{CC'})^*\Fc \,\,\simeq \,\, 
  R(l_{CC'})_* (\Oc_{Y_C})\otimes\Fc \,\,=\,\,\Fc. 
  \]

  Next, let us look at Condition (4)  (invertibility) of Definition \ref {def:H-flober}). 
  Apart from Conjecture \ref{conj:schub-flop}, the only case we need
  to consider is $\dim(C)=\dim(C')=2$, $\dim(D)=1$. In this case the diagram
  $Y_C\leftarrow Y_D \to Y_{C'}$ is, by Corollary \ref{cor:x-wall-sm}
 locally a product of a smooth manifold and of the diagram \eqref{eq:a-fib-prod} for the Atiyah flop
 (case $\gen=\sel_2$).  
   So the fact that the corresponding flop functor
   is an equivalence, follows from the case $\gen=\sel_2$.

  \paragraph{C. Collinear transitivity.} 
  
  We now turn to  verifying Condition (4) (collinear transitivity) of  Definition \ref{def:H-flober}. 
  That is, we define the flopping functor $\varphi_{CC'}: \Ec_C\to\Ec_{C'}$ as in \eqref{eq:flop-functors} and  verify that
  $\varphi_{BC}\circ\varphi_{AB} = \varphi_{AC}$ whenever $A,B,C$ are collinear faces of $\Hc$. 
  For this we distinguish several cases as to the possible dimensions of $A,B,C$. We will call Case $(i,j,k)$
  the situation when $\dim(A)=i$, $\dim(B)=j$, $\dim(C)=k$.  We only need to consider nonzero $i,j,k$. 
 Each such case may have several subcases
  as to the relative positions of $A,B,C$ with given dimensions.   
 
   \paragraph {  D. Case (2,2,2).}  Our arguments in this case are  parallel to those of \cite{riche}
   which establish the braid group relations between the  flop 
   functors acting on the $\Ec_C$ for $C$ running over chambers. 
   In fact, one has a description
   of the fundamental groupoid of the open stratum $\hen^*_\CC - \Hc_\CC$ in terms
   of morphisms $\varphi_{CC'}$ for  $\dim(C)=\dim(C')=2$ satisfying collinear transitivity, see  \cite{KS}, Appendix,
   so Case (2,2,2) also  proves these braid relations.    
   
   \vskip .2cm

    We consider several possibilities. First, we consider the {\em neighboring case}. That is, $A,B,C$
    are positioned as the faces $A_1, A_2, A_3$ in Fig. \ref{fig:4-cham} (in one of the two possible orientations). 
     By definition, $\varphi_{AB}$ and $\varphi_{BC}$ are
  defined by pullback and pushforward in the diagrams
  \[
  Y_A \lla Y_\alpha \lra Y_B, \quad Y_B \lla Y_\beta\lra Y_C,
  \]
   where $\alpha$ and $\beta$ are the rays between $A$ and $B$ and between $B$ and $C$, respectively.
   By Proposition \ref{prop:Y2}(d), 
    the composition $\varphi_{BC}\varphi_{AB}$ is defined by  pullback and pushforward in the diagram
  \[
  Y_A\buildrel p_A\over \lla Y_\alpha \times_{Y_B} Y_\beta \buildrel p_C\over \lra Y_C
  \]
 (the middle term $Y_\alpha \times_{Y_B} Y_\beta$ is the $2$-ray variety $Y_{(2)}$
 studied in \S \ref{sec:schubert}D and the square defining the fiber product is Tor-independent).

 On the other hand, the functor $\varphi_{AC}$ is  
 defined by  pullback and pushforward in the diagram
\[
Y_A \buildrel q_A\over \lla Y_0 \buildrel q_C\over \lra Y_C. 
\]  
 The identification of   $\varphi_{BC}\varphi_{AB}$ and  $\varphi_{AC}$ is implied by 
  Proposition \ref{prop:Y2} (which says that 
    $Y_{(2)}= Y_\alpha \times_{Y_B} Y_\beta$ is smooth)  which we can combine with the  commutative diagram
 \[
 \xymatrix{
  & Y_0  \ar[dl]_{q_A}  \ar[d]^{\sigma} \ar[dr]^{q_C}&
  \\ 
  Y_A&  \ar[l]_{p_A}  Y_\alpha\times_{Y_B}Y_\beta \ar[r]^{p_C}& Y_C
 }
 \]
 Indeed, smoothness implies that $R\sigma_*(\Oc)=\Oc$ and so $R\sigma_*\circ\sigma^* = \Id$. Therefore
 \[
 \varphi_{BC} \varphi_{AB} \,=\, R p_{C*} p_C^* \,=\,  R p_{C^*} R\sigma_* \sigma^* p_A^* \,= \,Rq_{C*} q_C^*
 \, = \, \varphi_{AC}.
 \]
 
 Next, we consider the {\em non-neighboring cases}. That is, $A,B,C$ are positioned as
 $A_0, A_2, A_3$ or as $A_0, A_1, A_3$ in  Fig. \ref{fig:4-cham} (in one of the two possible orientations). 
 For definiteness, let us consider the first of these possibilities. Then by the neighboring case above,
 $\varphi_{A_0, A_2}$ is given by pullback and pushforward through the $2$-ray variety, so
 $\varphi_{A_2, A_3} \varphi_{A_0, A_2}$ is, by Corollary \ref {cor:3-ray-sq}  (Tor-independence of the
Cartesian  square) given by pullback and pushforward through the $3$-ray variety $Y_{(3)}$.
Our statement then follows, similarly to the neighboring case above,  from 
 Proposition \ref{prop:Y3-rat} applied to the morphism $Y_0\to Y_{(3)}$. 
 
 \paragraph{ E. Case $(2,1,2)$ (neighboring).} Assume that $A,B,C$ are positioned as 
 $A_1, A_{12}, A_2$ in Fig.  \ref{fig:4-cham} (in one of the two possible orientations). 
 In this case the transitivity follows formally from idempotency. Indeed, by definition,
 $\varphi_{AC}=\varphi_{0C}\varphi_{A0}$, while substituting $\Id=\varphi_{0B}\varphi_{B0}$
 and re-arranging the brackets, we get 
 \be\label{eq:2-ray-diag1}
 \varphi_{BC}\varphi_{AB} \,\,=\,\, \varphi_{BC}\bigl(  \varphi_{0B}\varphi_{B0}\bigr) \varphi_{AB} 
 \,\, = \,\, 
\bigl( \varphi_{BC}  \varphi_{0B} \bigr) \bigl( \varphi_{B0} \varphi_{AB}\bigr) 
\,\,=\,\, 
 \varphi_{0C}\varphi_{A0}.  
 \ee
 
 \paragraph{F. Case $(1,2,1)$ (neighboring).}  Assume that $A,B,C$ are positioned as 
 $A_{01}, A_1, A_{12}$ in  Fig.  \ref{fig:4-cham} (in one of the two possible orientations). 
 Thus $Y_B$ is isomorphic to the Grothendieck resolution. We can assume that $Y_B=\wt\gen$
 is the standard Grothendieck resolution, i.e., $B$ is the dominant Weyl chamber. 
 
 Now, the functor $\varphi_{AC} = \varphi_{0C}\varphi_{A0}$ is defined by the pullback and push forward 
 through $Y_0$. This is the same as   pullback and pushforward
 through the $2$-ray variety $Y_A\times_{Y_B} Y_C$, as follows from a diagram similar to
 \eqref{eq:2-ray-diag1}: 
  \[
 \xymatrix{
  & Y_0  \ar[dl]  \ar[d]^{\sigma} \ar[dr]&
  \\
  Y_A&  \ar[l]^{\hskip -0.5cm s}  Y_A \times_{Y_B}Y_C \ar[r]_{\hskip 0.5cm p} & Y_C. 
 }
 \]
 On the other hand, the functor $\varphi_{BC}\varphi_{AB}$ is defined by pushforward and pullback
 through $Y_B=\wt\gen$, i.e., as $r^* Rq_*$ in the square
 \[
 \xymatrix{
 Y_A\times_{Y_B} Y_C\ar[d]_s  \ar[r]^{\hskip 0.5cm p} & Y_C\ar[d]^r
 \\
 Y_A\ar[r]_q &Y_B
 }
 \]
 So if we show that the square is Tor-independent, we will  identify $r^* Rq_*$ with $Rp_* s^* = \varphi_{AC}$. 
 By Corollary \ref{cor:x-wall-sm}, the varieties $Y_A$ and $Y_C$
 are the blowups of $Y_B$ along two codimension $2$ subvarieties $\Lambda_{AB}, \Lambda_{CB}$
 which are {\em transversal}. More precisely, $Y_B=\wt\gen$ is the total space of the tautological bundle 
 $\ul b\to F$ of Borels over $F$, and $\Lambda_{AB}, \Lambda_{CB}$ are total spaces of two
 transversal $SL_3$-equivariant subbundles whose fibers over the  standard base point of $F=G/B$
 are two subalgebras in the standard Borel $\ben$ of the form $\ben\cap s_{i}(\ben)$, where $s_i$, 
 $i=1,2$, are two simple reflections in the Weyl group. 
 
 Now, it is standard (and easily follows from Proposition
 \ref{prop:tor-ind}) that the fiber square of two transversal blowups is Tor-independent.

 \paragraph{G. Case $(1,1,1)$ (neighboring).} Assume that $A,B,C$ are positioned as
 $A_{01}, A_{12}, A_{23}$ in  Fig.  \ref{fig:4-cham} (in one of the two possible orientations). 
 The identity $\varphi_{AC}=\varphi_{BC}\varphi_{AB}$ in this case follows from the Tor-independence
 of the square defining the $3$-ray variety $Y_{(3)}$  (Corollary \ref{cor:3-ray-sq}) and from the fact
 that $Y_{(3)}$ has rational singularities (Proposition \ref{prop:Y3-rat}).

 \paragraph {H.  Case $A\geq B$ or $B\leq C$.}  If $A\leq B$, then for any $C$ we have 
 \[
 \varphi_{BC}\varphi_{AB} \,\, = \,\, \varphi_{0C}\varphi_{B0}\varphi_{AB} \,\, \buildrel A\geq B\geq 0 \over =
 \,\, \varphi_{0C}\varphi_{A0}\,\,=\,\,\varphi_{AC}. 
 \]
 We have a similar statement and argument when $B\geq C$. 
 
 \paragraph{ I. Reduction to neighboring cases.}  We now give an inductive argument which reduces
 all the remaining cases to the neighboring cases that have been already considered. 
 
 For any two faces $A,B$ of $\Hc$ of dimensions $\neq 0$ we denote by $d(A,B)$ the
 {\em incidence distance} between $A$ and $B$, i.e., the  minimal $m$ such that there exists a chain
 of faces of dimensions $\neq 0$
 \[
 D_0=A, \,D_1,\, \cdots,\,  D_m=B
 \]
 such that for any $i$ we have $D_i<D_{i+1}$ or $D_i>D_{i+1}$. For example, if $A$ and $B$ are two
 adjacent chambers, then $d(A,B)=2$. 
 
 \vskip .2cm
 
 Suppose we know transitivity $\varphi_{AC}=\varphi_{BC}\varphi_{AB}$ for all collinear $A,B,C$ with $d(A,C)<L$.
 Let us prove it in the case $d(A,C)=L$. 
 
 \vskip .2cm
 
 Assume first that $\dim(A)=2$. In this case we take the $1$-dimensional face $D<A$ in the direction of $B,C$. 
 We can then write by inductive assumption
 \[
 \begin{gathered}
 \varphi_{BC}\varphi_{AB} \,\, \buildrel \on{ind.} \over =\,\, (\varphi_{BC} \varphi_{DB})\varphi_{AD}
 \,\, = \,\, \varphi_{BC} (\varphi_{DB} \varphi_{AD})  \,\, \buildrel \on{ind.} \over =\,\, 
 \\
  \,\, \buildrel \on{ind.} \over =\,\, \varphi_{DC}\varphi_{AD}  \buildrel \on{(H)} \over =\,\, \varphi_{AC},
 \end{gathered}
 \]
where the last equality is an instance of Case H above. 

A similar argument works when $\dim(C)=2$. 

\vskip .2cm 
  
  So our statement reduces to the case $\dim(A)=\dim(C)=1$ which we now assume. Now, if $\dim(B)=2$,
  then we take a $1$-dimensional face $D\leq B$ in the direction of $A$ or $C$ (there must be a gap in
  one of the directions, otherwise we are in a neighboring case).   Suppose $D$ is in the direction of $A$.
  Then, by geometry of our arrangement ($3$ lines only), the faces $D,B,C$ must be neighboring, i.e.,  fit into  Case F  and 
  the faces $A,D,C$ fit into Case G above, so
   we write
  \[
  \begin{gathered}
  \varphi_{BC}\varphi_{AB} \,\,\buildrel \on{ind.}\over =\,\, \varphi_{BC}( \varphi_{DB}\varphi_{AD}) \,\,=\,\,
   (\varphi_{BC} \varphi_{DB})\varphi_{AD}  \,\,\buildrel (F) \over =\,\,
   \\
     \,\,\buildrel (F)\over =\,\,
   \varphi_{DC}\varphi_{AD} \,\,\buildrel (G)\over = \,\, \varphi_{AC}. 
   \end{gathered}
  \]
This finishes the proof of Theorem \ref{thm:main}.

  \vfill\eject

  \let\thefootnote\relax\footnote {A.B.: 
  Kavli Institute for Physics and Mathematics of the Universe (WPI), 5-1-5 Kashiwanoha, Kashiwa-shi, Chiba, 277-8583, Japan, Steklov Institute of Mathematics, Moscow, Russia, National Research University Higher School of Economics, Russian Federation
Email: {\tt bondal@mi.ran.ru}
\\
M.K.: 
  Kavli Institute for Physics and Mathematics of the Universe (WPI), 5-1-5 Kashiwanoha, Kashiwa-shi, Chiba, 277-8583, Japan.
Email: {\tt mikhail.kapranov@ipmu.jp}
\\
V.S.:   
Universit\'e Paul Sabatier, 
Institut de Mathématiques de Toulouse, 
118 Route de Narbonne, 
31062 Toulouse, France.   
Email: {\tt schechtman@math.ups-tlse.fr }
}


\begin{thebibliography}{100}


 
  
 \bibitem{A1}  R. Anno Spherical functors, arXiv:0711.4409. 
 
  \bibitem{AL2}     R. Anno, T. Logvinenko.  Spherical DG-functors, 
arXiv:1309.5035.

\bibitem{babson1} E. Babson, P.  E. Gunnells, R. Scott. A smooth space of tetrahedra.  
 {\em Adv. Math.}  {\bf 165} (2002) 285–312.

\bibitem{babson2} E. Babson, P.  E. Gunnells, R. Scott. Geometry of the tetrahedron space.
  {\em Adv. Math. }  {\bf 204} (2006) 176–203.

\bibitem {beil-gluing} A.~Beilinson.
\newblock How to glue perverse sheaves. In: $K$-theory, arithmetic and geometry (Moscow, 1984), 
\newblock {\em  Lecture Notes in Math.} {\bf 1289}, 
Springer-Verlag, 1987, 42 - 51.

\bibitem {BBD}  A.~Beilinson, I.~Bernstein, P.~Deligne.
\newblock Faisceaux Pervers, {\em Ast\'erisque} {\bf 100}, 1982. 

 
 \bibitem{bez} R.  Bezrukavnikov. Noncommutative counterparts of the 
 Springer resolution. arXiv
math/0604445.

\bibitem{bez-talk}  R.  Bezrukavnikov. Commutative and noncommutative
symplectic resolutions and perverse sheaves. Lecture May 18 2015. 

 
 \bibitem{BR}
 R.~Bezrukavnikov, S.~ Riche.
  \newblock Affine braid group actions on derived categories of Springer resolutions.
  \newblock {\em  Ann. Sci. \'Ec. Norm. Sup\'er.} {\bf 45}  (2012), 535-599. 
  
  \bibitem{BB} A. Bodzenta, A. Bondal.	 Flops and spherical functors.  arXiv:1511.00665.
  
  \bibitem{BB18} A. Bodzenta, A. Bondal. Canonical tilting relative generators, {\em Advances in Mathematics}, {\bf 323} (2018) 226-278. 
  
  \bibitem{bondal-kapranov} A. I. Bondal, M. M. Kapranov. Representable functors,
  Serre functors and mutations.  {\em Math. USSR Izv.} {\bf 35}  (1990) 519–541.  

    \bibitem{BK-enhanced} A. I. Bondal, M. M. Kapranov. Enhanced triangulated categories.
   {\em Math. USSR Sbornik},  {\bf 70} (1991) 93-107. 
   
   \bibitem{BO} A. Bondal, D. Orlov. Semiorthogonal decomposition for algebraic varieties. 
   arXiv  alg-geom/9506012. 
   
   \bibitem{BO2} A. I. Bondal, D. O. Orlov. Derived categories of coherent sheaves. {\em Proceedings of the
International Congress of Mathematicians, Vol. II (Beijing, 2002)}, (2002) 47-56, Beijing, Higher Ed. Press.
   
   \bibitem{BvdB} A. I.  Bondal, M. van den Bergh. Generators and representability of functors in commutative and
    noncommutative geometry. {\em Mosc. Math. J. } {\bf 3}  (2003)  1-36. 
   
   \bibitem{bridg-flop} T. Bridgeland. Flops and derived categories. 
   {\em Invent. Math.} {\bf 146} (2002) 613-632. 
   
   \bibitem{brieskorn-german} E. Brieskorn. Die Aufl\"osung der rationalen Singulari\"aten holomorpher 
   Abbildungen, {\em Math. Ann.}  {\bf 178} (1968), 255–270.
   
   \bibitem{brieskorn} E. Brieskorn. Singular elements of semi-simple algebraic groups.
   {\em Actes Congr\'es Intern. Math. Nice}, 1970, Tome 2, 279-284. 
   
   \bibitem{collino-fulton} A. Collino, W. Fulton. Intersection rings of spaces of
   triangles. {\em M\'em. Soc. Math. France}, {\bf 38} (1989) 75-117. 
   
   \bibitem{chriss-ginz} N. Chriss, V. Ginzburg. Representation Theory and Complex Geometry. 
   Birhka\"user, Boston, 2010. 
   
   \bibitem{donovan} W. Donovan. Perverse schobers and wall crossing.  arXiv:1703.00592. 
   
   
   \bibitem{DW} W. Donovan, M. Wemyss. 
  Twists and braids for general 3-fold flops.   	arXiv:1504.05320.
  
  \bibitem{Dri}
  V. Drinfeld.  DG quotients of DG categories. {\em  J.  Algebra}  {\bf 272}  (2004) 643-691. 
  
  \bibitem{DHKS} W. G. Dwyer, P. S. Hirschhorn, D. M. Kan, J. H. Smith.
  Homotopy Limit Functors on Model Categories and Homotopical Categories. AMS Publ. 2004. 
  
  \bibitem{DK-triang} T. Dyckerhoff, M. Kapranov. Triangulated surfaces in triangulated categories. 
  arXiv:1306.2545. 
  
  \bibitem{DKSS}  T. Dyckerhoff, M. Kapranov, V. Schechtman, Y. Soibelman.
  Fukaya categoris with coefficients, in preparation.  
  
  
     \bibitem{GGM} A.~Galligo, M.~Granger, P.~ Maisonobe.
\newblock $\Dc$-modules et faisceaux pervers dont le support singulier
est un croisement normal. \newblock {\em Ann. Inst. Fourier Grenoble},
{\bf 35} (1985), 1-48. 

\bibitem{harder} A. Harder, L. Katzarkov. Perverse sheaves of categories and some applications.
arXiv 1708.01181. 

\bibitem{hirschhorn} P. S. Hirschhorn. Model Categories and Their Localizations.
AMS Publ. 2003. 

\bibitem{hotta-kash} R. Hotta, M. Kashiwara. The invariant holonomic system on a semisimple Lie algebra.
{\em Invent. Math.} {\bf 75} (1984) 327-358. 

\bibitem{iliev1} A. Iliev, L. Manivel. Severi varieties and their varieties of reductions.
  {\em J. Reine Angew. Math.}  {\bf 585}  (2005), 93–139. 

\bibitem{iliev2} A. Iliev, L. Manivel. Varieties of reductions for $\gen\len_n$.
arXiv:math/0501329. 

    \bibitem{magyar}  W. van der Kallen, P.  Magyar.
 The space of triangles, vanishing theorems, and combinatorics.  
 {\em  J. Algebra}, {\bf  222} (1999) 17-50.

 
 \bibitem{KS} M. Kapranov, V. Schechtman, Perverse sheaves over real hyperplane 
arrangements,  {\em Ann. Math.} {\bf 183} (2016) 617-679. 

\bibitem {KS-schobers} M. Kapranov, V. Schechtman, Perverse Schobers, arXiv:1411.2772. 


   
 \bibitem{kaw} Y. Kawamata. On the cone of divisors of Calabi-Yau fiber spaces. 
 {\em  Int. J.  Math.}    {\bf 8}  (1997), 665-687. 
 
 \bibitem{kostant} B. Kostant. Lie group representations on polynomial rings. 
{\em Bull. Amer. Math. Soc.} {\bf 69} (1963) 518-526. 

\bibitem{Kuz1} A. Kuznetsov. Hyperplane sections and derived categories
 {\em Russ.  Math. Izv.}  {\bf 70} (2006), 447–547.  

\bibitem{Kuz2} A. Kuznetsov. Homological projective duality.
 {\em Publ. Math. Inst. Hautes \'Etudes Sci.}  {\bf 105} (2007), 157–220.
 
 
 \bibitem{namikawa} Y. Namikawa. Mukai flops and derived categories. II. In:
  ``Algebraic structures and moduli spaces",  {\em CRM Proc. Lecture Notes} 
  {\bf 38}, p.149–175, AMS Publ. 2004.
 
 \bibitem{pinkham}   H. Pinkham. 
Factorization of birational maps in dimension 3,
in: {\em Singularities} (P. Orlik, ed.), {\em Proc.
Symp. Pure Math.}  {\bf 40} pt.  2,  p. 343-371, American 
Math. Soc.,  Providence, 1983.

\bibitem{reid} M. Reid.  Young person’s guide to canonical singularities, in:  ``Algebraic geometry'' (Bowdoin, 1985), {\em Proc. Sympos. Pure Math.}  {\bf 46} Part 1, AMS 1987, p. 345–414.

\bibitem{riche} S. Riche. Geometric braid group actions on derived categories
of coherent sheaves (with a joint appendix with R. Bezrukavnikov).
{\em Representation Theory} {\bf 12} (2008) 131-169. 
 
 \bibitem{roberts} J. Roberts. Old and new results about the triangle varieties. 
 {\em Lecture Notes in Math.} {\bf 1311}, p. 197-219, Springer, Berlin, 1988. 
 
 \bibitem{schubert} H. Schubert. Anzahlgeometrische Behandlung des Dreiecks.
 {\em Math. Ann.} {\bf 17} (1880) 1213-1255. 
 
   \bibitem{semple} J. G. Semple. The triangle as a geometric variable.
 {\em Mathematika} {\bf 1} (1954) 80-88. 
 
 \bibitem{tabuada} G. Tabuada. Th\'eorie homotopique des DG-categories.
 arXiv:0710.4303. 
 
 \bibitem{tabuada2} G. Tabuada. Homotopy theory of dg categories via localizing pairs and Drinfeld's dg quotient.
 {\em Homology, Homotopy and Appl.} {\bf 12} (2010) 187-219. 
 
 \bibitem{toen} B. To\"en.
 {The homotopy theory of dg-categories and derived Morita theory. }
 {\em Invent. Math.}  {\bf 167}  (2007)  615–667. 

 \bibitem{vdB} M. van den Bergh.  Three-dimensional flops and noncommutative rings.
 {\em Duke Math. J.}  {\bf 122} (2004) 423–455.


   

  \end{thebibliography}
\end{document}